\documentclass[12pt]{article}
\pdfoutput=1
\usepackage{latexsym, graphicx,cite}
\usepackage{amsmath}
\usepackage{amssymb}
\usepackage{amsthm}
\usepackage{mathtools}
\usepackage{graphicx}
\usepackage{caption}
\usepackage{subcaption}
\usepackage{titlesec} 
\usepackage{titletoc} 
\usepackage{tocloft} 

\usepackage{color, colortbl}
\definecolor{LightCyan}{rgb}{0.88,1,1}
\usepackage{protosem}
\newcommand{\reye}{r_{\hspace{-0.08cm}\textproto{\Aayin}}}
\usepackage[top = 1. in,bottom = 1. in, left = 1 in, right=1 in]{geometry}

\newcommand{\arXiv}[1]{\href{http://www.arXiv.org/abs/#1}{arXiv:#1}}
\usepackage[colorlinks=true, linkcolor=blue, bookmarks=true]{hyperref}

\titleformat{\section}[hang]{\normalfont\large\bfseries}{\thesection.}{1pc}{}
\titlespacing{\section}{0pc}{-0.ex plus .1ex minus .2ex}{-1.pc}
\titleformat{\subsection}[hang]{\normalfont\normalsize\bfseries}{\thesubsection.}{1pc}{}
\titlespacing{\subsection}{0pc}{-0.5ex plus .1ex minus .2ex}{-1pc}
\titleformat{\subsubsection}[hang]{\normalfont\normalsize\bfseries}{\thesubsubsection.}{1pc}{}
\titlespacing{\subsubsection}{0pc}{-0.5ex plus .1ex minus .2ex}{-1pc}

\titlecontents{section}[1.72em]{\addvspace{-0.5pc}\bfseries}{\thecontentslabel.\enspace}{}{\titlerule*[0.6pc]{.}\contentspage}

\newcommand{\beq}{\begin{equation}}
\newcommand{\eeq}{\end{equation}}
\newcommand{\beqnn}{\begin{equation*}}
\newcommand{\eeqnn}{\end{equation*}}
\newcommand{\ber}{\begin{array}}
	\newcommand{\eer}{\end{array}}

\newcommand{\LL}{\Omega}
\newcommand{\s}{\sigma}

\newcommand{\w}{\omega}

\newcommand{\ena}{\end{eqnarray}}
\newcommand{\beqa}{\begin{eqnarray}}
\newcommand{\eeqa}{\end{eqnarray}}
\newcommand{\bea}{\begin{eqnarray}}
\newcommand{\eea}{\end{eqnarray}}

\theoremstyle{proposition}

\theoremstyle{remark}

\begin{document}
	\begin{titlepage}
		\begin{flushright}
			\phantom{arXiv:yymm.nnnn}
		\end{flushright}
		\vspace{1cm}
		\begin{center}
			\begin{changemargin}{-0.6cm}{-0.6cm} 
				\begin{center}
					{\Large\bf
					Self-Similar Solutions to the Compressible Euler Equations and their Instabilities
				}\\
				\end{center}
			\end{changemargin}
			\vspace{0.5cm}
			{\large Anxo Biasi}
			\vskip -2mm
			{\em  Institute of Theoretical Physics, Jagiellonian University, Krakow, Poland}
			\vskip -2mm
			{\small\noindent {\tt anxo.biasi@gmail.com}}
			\vskip 10mm
		\end{center}
		\vspace{0cm}
		\begin{center}
			{\large\bf Abstract}\vspace{-0.2cm}
		\end{center}	
		
			\noindent
			This paper addresses the construction and the stability of self-similar solutions to the isentropic compressible Euler equations. These solutions model a gas that implodes isotropically, ending in a singularity formation in finite time. The existence of smooth solutions that vanish at infinity and do not have vacuum regions was recently proved and, in this paper, we provide the first construction of such smooth profiles, the first characterization of their spectrum of radial perturbations as well as some endpoints of unstable directions. Numerical simulations of the Euler equations provide evidence that one of these endpoints is a shock formation that happens before the singularity at the origin, showing that the implosion process is unstable.

			
		\hypersetup{linkcolor=black}
		\tableofcontents

	\end{titlepage}
	

	\section{Introduction\vspace{3mm}}\label{sec:Intro}
	
	Singularity formation in partial differential equations (PDEs) is a fundamental problem of great interest in physics and mathematics. A particular materialization of these phenomena are self-similar solutions. They commonly follow the scaling law of the PDE that they solve; however, there are exceptions to this rule with solutions following anomalous scaling. This paper is concerned with the latter self-similar solutions in the scenario of the {\em isentropic compressible Euler equations}  
	\beq
		\begin{cases}
			\partial_t\rho + \nabla\cdot(\rho u) = 0\\
			\rho\partial_tu + \rho u \cdot \nabla u + \nabla p = 0\\
			p = \frac{\gamma-1}{\gamma}\rho^\gamma\\
			\rho(t,y)>0
		\end{cases}
		\label{eq:Euler_eq}
	\eeq
	with $y\in\mathbb{R}^d$, $d\geq2$ and $\gamma>1$. The equation of state $p(\rho)$ is associated with an ideal gas where $\gamma$ is the {\em heat capacity ratio}.
	
	This setup is a particular case where, in 1942, Guderley found spherically symmetric self-similar solutions \cite{Guderley} that model an imploding gas that collapses to the center of symmetry (see \cite{Meyer} for an excellent review). In the context of fluid dynamics this problem belongs to the family of classical self-similar solutions: Guderley problem (converging/diverging shock waves) \cite{Guderley,Landau,Stanyukovich}, Noh problem  (blast-waves) \cite{Noh}, Sedov-Taylor-von Neumann problem (blast-waves) \cite{Taylor,Sedov}, Larson-Penston problem (self-gravitating fluid) \cite{Larson,Penston}. Despite these solutions were discovered decades ago, they have become relevant to current research. Some examples are the extension of these solutions \cite{BoydNN, Velikovich1}, code verification \cite{GuderleyReview, CodeVerification}, rigorous construction of smooth solutions \cite{MerleEuler, LarsonPenston2020}, etc. Other recent results in singularity formation in the Euler equations are \cite{Eggers,Christodoulou1,Christodoulou2,Christodoulou3,Vlad1,Vlad2,Vlad3,Alexis,Buckmaster,Olga,Laurent}. See also \cite{Rod,Zoe,Oli} and references therein for stabilization of relativistic fluids.
	
\begin{figure}[h!]
			\vspace{-0.8cm}
		\centering	
		\includegraphics[width=8.3cm]{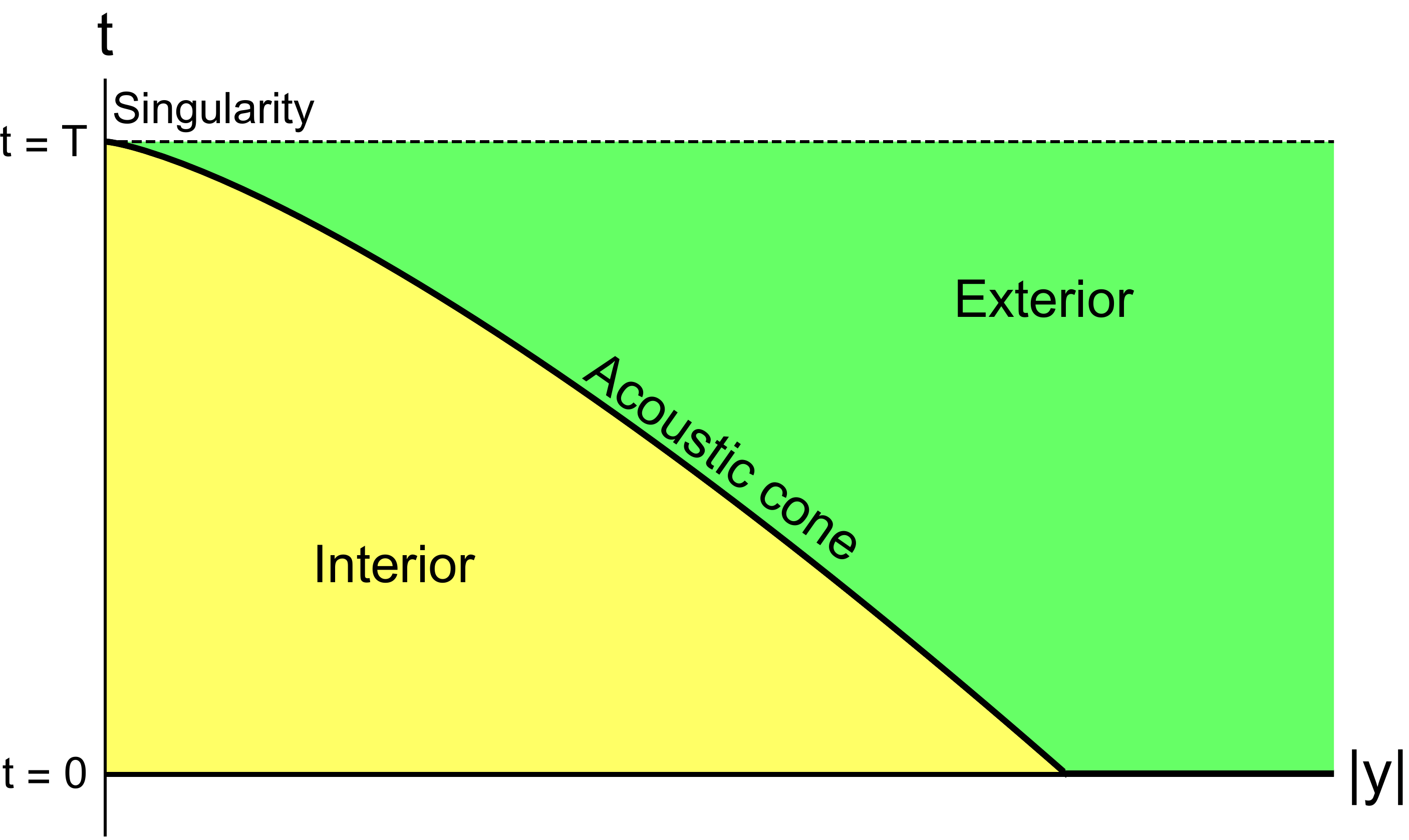}
		\caption{\small $(|y|,t)$-diagram of self-similar solutions (radius,time). At $t=0$ we have the initial data that implodes toward the origin forming a singularity at $t=T$. These solutions extend beyond the backward acoustic cone of the singularity, represented by a black solid line.}
		\label{fig:Acoustic_cone}
	\end{figure}

	In our setup, the isentropic compressible Euler equations, Guderley's construction  leads to three kinds of self-similar profiles. They extend beyond the backward acoustic cone of the singularity (see fig.~\ref{fig:Acoustic_cone}), and their regularity on this cone allows us to classify them as follows: solutions containing a shock, solutions with higher regularity but non-smooth (we will denominate them {\em non-smooth solutions} (NSSs)) and smooth solutions (SSs). 
	The interest in sharp transitions in gases to model explosions made shock-waves the most studied case, being the central object of an extensive amount of literature. The existence of SSs was recently proved in \cite{MerleEuler} almost eighty years since the publication of \cite{Guderley}. They arise after a clever fine-tuning of NSSs; unfortunately, the proof developed in \cite{MerleEuler} does not provide the specific values of the parameters that determine these profiles. This is one of the main results in our paper.

	The interest in the family of self-similar solutions that we are going to describe is further increased by the fact that, in the appropriate regimes, these solutions emerge as the leading dynamics of the energy-supercritical defocusing nonlinear Schr\"odinger equation (NLS) \cite{MerleNLS} and compressible Navier-Stokes equations (NS) \cite{MerleNS}. Despite this fact, the singularity formation of self-similar profiles is not guaranteed in these setups. Instabilities in the context of (\ref{eq:Euler_eq}) and/or the subleading dynamics of NLS or NS may prevent the process. Nevertheless, this is a heuristic idea that considers both obstacles separately. In \cite{MerleNLS, MerleNS}, the authors made use of the smoothness of the solutions that they discovered in \cite{MerleEuler} to prove the existence of smooth solutions to NLS and NS that blow-up in finite time. They also showed that they are stable up to a finite number of unstable directions to be determined.

		 In this paper we address the construction of SSs as well as the stability of both SSs and NSSs with respect to radial perturbations, preserving the symmetry of the problem. The stability of shocks was extensively addressed in the literature \cite{Hafele,AxfordHolm,Brushlinskii,Chen,Murakami}. Our interest in the characterization of perturbations with different levels of regularity is motivated by the fact that, in physics, we cannot guarantee the smoothness of a phenomenon even when the initial data is smooth. In fluid dynamics the interaction with the environment is a source of non-smooth perturbations coming from imperfections on the boundaries or from the dynamical loss of regularity. Furthermore, the study of non-smooth perturbations has a practical utility when numerical methods are involved; a practice that is rapidly increasing. The drawback of these methods, the numerical noise, excites the full spectrum of perturbations without distinguishing between smooth or non-smooth solutions or modes. For this reason, to quantify the limitations of numerical methods we require the characterization of the full spectrum.
	 
	 \subsection{Main Results}
	 
	We provide the first construction of SSs and their smooth linear perturbations. We have developed a numerical strategy to construct them in the regime where NSSs have low regularity; here the number of unstable directions is small, giving the possibility of stable SSs. However, all SSs have at least an unstable smooth mode (excluding artificial\footnote{This kind of problems has unstable modes associated with a choice of self-similar coordinates.} instabilities). Furthermore, performing a numerical time-evolution we find that the endpoint of this unstable direction is the formation of a shock before the blow-up at the origin. Therefore, the singularity formation associated with SSs is unstable\footnote{We say that a singularity formation is stable if small smooth perturbations to the initial data do not prevent the process or change the structure of the singularity, up to small deviations in its parameters; otherwise it is said to be unstable.}. 
	
	
	\section{Self-Similar Solutions\vspace{3mm}}\label{sec:Setup_Summary_Literature}
	
	This section provides a general derivations of self-similar solutions to (\ref{eq:Euler_eq}), referring the reader to \cite{MerleEuler} for further and technical details. Our notation mainly follows \cite{MerleEuler} to ease the read of both papers. Alternative introductions to the Guderley problem can be found in \cite{Guderley,Meyer, Landau, GuderleyReview, Ponchaut}. These solutions are spherically symmetric self-similar profiles of the form
		\beq
	\begin{cases}	
		\rho(t,y) = \left(\frac{1}{T-t}\right)^{\frac{2}{\gamma-1}\left(\frac{r-1}{r}\right)}\hat\rho\left(\frac{|y|}{r\left(T-t\right)^{1/r}}\right),\\
		u(t,y) = \left(\frac{1}{T-t}\right)^{\left(\frac{r-1}{r}\right)}\hat u\left(\frac{|y|}{r\left(T-t\right)^{1/r}}\right),
	\end{cases}
	\eeq
	that blow-up at the center of symmetry at time $T$ and $r>1$ is a continuous parameter denominated the blow-up speed. These profiles extend beyond the backward acoustic cone of the singularity (see fig.~\ref{fig:Acoustic_cone}); this fact will play a crucial role in this paper. The cone is given by the equation
	\beq
		\frac{d|y|}{dt} = - (u\pm c), \qquad c = \sqrt{\frac{\partial p}{\partial\rho}},
	\eeq
	with $c$ being the speed of sound.
	
	To construct self-similar solutions we introduce the self-similar variables $(\tau,Z)$
	\beq
	Z := \frac{|y|}{r(T-t)^{\frac{1}{r}}}, \qquad \left(T-t\right) := T e^{-r\tau},
	\label{eq:Self-Similar_Variables}
	\eeq
	(note the synchronization $t=0\leftrightarrow \tau=0$), where the blow-up time $T$ corresponds to $\tau = \infty$ and equations for $\hat\rho$ and $\hat u$ take the form ($\ ' := \partial_Z$)
	\begin{align}
	&\partial_\tau\hat\rho + \ell(r-1)\hat\rho+ Z \hat{\rho}' + \left(\hat\rho \hat u\right)'+\frac{d-1}{Z}\hat\rho \hat u = 0, \label{eq:Rho_dtau}\\
	&\partial_\tau \hat u + (r-1)\hat u+ Z\hat{u}'+ \hat{u}\hat{u}'+\left(\hat{\rho}^{\gamma-1}\right)'=0,
	\label{eq:U_dtau}
	\end{align}
	with  $\ell=2/(\gamma-1)$. Self-similar solutions arise in the case that $\hat\rho$ and $\hat u$ do not depend on $\tau$, reducing the problem to a system of ODEs
	\begin{align}
		&\ell(r-1)\hat\rho+ Z \hat{\rho}' + \left(\hat\rho \hat u\right)'+\frac{d-1}{Z}\hat\rho \hat u = 0, \label{eq:Rho_static}\\
		&(r-1)\hat u+ Z\hat{u}'+ \hat{u}\hat{u}'+\left(\hat{\rho}^{\gamma-1}\right)'=0.
	\label{eq:U_static}
	\end{align}
	After applying the Emden transform \cite{MerleEuler}
	\beq
	Z = e^x, \qquad \hat\rho = \left(\sqrt{\frac{\ell}{2}}Z\s(x)\right)^{\ell}, \qquad \hat u = - Z\w(x),
	\label{eq:Emden_transform}
	\eeq
	equations (\ref{eq:Rho_static})-(\ref{eq:U_static}) are written in terms of an {\em autonomous}~system~of~ODEs~determined~by~$(d,\ell,r)$
	\beq
	\frac{d\w}{dx} = - \frac{\Delta_1(\w,\s)}{\Delta(\w,\s)}, \qquad 	\frac{d\s}{dx} = - \frac{\Delta_2(\w,\s)}{\Delta(\w,\s)},
	\label{eq:System_W_S}
	\eeq 
	with
	\begin{align}
	& \Delta = \left(\w-1\right)^2-\s^2,\\
	& \Delta_1 = \w (\w -1)(\w-r)-d(\w-\frac{\ell}{d}(r-1))\s^2,\\
	& \Delta_2 = \frac{\s}{\ell}\left[(\ell+d-1)\w^2-\left(\ell+d+(\ell-1)r\right)\w+\ell r-\ell \s^2\right].
	\end{align}
	Note that $\Delta,\Delta_1,\Delta_2$ only depend on $\sigma,\omega$ and the parameters $(d,\ell,r)$.
	
\begin{figure}[h!]
		\vspace{0.5cm}
		\centering	
		\includegraphics[width=13cm]{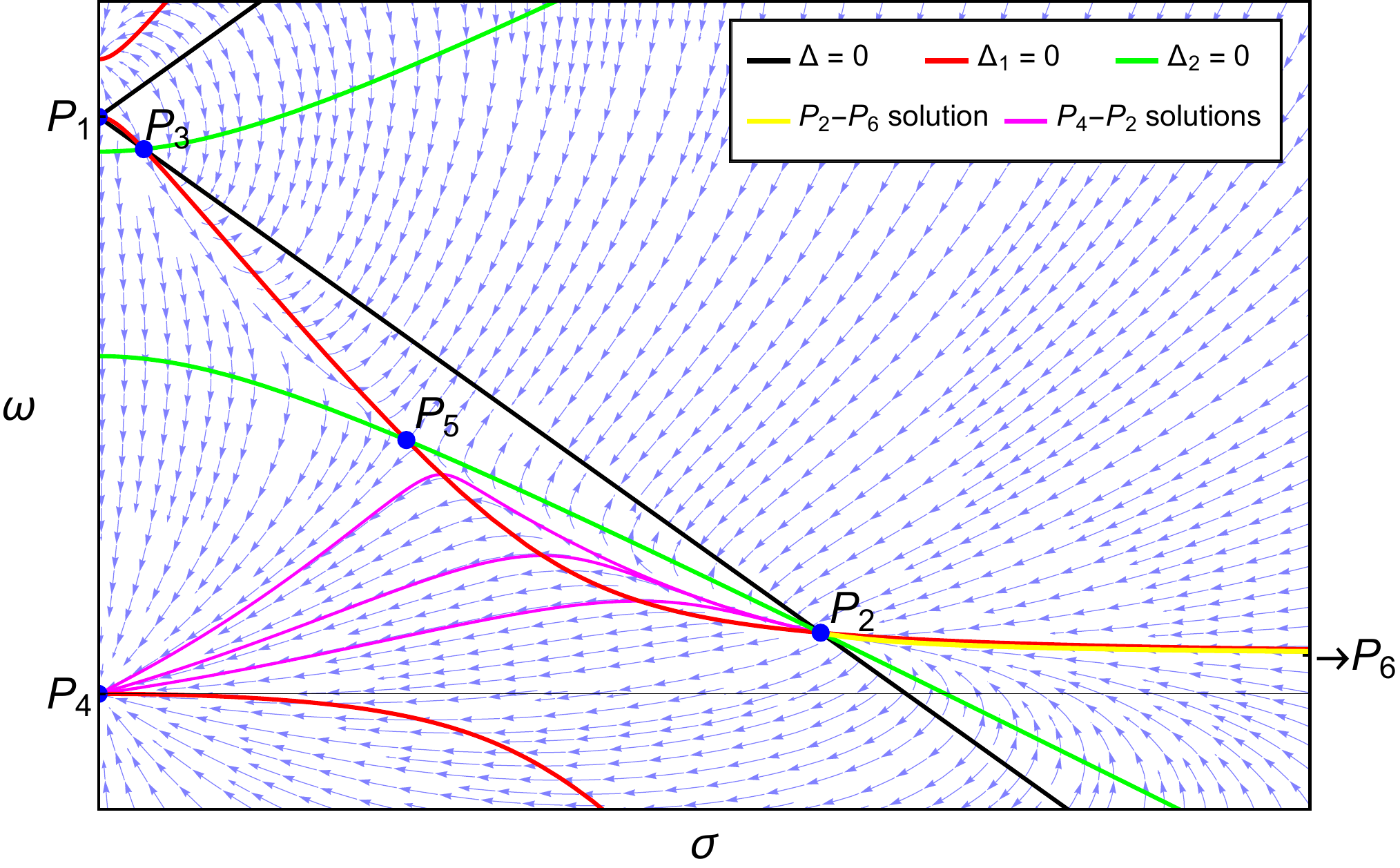}
		\caption{\small Phase portrait of system (\ref{eq:System_W_S}) with parameters $(d,\ell,r)=(3,2,1.1)$. Blue arrows represent the vector field, the yellow line the unique solution connecting $P_2$-$P_6$ and the purple ones examples of solutions connecting $P_2$-$P_4$.\\}
		\label{fig:Phase_Portrait}
	\end{figure}

	With the autonomous system (\ref{eq:System_W_S}) on our hands the construction of self-similar solutions is reduced to the analysis of the phase-portrait $(\s,\w)$. An example of this diagram is provided in fig.~\ref{fig:Phase_Portrait}. First we have to identify the location of singular points of system (\ref{eq:System_W_S}) because trajectories $\omega(\sigma)$ only begin or end at these points. Once they are identified we need to determine the set of these points that are associated with the radial origin ($Z=0$) and the radial infinity ($Z=\infty$). It allows us to locate the set of singular points where our trajectories begin or end in relation to the conditions that we want to impose on our solutions at the origin and infinity. After that, we have to find trajectories $\omega(\sigma)$ that connect both edges. Fig.~\ref{fig:Phase_Portrait} shows that this task is not simple because the phase portrait has a complex structure; fortunately, it was deeply analyzed in \cite{MerleEuler} restricted to the range of parameters
	\beq
	d\geq 2,\qquad \ell>0,\qquad 1<r<\reye = \begin{cases}
		r^* & \text{for}\  \ell < d\\
		r^+ & \text{for}\  \ell > d
	\end{cases}
	\label{eq:parameters}
	\eeq
	with
	\beq
		 r^*=\frac{d+\ell}{\ell+\sqrt{d}}, \qquad r^+=1+\frac{d-1}{ \left(1+\sqrt{\ell}\right)^2}, \qquad 1<r^*\leq r_+ \quad (r^*=r_+\ \text{only for}\ \ell=d).
	\eeq
	
	 In this range of parameter the phase-portrait has six singular points denoted by $P_{1\leq i\leq 6}$ in relation to three curves, $\Delta=0$, $\Delta_1=0$ and $\Delta_2=0$ (fig.~\ref{fig:Phase_Portrait} provides a visual representation):
	\begin{itemize}
		\item Sonic lines: $\w=\pm\s+1\ \Rightarrow\ \Delta=0$. These lines are the materialization of the backward acoustic cone of the singularity. Along them the RHS of system (\ref{eq:System_W_S}) is singular. It means that trajectories $(\sigma(Z),\omega(Z))$ lose their regularity when they reach these lines. The only opportunities to cross sonic lines preserving some level of regularity are points where the RHS is regularized, namely points where $\Delta=\Delta_1=\Delta_2 =0$.
		\item $P_4$: $(\s,\w)=(0,0)$. We want that our trajectories $\w(\s)$ end at this point because this is an attractor that corresponds to  $Z=\infty$. From (\ref{eq:Emden_transform}) we can see that $P_4$ is the only point where the density and velocity vanish or go to constant values at the spatial infinity. Furthermore, trajectories that end at this point behave like (given $(d,\ell,r)$)
		\beq
		\w(Z) \sim \kappa\frac{\eta}{Z^r} + C(\kappa) \frac{\eta^2}{Z^{2r}} + ... \qquad
		\s(Z) \sim \frac{\eta}{Z^r} + \tilde{C}(\kappa) \frac{\eta^2}{Z^{2r}} + ...
		\label{eq:expansion_ws_boundary}
		\eeq
		where $\kappa$ labels curves ending at $P_4$ and $\eta$ represents the symmetry $Z\to Z/\eta$. Hence, for $r>1$ we have our desired boundary conditions $\rho(\infty)=u(\infty)=0$.
		\item $P_6$:  $(\s,\w)=(\infty,\frac{\ell}{d}(r-1))$. We want that our trajectories $\w(\s)$ start at this point because this is a saddle point that corresponds to the origin $Z=0$. The structure of these curves close to $P_6$ is of the form
		\beq
			\s(Z) \sim \frac{s_{-1}}{Z} + s_1 Z+... \qquad \w(Z) \sim \frac{\ell}{d}(r-1) + \w_2 Z^2+...
		\eeq
		 then, from (\ref{eq:Emden_transform}) we see that trajectories that begin at $P_6$ satisfy $\rho(0)>0$ and $u(0)=0$.
		\item $P_1 = (0,1)$, $P_2 = (\s_2,1-\s_2)$ and $P_3=(\s_3,1-\s_3)$. These singular points are on the sonic line $\w=1-\s$ and satisfy ${\Delta=\Delta_1=\Delta_2=0}$. $\s_3$ and $\s_2$ are determined by these three conditions and $\s_3\leq\s_2$. Trajectories that either end or begin at these points are associated with a finite radius $0<Z<\infty$. For this reason, these points are the only options to construct trajectories that connect the spatial origin at $P_6$ with the spatial infinity at $P_4$ crossing the sonic line but preserving some level of regularity.  In particular $P_2$ will be this intermediate point as we will explain later, playing a central role in this paper. We denote it by $P_2,\ \s_2$ (in $\sigma$) and $Z_2$ (in Z).
		\item $P_5$: $(\s,\w) = \left(\frac{r\sqrt{d}}{\ell+d},\frac{r\ell}{\ell+d}\right)$. At this point $\Delta_1=\Delta_2=0$ but $\Delta\neq0$ ($\Delta = 0$ only when $r=r^*$ and $l\leq d$). The position of this point relative to the sonic line depends on the values of the parameters. In the range where we will construct SSs $(1<r<r^*)$ it is always below this singular line. See \cite{MerleEuler} for a detailed description.
	\end{itemize}

	After this analysis of the phase portrait the problem of constructing self-similar solutions, that vanish at infinity and with no vacuum regions, is reduced to the construction of trajectories $\w(\s)$ that connect the origin ($P_6$) with the infinity ($P_4$) crossing the acoustic cone (sonic line) through either $P_1,P_2$ or $P_3$. From fig.~\ref{fig:Phase_Portrait} we can see that there are no trajectories that connect $P_6$ with either $P_1$ or $P_3$. Then, our only option is $P_2$. The idea consists of constructing trajectories $P_2$-$P_6$ and $P_2$-$P_4$, gluing them at the common point $P_2$ and studying how this connection depends on the parameters $(d,\ell,r)$. This analysis can be found in \cite{MerleEuler} (Lemma 1.1):
	\begin{enumerate}
		\item There is a unique (up to scaling) trajectory connecting $P_2$-$P_6$ which is $\mathcal{C}^\infty$ in $Z\in[0,Z_2)$.
		\item There exists a one-parameter family of trajectories connecting $P_2$-$P_4$ which is $\mathcal{C}^\infty$ in $Z\in (Z_2,\infty)$. 
		\item Trajectories $P_2$-$P_4$ can be glued with the trajectory $P_2$-$P_6$ at $P_2$ with $0<Z_2<\infty$.
	\end{enumerate}
	Now we are going to describe the process of connecting the unique curve $P_2$-$P_6$ with some trajectories $P_2$-$P_4$. First one needs to know that the slope of any trajectory that ends/begins at $P_2$ is
	\beq
		\w_{\pm} = \frac{c_4-c_1\pm\sqrt{(c_1-c_4)^2+4c_2c_3}}{2|c_2|}	
		\label{eq:slope_w_pm}
	\eeq
	where $c_i$ are the entries of the Jacobian matrix of the vector field of system (\ref{eq:System_W_S})
	\beq
	J = \begin{pmatrix}
		c_1 & c_3\\
		c_2 & c_4
	\end{pmatrix},\qquad \begin{cases}
		c_1 = \partial_\w \Delta_1\big{|}_{P_2}, &	c_3 = \partial_\w \Delta_1\big{|}_{P_2},\\ 
		c_2 = \partial_\s \Delta_2\big{|}_{P_2}, & c_4 = \partial_\s \Delta_2\big{|}_{P_2}.
		\label{eq:c_i_coefficients}
	\end{cases}
	\eeq
	In particular, the slope of the unique trajectory $P_2$-$P_6$ is $\w_-$ and therefore, in our future construction of smooth solutions we need curves $P_2$-$P_4$ that end at $P_2$ with the same slope. To explain the particularities of these curves that can be glued to $P_2$-$P_6$ with the appropriate slope, we are going to use fig.~\ref{fig:Phase_Portrait} restricted to $r<r^*$, providing a visual description of our region of interest. A technical and detailed explanation can be found in \cite{MerleEuler}. Trajectories that begin at $P_2$ (with the appropriate slope to be glued to $P_2$-$P_6$) and end at $P_4$ emerge from $P_2$ between the green and red lines in fig.~\ref{fig:Phase_Portrait} (the ``eye" structure formed by these lines between $P_5$ and $P_2$). To escape from this region trajectories must cross either the green or the red line. In the former case they do not end at $P_4$ and in the latter case they constitute our desired trajectories $P_2$-$P_4$. When $r$ increases $P_5$ moves toward $P_2$ and this window between them shrinks. Moreover, for $\ell\leq d$ when $r=r^*$ the window is closed ($P_5=P_2$). In case of $\ell>d$ a similar process happens when $r\sim r_+$, but in this case $P_3$ plays the role of $P_5$ (see \cite{MerleEuler}). From this description we see that the family of curves $P_2$-$P_4$ that we will use belong to a finite interval of the parameter $\kappa_{\min}<\kappa<\kappa_{\max}$ given in (\ref{eq:expansion_ws_boundary}).

	The final products of the process described above are spherically symmetric self-similar profiles $(\hat\rho(Z),\hat u(Z))$ that vanish at infinity with no vacuum regions. They have regularity $\mathcal{C}^{\infty}\left(\mathbb{R}\setminus\{Z_2\}\right)$ but are generically non-smooth at $Z_2$ (on the acoustic cone) because the structure of curves $P_2$-$P_6$ and $P_2$-$P_4$ at this singular point is of the form ($\xi := \sigma - \sigma_2$)
	\beq
		\omega(\xi) = \underbrace{\sum \tilde{\omega}_n\xi^n }_{\text{integer powers}} + \underbrace{c \ |\xi|^{\nu}\sum\left(...\right)}_{\text{non-integer powers}}
		\label{eq:expansion_w_in_powers_sigma_Z2}
	\eeq
	where we define $\nu := \lambda_{-}/\lambda_+$ with $\lambda_\pm$ the eigenvalues of the Jacobian matrix
	\beq
		\lambda_{\pm} = \frac{c_1+c_4\pm\sqrt{(c_1-c_4)^2+4c_2c_3}}{2}.
		\label{eq:lambda_eigenvalues}
	\eeq
	From the definitions of $c_i$ given in (\ref{eq:c_i_coefficients}) we can see that the regularity depends on the parameters $\nu=\nu(d,\ell,r)$. It can be arbitrarily large, but finite, because
	\beq
		\lambda_- < \lambda_+<0 \qquad \text{and} \qquad \nu = \frac{\lambda_-}{\lambda_+} \underset{r\sim \reye}{\sim} \frac{C(\ell,d)}{\reye-r}.
	\eeq
	Therefore, we see that manipulating the parameters of the problem we can construct NSSs of any finite regularity. One may think that for integer values of $\nu$ these self-similar solutions would be smooth. However, this is not guaranteed and logarithms appear. The only option to have a smooth connection at $P_2$ is that the coefficient $c$ vanishes. In \cite{MerleEuler} the authors proved the existence of discrete values of $r$ (and $\kappa$), that accumulate at the critical speed $\reye$, where the factor $c$ vanishes for both trajectories $P_2$-$P_6$ and $P_2$-$P_4$; proving the existence of SSs. Nevertheless, in order to construct a robust theorem based on a rigorous proof this result is valid for $r$ very close to $\reye$ ($\nu\to\infty$) and the specific values of $r$ such that $c(r)=0$ were not known. In section~\ref{sec:Blow_Up_Profiles_Linear_Perturbations} we construct SSs and provide the first known values of $r$ (and $\kappa$) associated with them. It will be done far from any limit ($1<\nu<8$), where numerical methods work appropriately. Another motivation to explore this range of parameters is that we expect that the number of unstable directions of SSs grows with $\nu$, therefore, it provides the opportunity to construct a stable SSs (unfortunately all SSs are unstable as we show in section~\ref{sec:Blow_Up_Profiles_Linear_Perturbations}).

	
	\section{Smooth Solutions and Linear Perturbations\vspace{3mm}}\label{sec:Blow_Up_Profiles_Linear_Perturbations}
	
	 In this section we construct SSs and their linear modes, leaving the characterization of endpoints for section~\ref{sec:Quasi_Shock_Formation}. We warn the reader that, in order to keep the attention on the main results, technical and numerical details are relegated to appendices~\ref{sec:Appendix_Regularity_Linear_Modes_Z2},~\ref{sec:Appendix_Numerical_Methods}.
	
		
		\begin{figure}[h!]
		\centering	
		\begin{subfigure}[b]{0.5\textwidth}
			\centering
			\includegraphics[width=7.5cm]{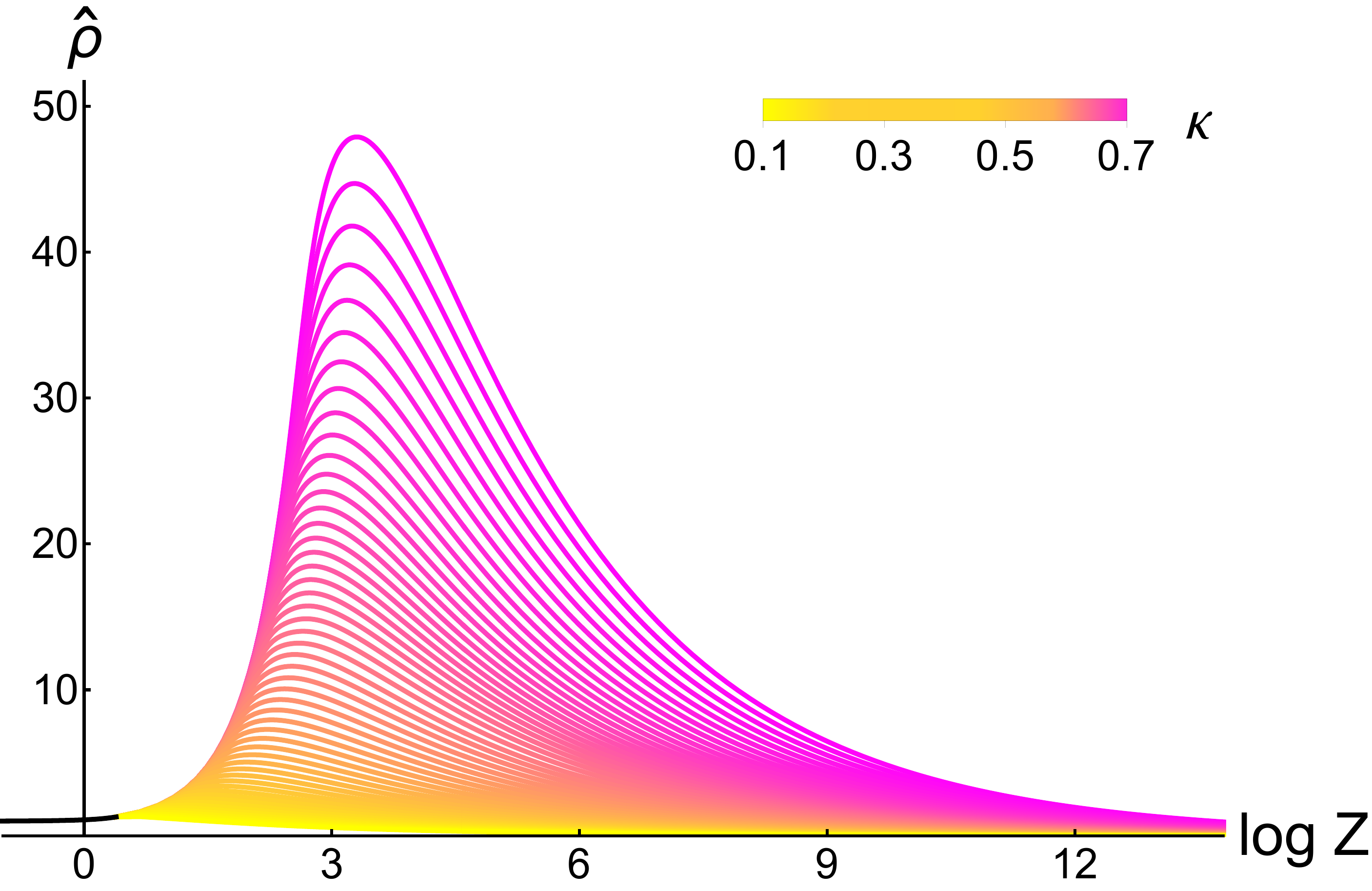}
			\label{fig:Rho_profiles_family}
		\end{subfigure}%
		\begin{subfigure}[b]{0.5\textwidth}
			\centering		\hspace{2cm}\includegraphics[width=7.5cm]{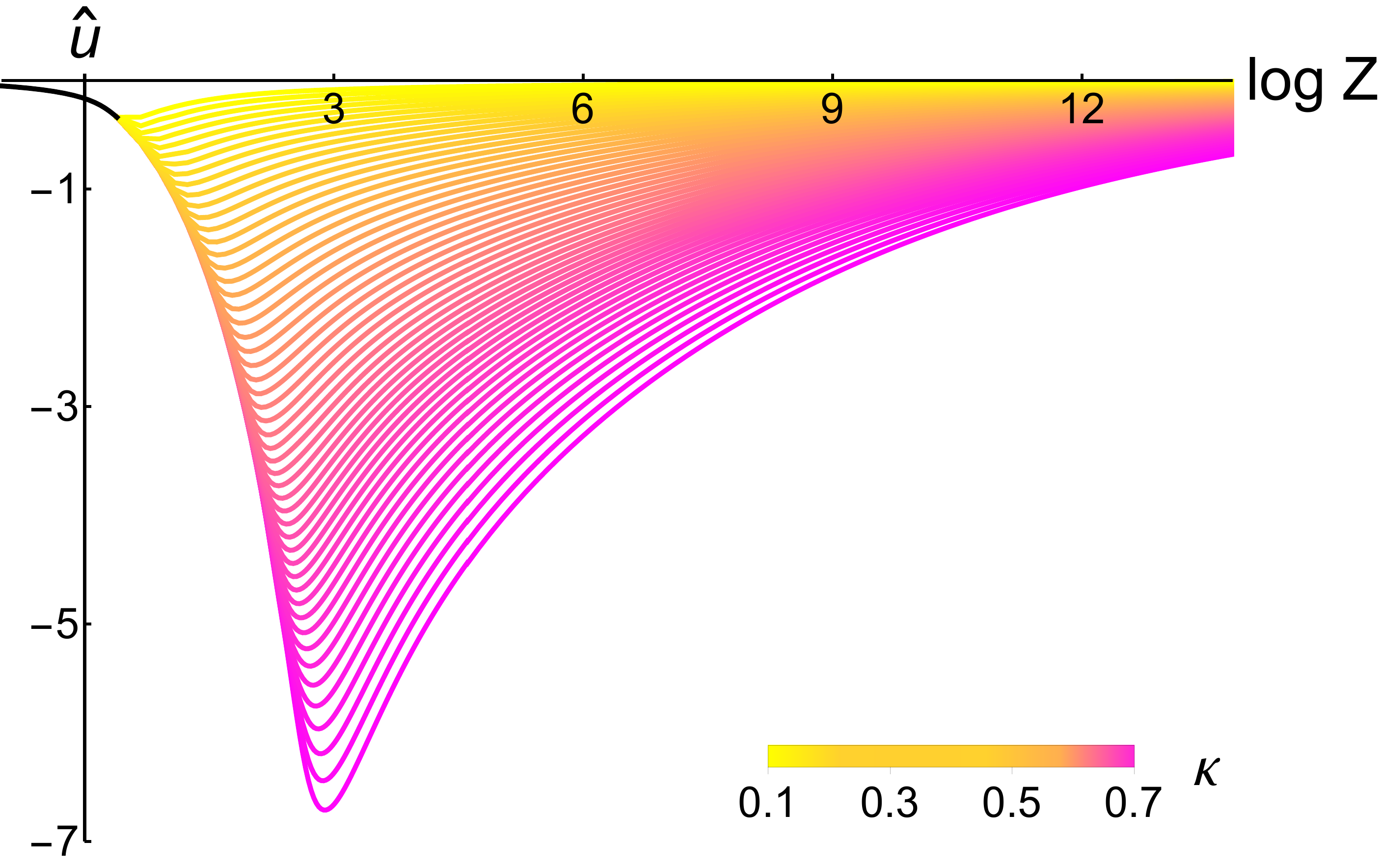}
			\label{fig:U_profiles_family}
		\end{subfigure}%
		\caption{\small Family of self-similar solutions associated with $(d,\ell,r) = (3,2,1.2)$ and parameterized by $\kappa$. The black line represents the unique solution in $Z\in[0,Z_2)$ and colored lines solutions that depend on $\kappa$ in $Z\in(Z_2,\infty)$.}
		\label{fig:Profiles_family}
	\end{figure}

		
				


	\subsection{Smooth Solutions}
		\label{sec:Smooth_Solutions}
		 
	We have seen that once $(d,\ell,r)$ are fixed, self-similar solutions belong to a one-parameter family (see fig.~\ref{fig:Profiles_family}).
	They consist of a trajectory from the spatial infinity to the acoustic cone ($P_2$-$P_4$) parameterized by $\kappa$, (\ref{eq:expansion_ws_boundary}), and the unique trajectory from the origin to the acoustic cone ($P_2$-$P_6$). Hence, a specific self-similar profile is fully determined by $(d, \ell,r,\kappa)$. At $P_2$ these trajectories $\omega(\s)$ have an expansion of the form (non-integer $\nu$, $\xi:= (\s - \s_2)$)
		\beq
		\omega(\xi) = \underbrace{\sum \tilde{\omega}_n\xi^n }_{\text{integer powers}} + \underbrace{c_\pm \ |\xi|^{\nu} \sum\left(...\right)}_{\text{non-integer powers}},
		\label{eq:expansion_w_in_powers_sigma_Z2_plusminus}
		\eeq
		where $c_+$ ($c_-$) is associated with $P_2$-$P_6$ ($P_2$-$P_4$); namely, the expansion from the right (left) of $P_2$. Then, $c_{+}$ and $c_-$ are the cause of the non-smoothness at $P_2$. One may think that for integer values of $\nu$ self-similar solution are smooth. However, this is not guaranteed and extra terms appear. For this reason we avoid these values in our construction of smooth profiles, focusing on non-integer $\nu$, where \cite{MerleEuler} provides theorems.
		The deep analysis performed in this paper guarantees that $|\xi|^{\nu}$ is the lowest non-integer power and that if $c_{\pm} = 0$ for a given self-similar profile, then this solution is smooth. 		
		
\begin{figure}[h!]
		\vspace{2cm}
			\centering	
			\begin{subfigure}[b]{0.5\textwidth}
				\centering
				\includegraphics[width=8.1cm]{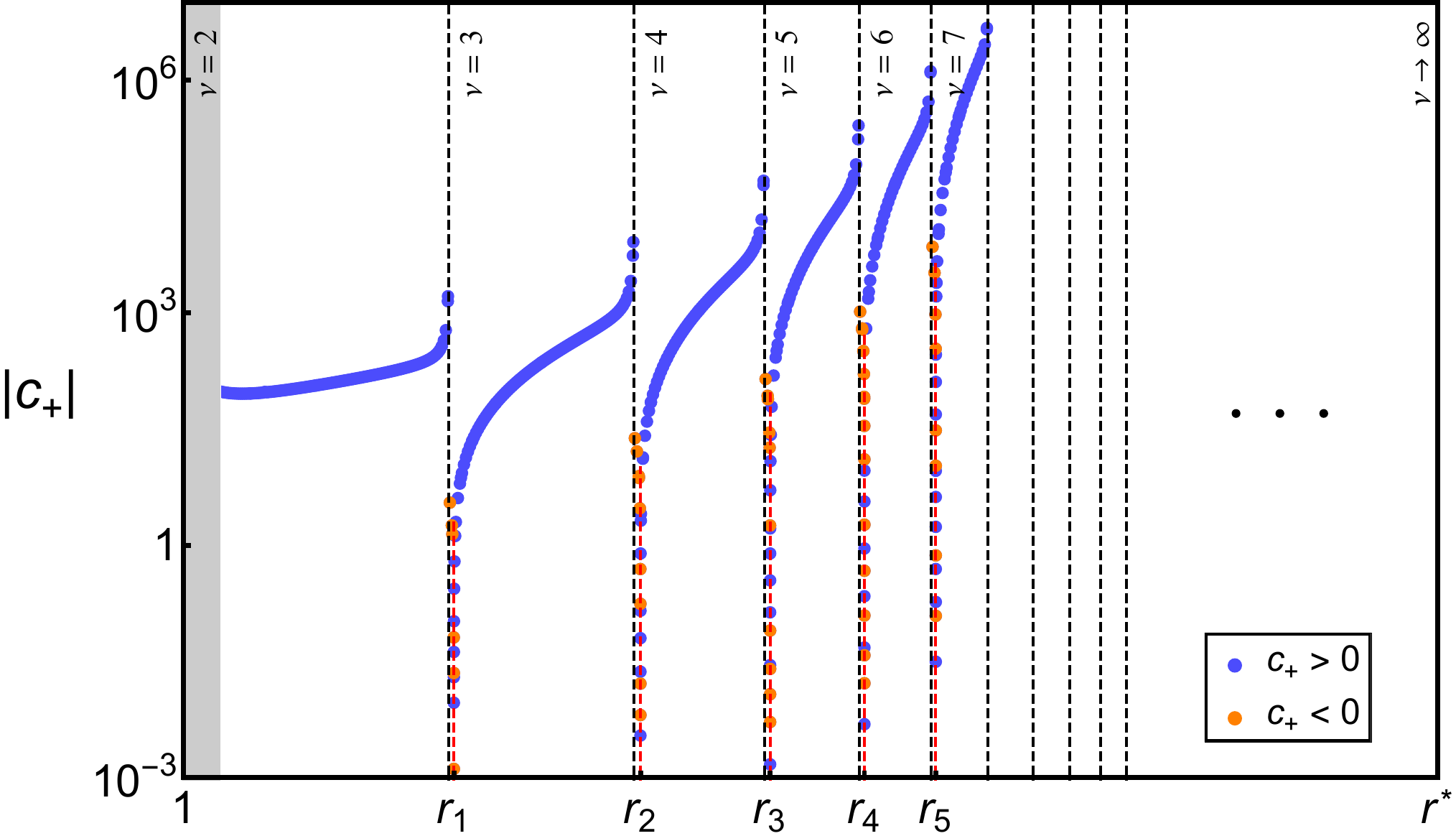}
				\caption{$(d,\ell) = (2,2)$}
				\label{fig:Search_rk_A_d2}
			\end{subfigure}%
			\begin{subfigure}[b]{0.5\textwidth}
				\centering
				\includegraphics[width=8.1cm]{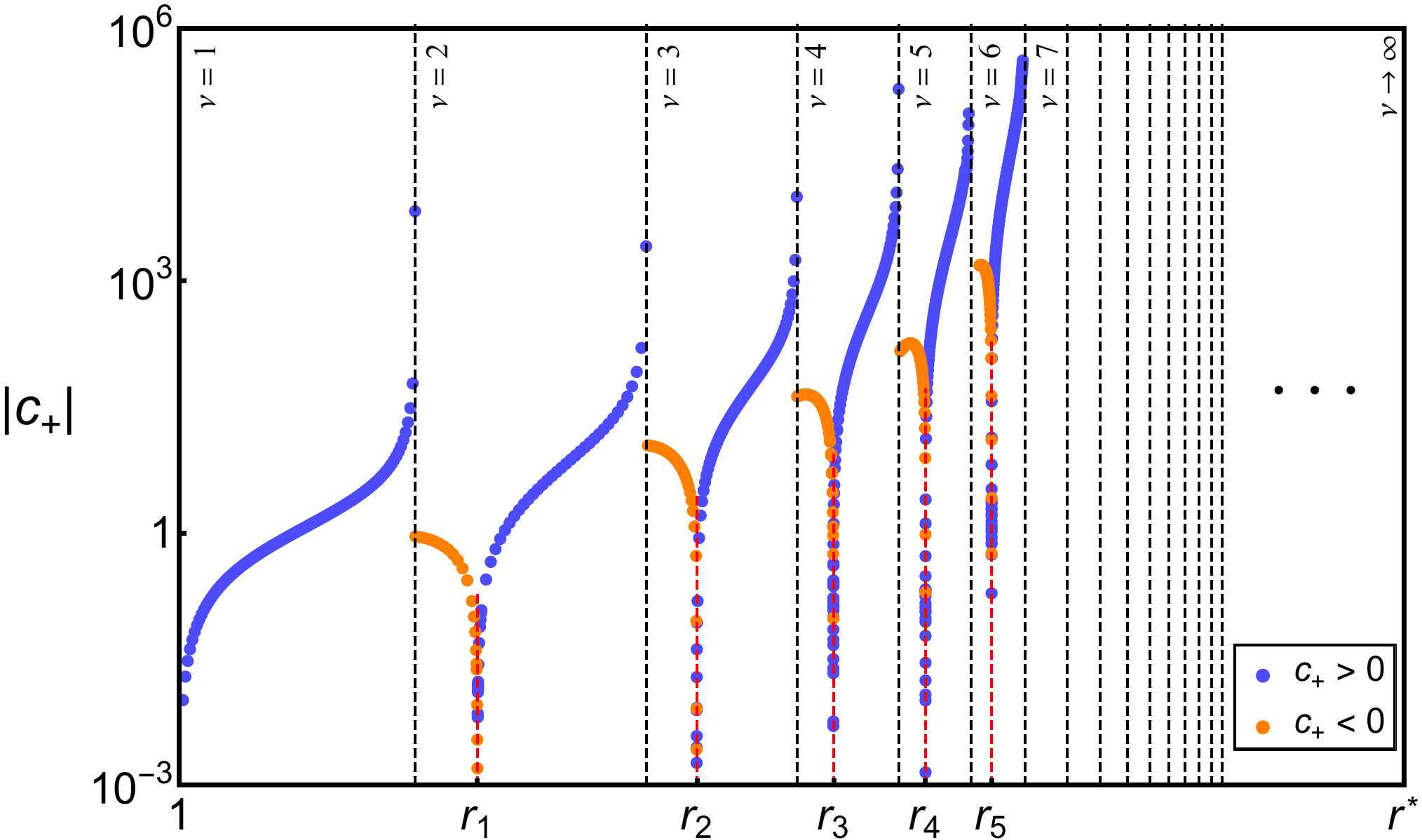}
				\caption{$(d,\ell) = (3,2)$}
				\label{fig:Search_rk_A_d3}
			\end{subfigure}%
			
			\vspace{0.5cm}
			\begin{subfigure}[b]{0.5\textwidth}
				\centering
				\includegraphics[width=8.1cm]{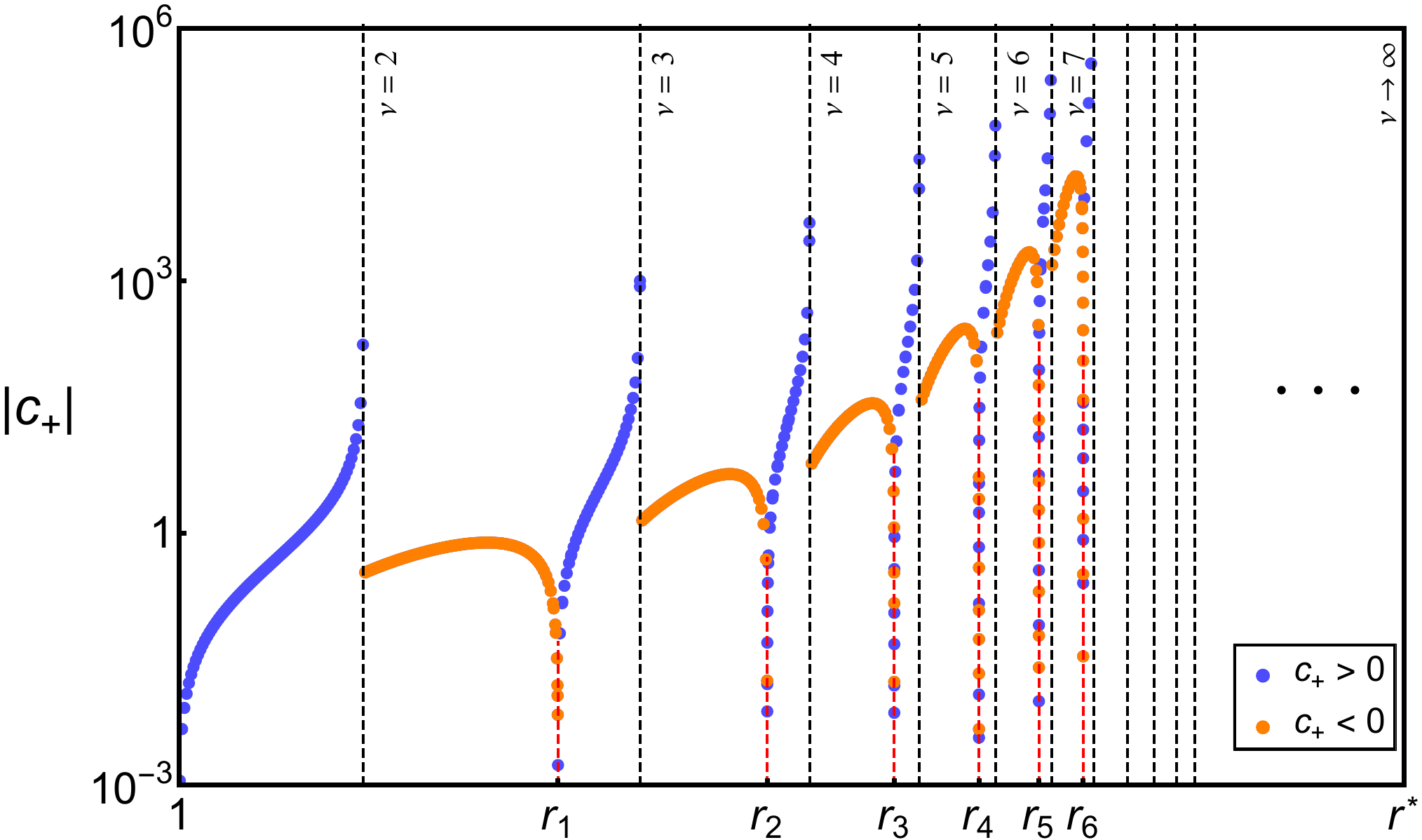}
				\caption{$(d,\ell) = (4,4)$}
				\label{fig:Search_rk_A_d4}
			\end{subfigure}%
			\begin{subfigure}[b]{0.5\textwidth}
				\hspace{0.5cm}
				\centering
				\includegraphics[width=8.1cm]{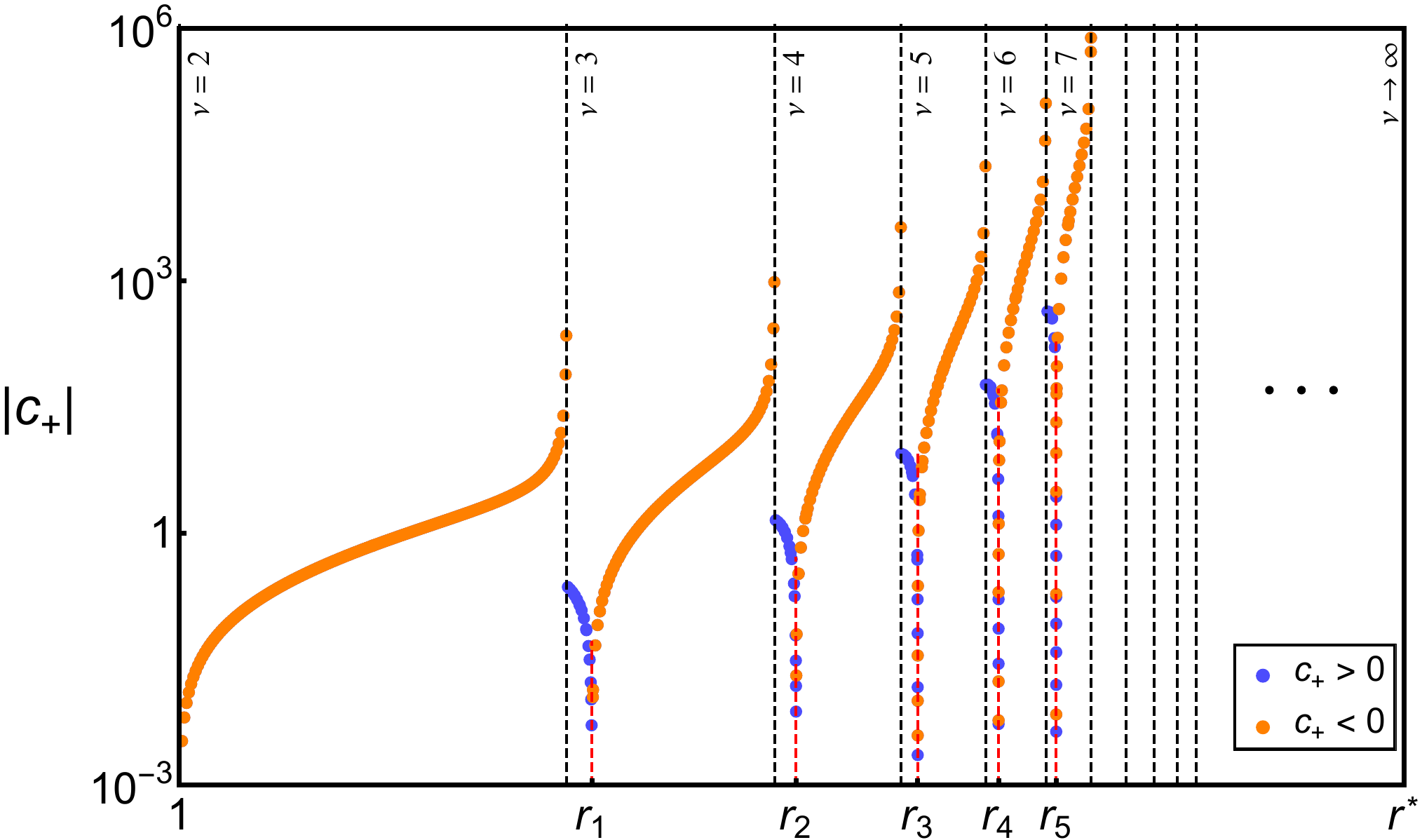}
				\caption{$(d,\ell) = (5,5)$}
				\label{fig:Search_rk_A_d5}
			\end{subfigure}%
			
			\caption{\small These plots show representative examples of the dependence of coefficient $|c_+|$ on $r$ for fixed $(d,\ell)$ (and $\nu<8$). Blue points represent $c_+>0$ and orange points $c_+<0$. Transitions between colors, red dashed lines, mark values of $r$ where $c_+$ vanishes ($r_n$), while black dashed lines mark values where $\nu$ is integer. Ellipsis mean that the structure of these plots should be extended to $\nu\to\infty$. We can observe that in $d\geq3$ when $r$ approaches $1$ coefficient $|c_+|$ goes to zero. However, in $d=2$ this is not what happens and numerical methods have troubles to describe this region, the gray area means that the quality of our results for $r\sim1$ is not good enough to report results. This difference between $d=2$ and $d\geq3$ comes from the fact that the slope $\w_-$ of trajectories $P_2$-$P_6$ (given in (\ref{eq:slope_w_pm})) goes to zero for $d\geq3$ when $r\to1$, but in $d=2$ it goes to a nonvanishing value while $w(\sigma)\to0$, producing a region close to $P_2$ with an abrupt transition. An extra zero of $c_+$ may be hidden in this region. \\ \\}
			\label{fig:Search_rk_dimensions}
		\end{figure}

\begin{figure}[h!]
		\vspace{2.5cm}
		\centering	
		\begin{subfigure}[b]{0.5\textwidth}
			\centering
			\includegraphics[width=8.1cm]{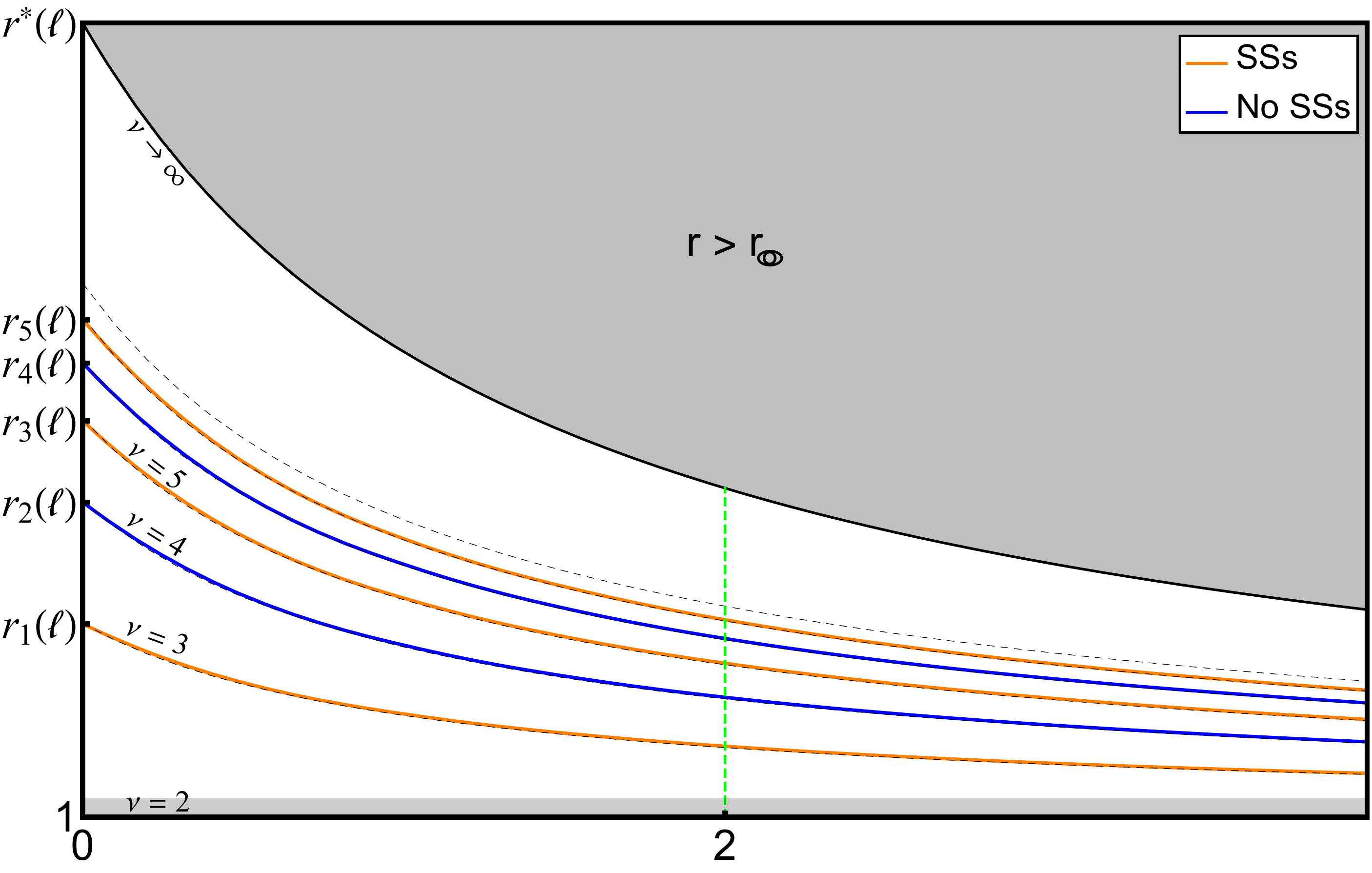}
			\caption{Diagram $(\ell,r)$ in $d=2$}
			\label{fig:Search_rk_B_d2}
		\end{subfigure}%
		\begin{subfigure}[b]{0.5\textwidth}
			\centering
			\includegraphics[width=8.1cm]{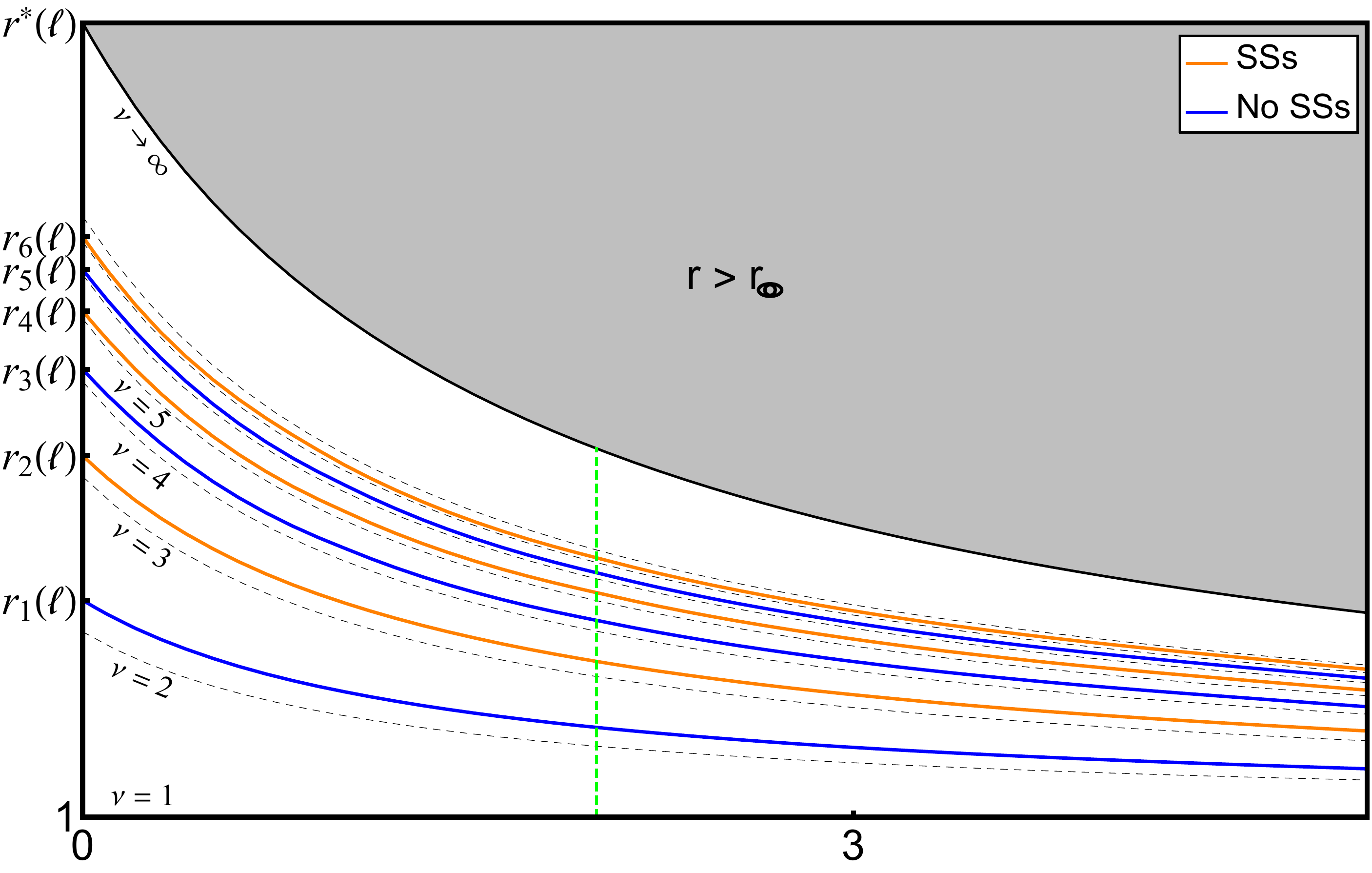}
			\caption{Diagram $(\ell,r)$ in $d=3$}
			\label{fig:Search_rk_B_d3}
		\end{subfigure}%
		
		\vspace{0.5cm}
		\begin{subfigure}[b]{0.5\textwidth}
			\centering
			\includegraphics[width=8.1cm]{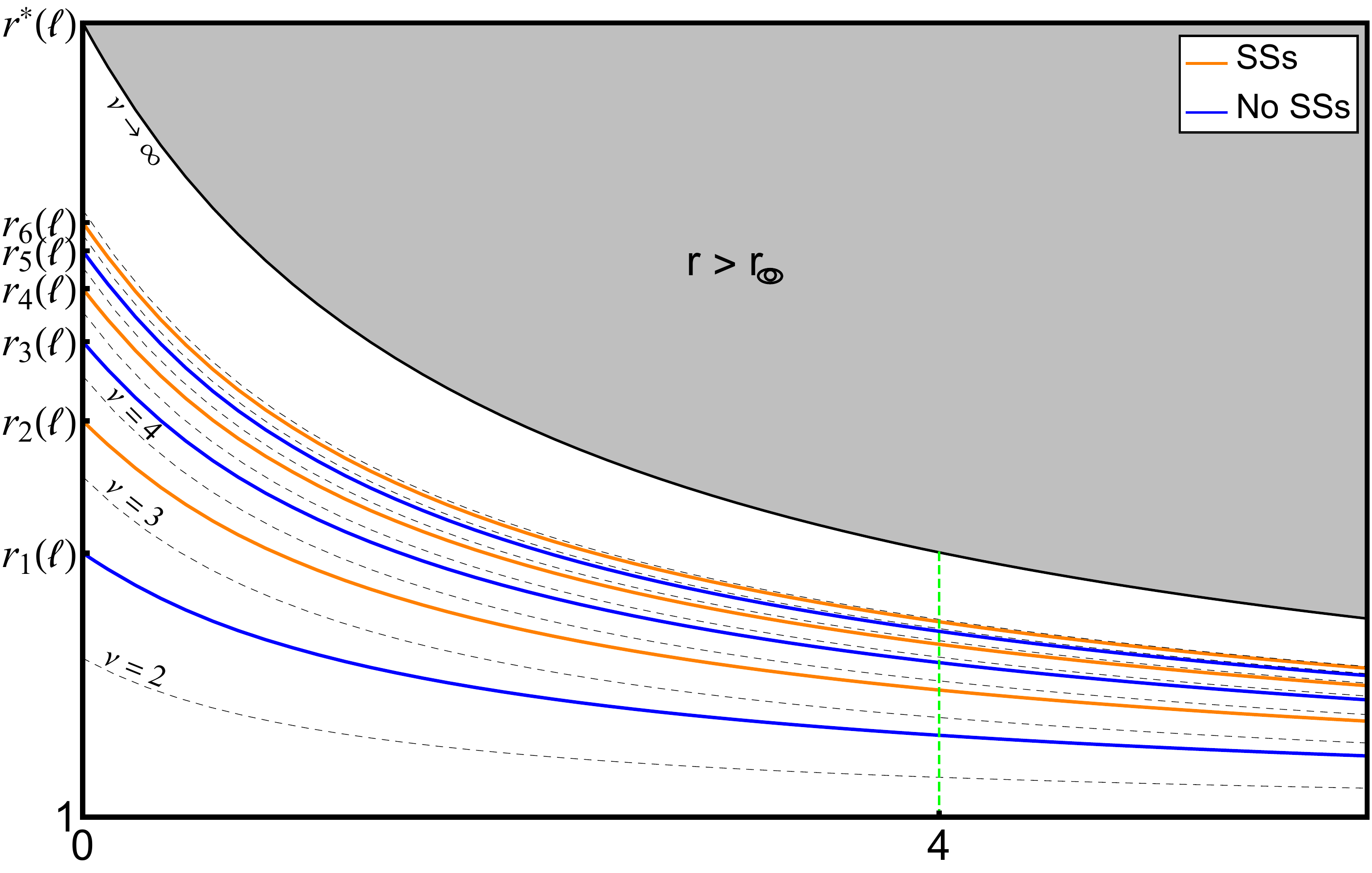}
			\caption{Diagram $(\ell,r)$ in $d=4$}
			\label{fig:Search_rk_B_d4}
		\end{subfigure}%
		\begin{subfigure}[b]{0.5\textwidth}
			\hspace{0.5cm}
			\centering
			\includegraphics[width=8.1cm]{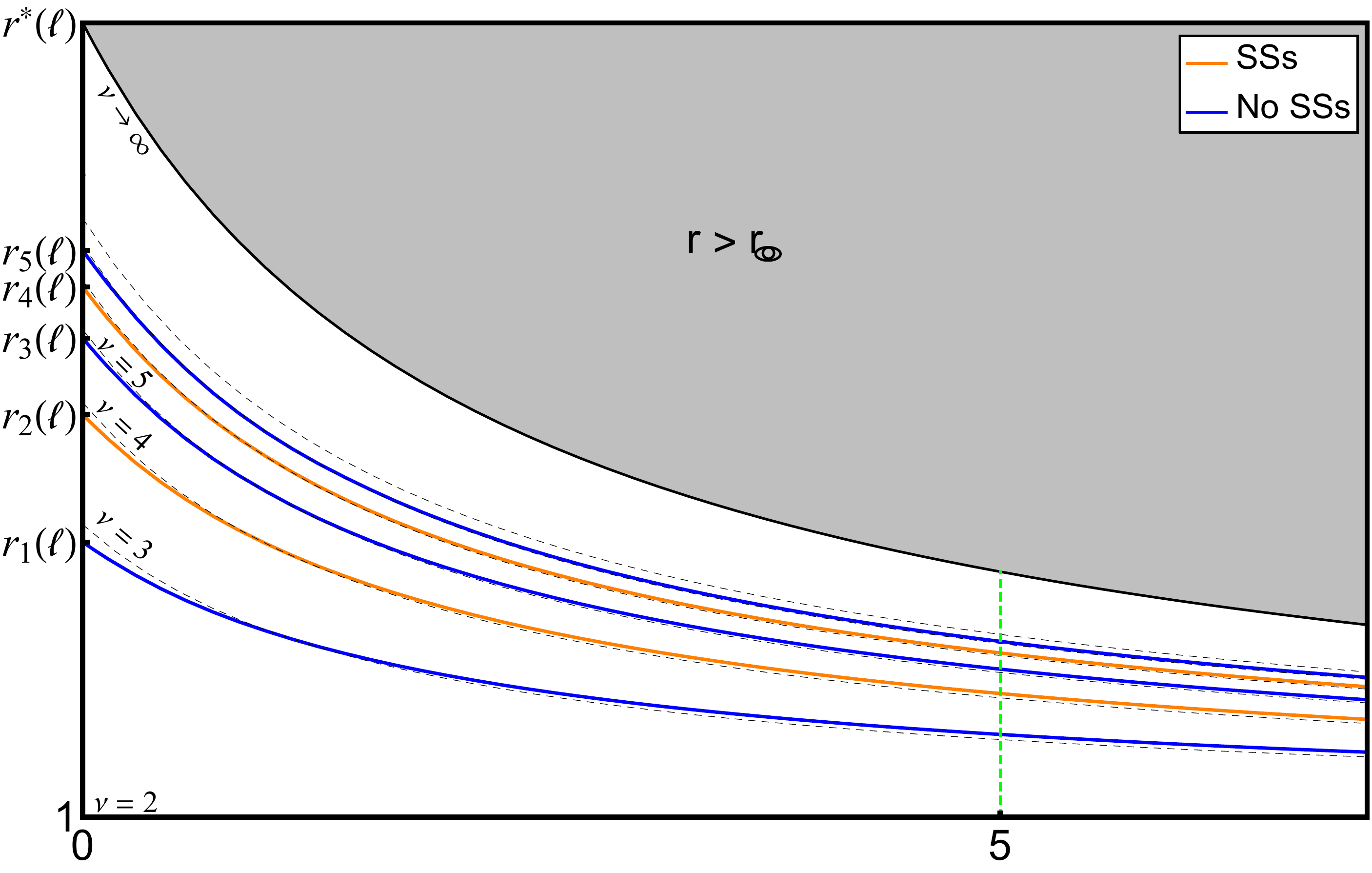}
			\caption{Diagram $(\ell,r)$ in $d=5$}
			\label{fig:Search_rk_B_d5}
		\end{subfigure}%
		
		\caption{\small Diagrams $(\ell,r)$ in different dimensions. They show the evolution of the first zeros of $c_+$ with $\ell$; namely, curves $r_n(\ell)$ (orange and blue lines). Orange and blue colors represent the parity condition described in the main text, there are no SSs for blue curves and there are SSs for orange curves. Black dashed lines are values of $(\ell,r)$ where the exponent $\nu(r,\ell)$ is integer. The gray area is the region $r>\reye$ out of our range of parameters and the green line is the location of fig.~\ref{fig:Search_rk_dimensions} in this diagram. From these plots we see that curves $r_n(\ell)$ do not intersect each other and that these lines can cross black dashed lines (see $d=5$), meaning that they transit from regions with different integer $[\nu]$. Finally, in the first plot ($d=2$) the gray area close to $r=1$ represents the region where our methods do not provide accurate access as it was explained in fig.~\ref{fig:Search_rk_dimensions}.\\ \\ \\ \\ \\}
		\label{fig:Search_rk_B}
	\end{figure}

\begin{figure}[h!]
		\vspace{1.5cm}
		\centering	
		\begin{subfigure}[b]{0.35\textwidth}
			\centering
			\includegraphics[width=5.5cm]{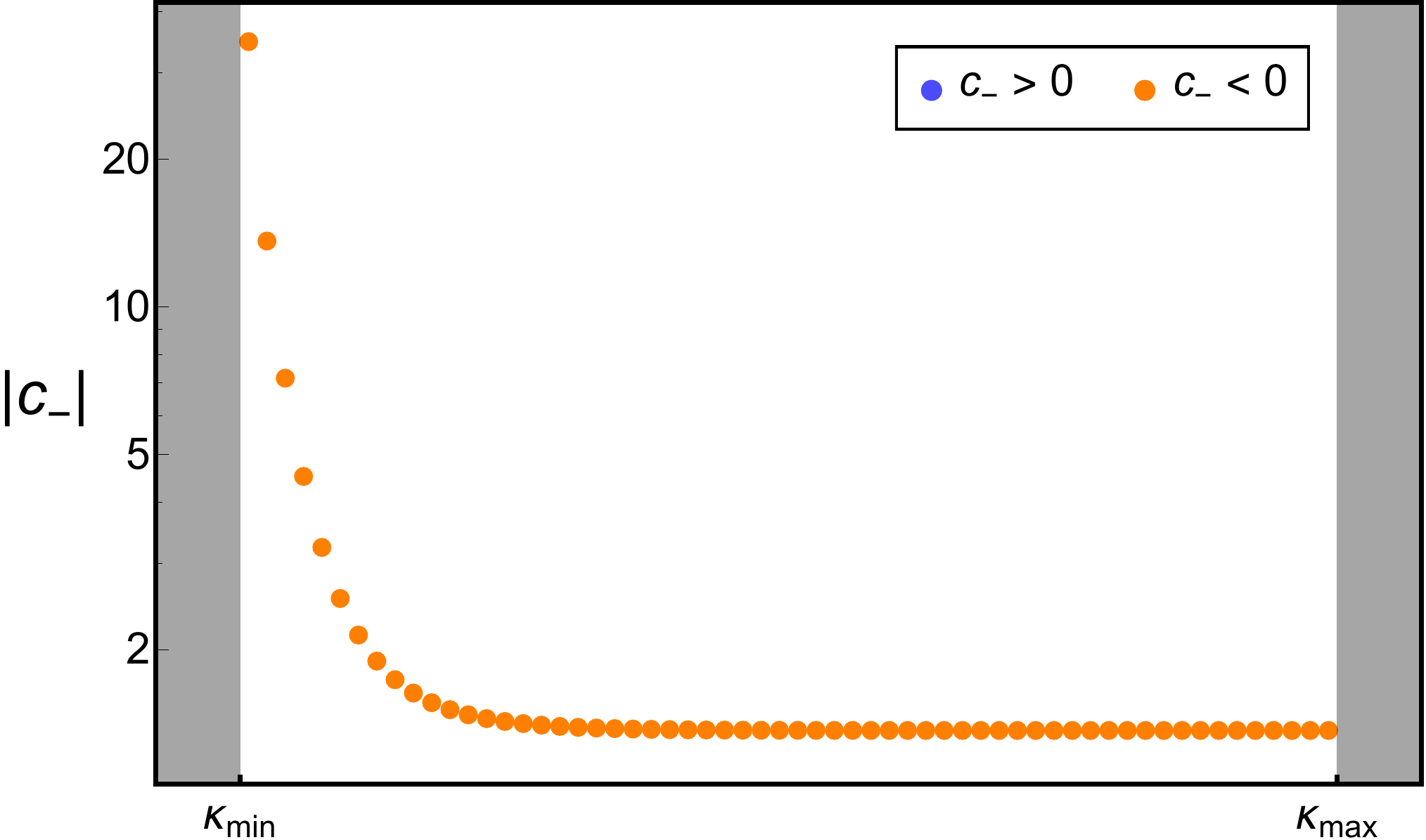}
			\caption{$(d,\ell,r) = (3,3,r_1)$}
		\end{subfigure}%
		\begin{subfigure}[b]{0.35\textwidth}
			\centering
			\includegraphics[width=5.5cm]{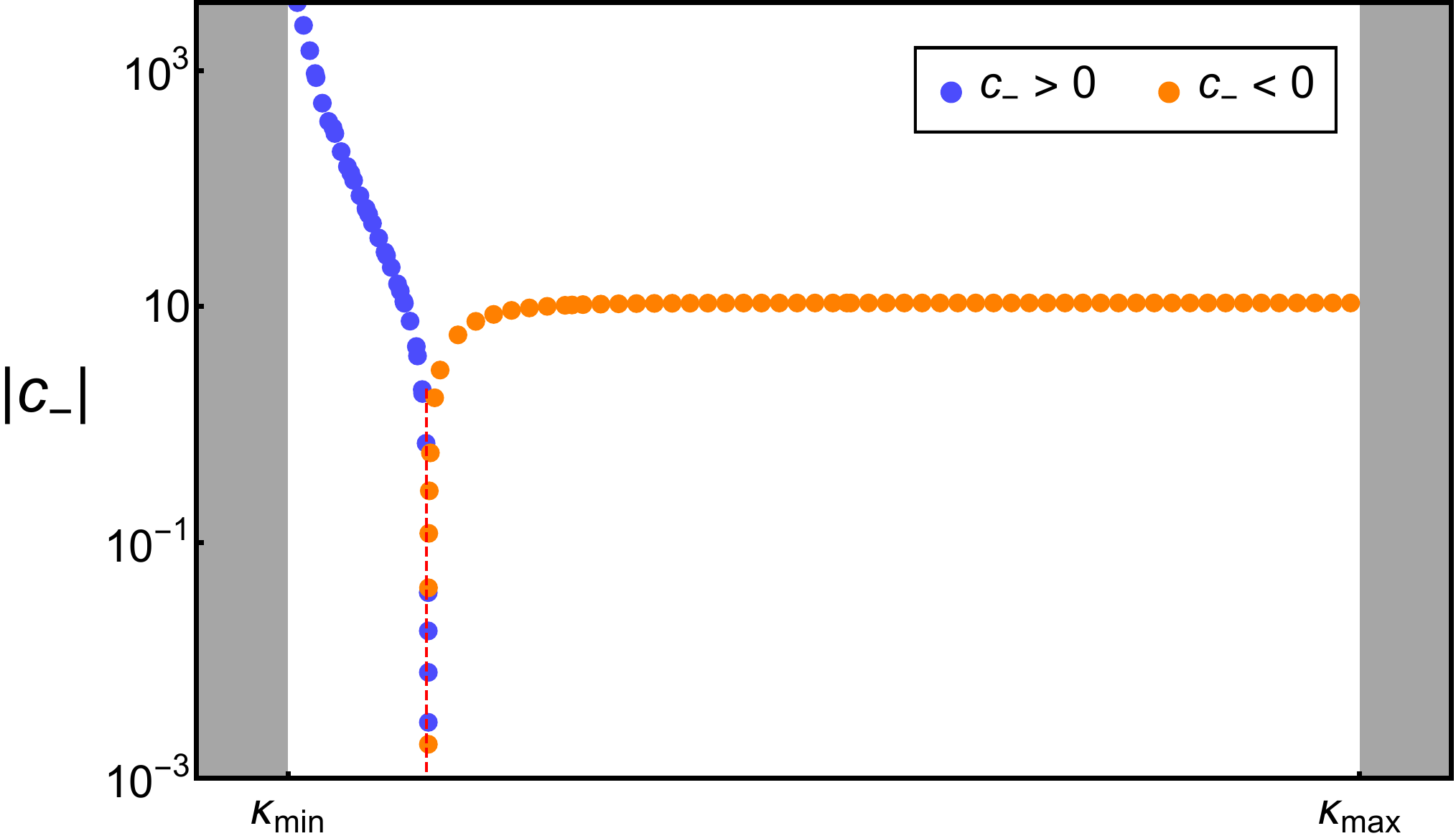}
			\caption{$(d,\ell,r) = (3,3,r_2)$}
		\end{subfigure}%
		\begin{subfigure}[b]{0.35\textwidth}
			\centering		
			\includegraphics[width=5.5cm]{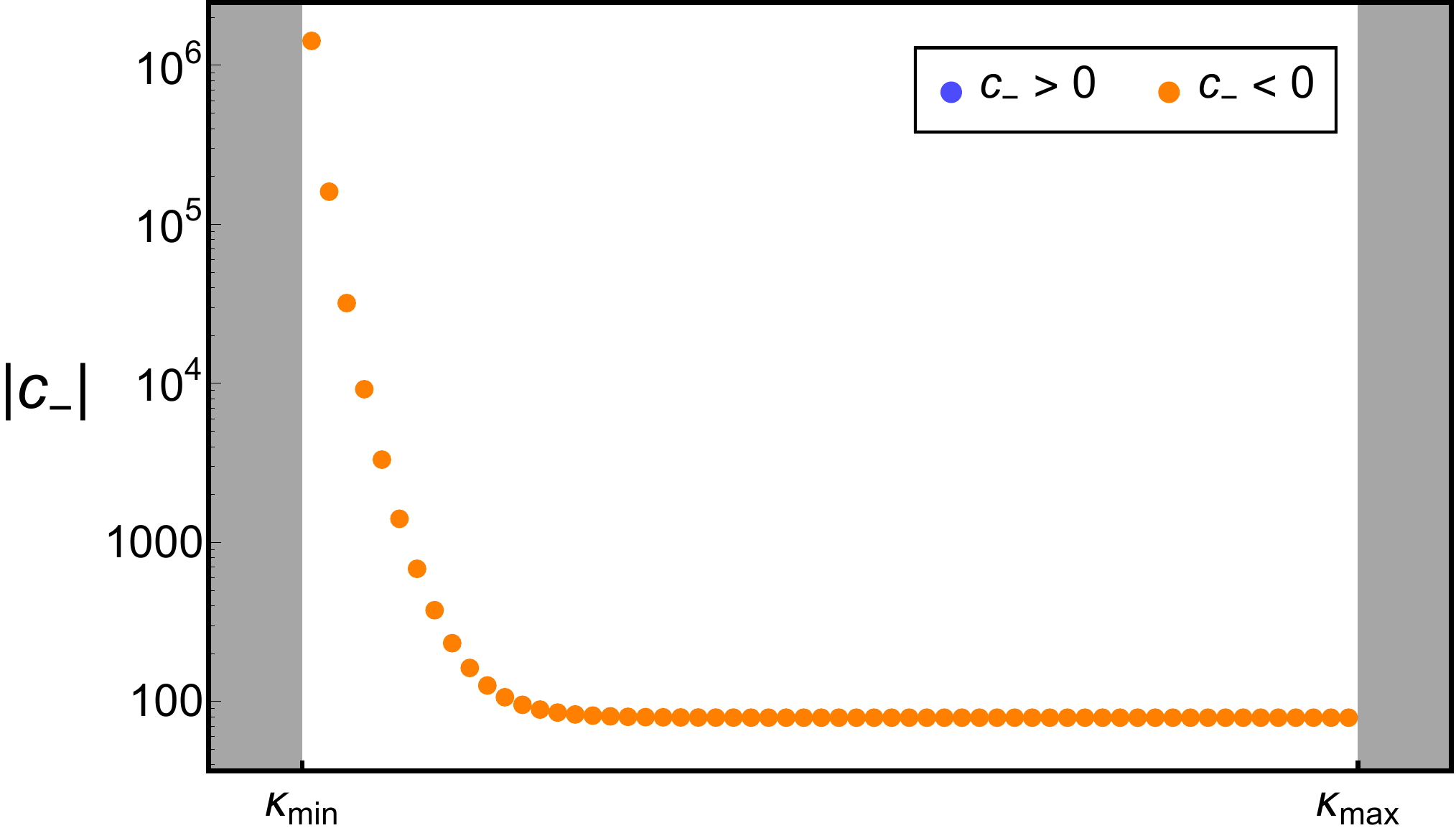}
			\caption{$(d,\ell,r) = (3,3,r_3)$}
		\end{subfigure}%
		
		\vspace{0.5cm}
		
		\begin{subfigure}[b]{0.35\textwidth}
			\centering		
			\includegraphics[width=5.5cm]{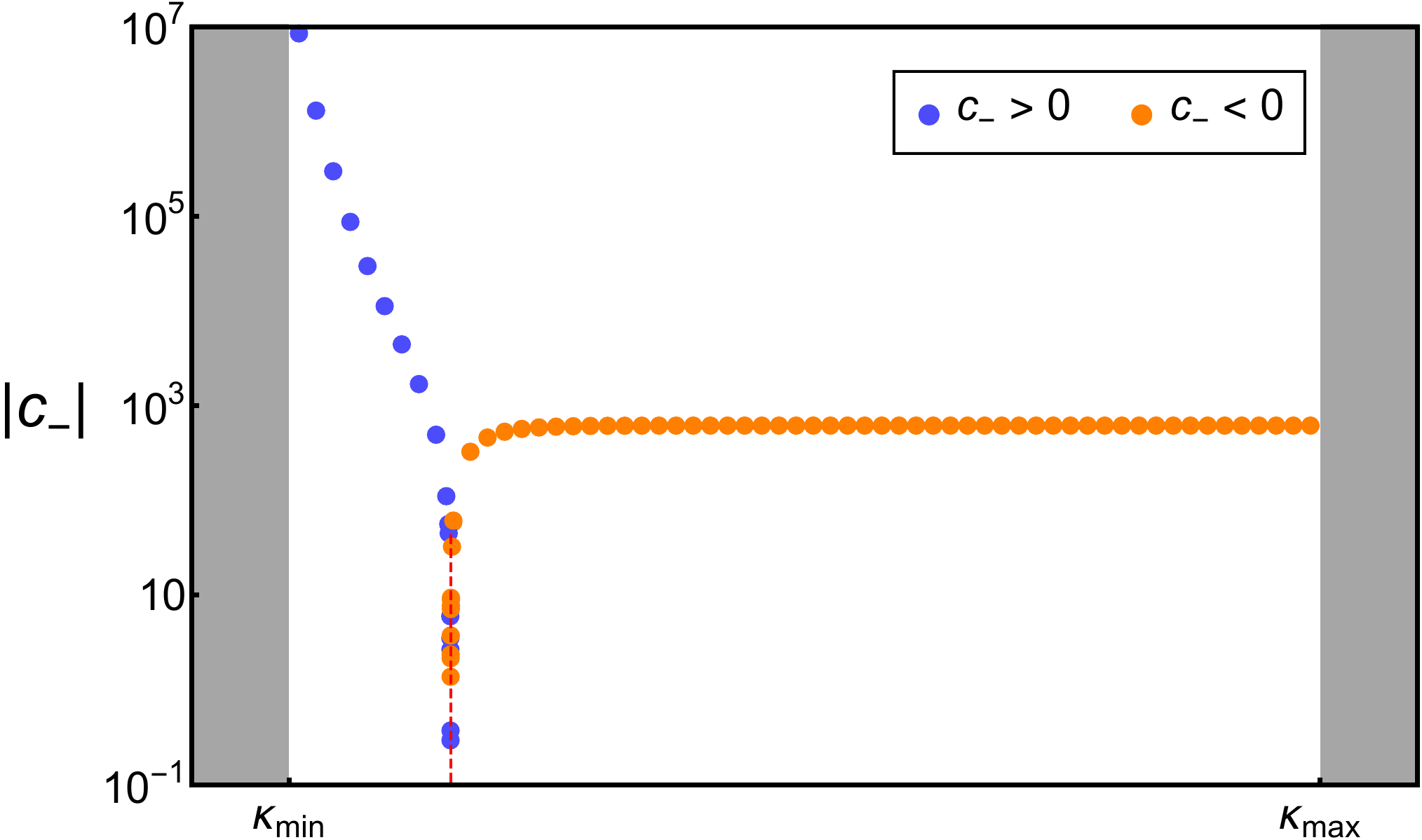}
			\caption{$(d,\ell,r) = (3,3,r_4)$}
		\end{subfigure}%
		\begin{subfigure}[b]{0.35\textwidth}
			\centering		
			\includegraphics[width=5.5cm]{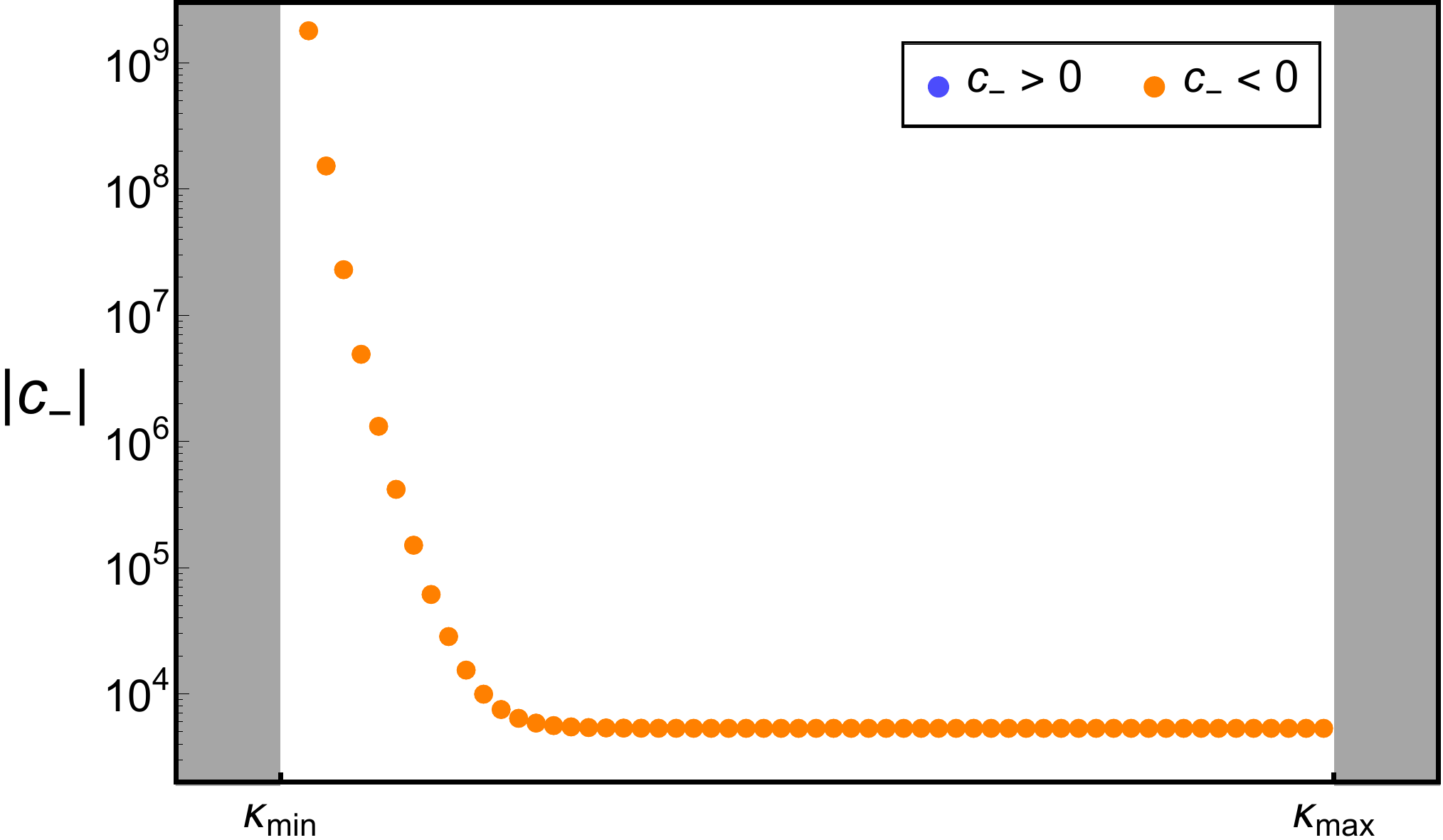}
			\caption{$(d,\ell,r) = (3,3,r_5)$}
		\end{subfigure}%
		\begin{subfigure}[b]{0.35\textwidth}
			\centering		
			\includegraphics[width=5.5cm]{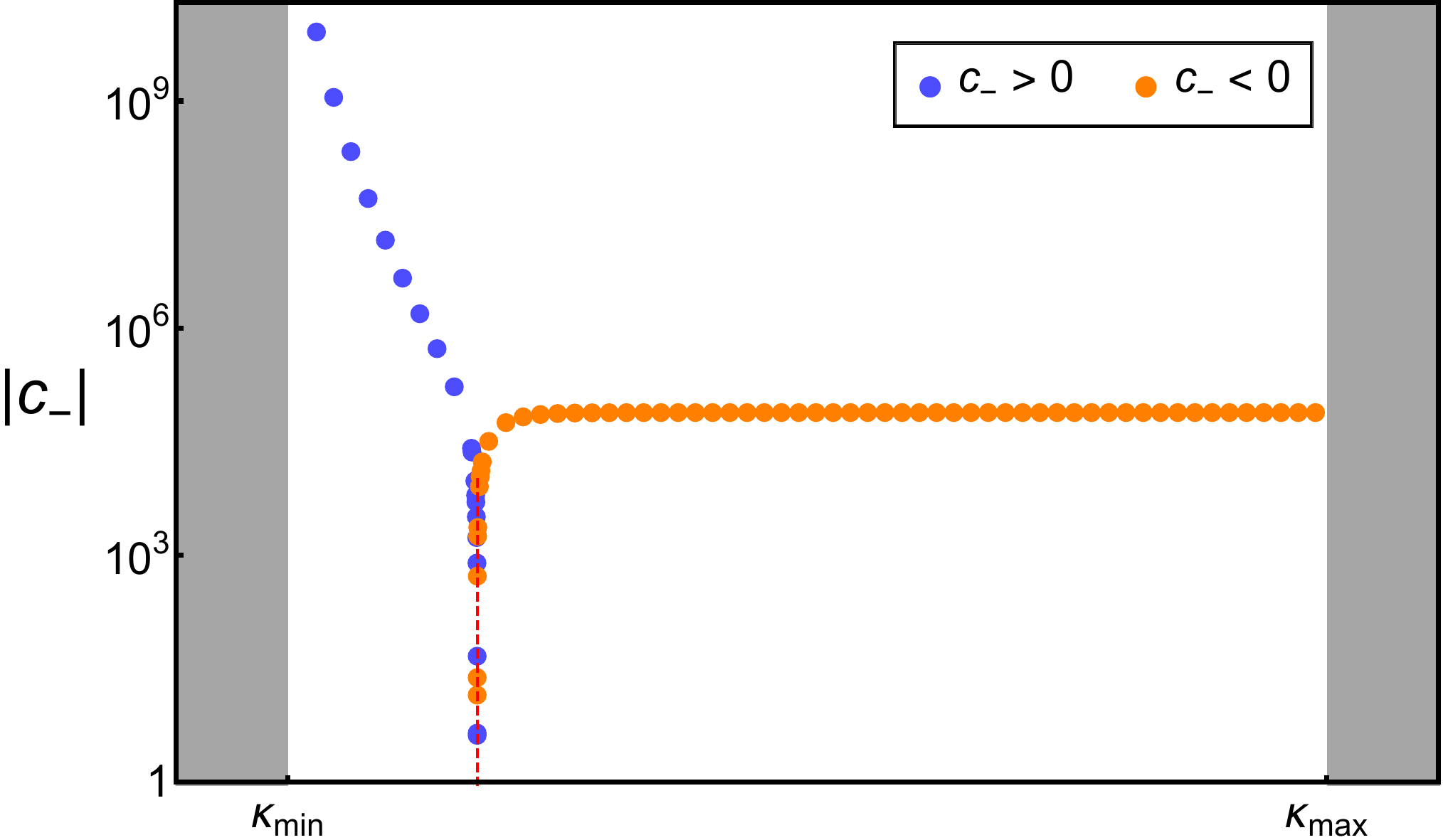}
			\caption{$(d,\ell,r) = (3,3,r_6)$}
		\end{subfigure}%
	
		\vspace{0.5cm}
		
		\begin{subfigure}[b]{0.35\textwidth}
			\centering		
			\includegraphics[width=5.5cm]{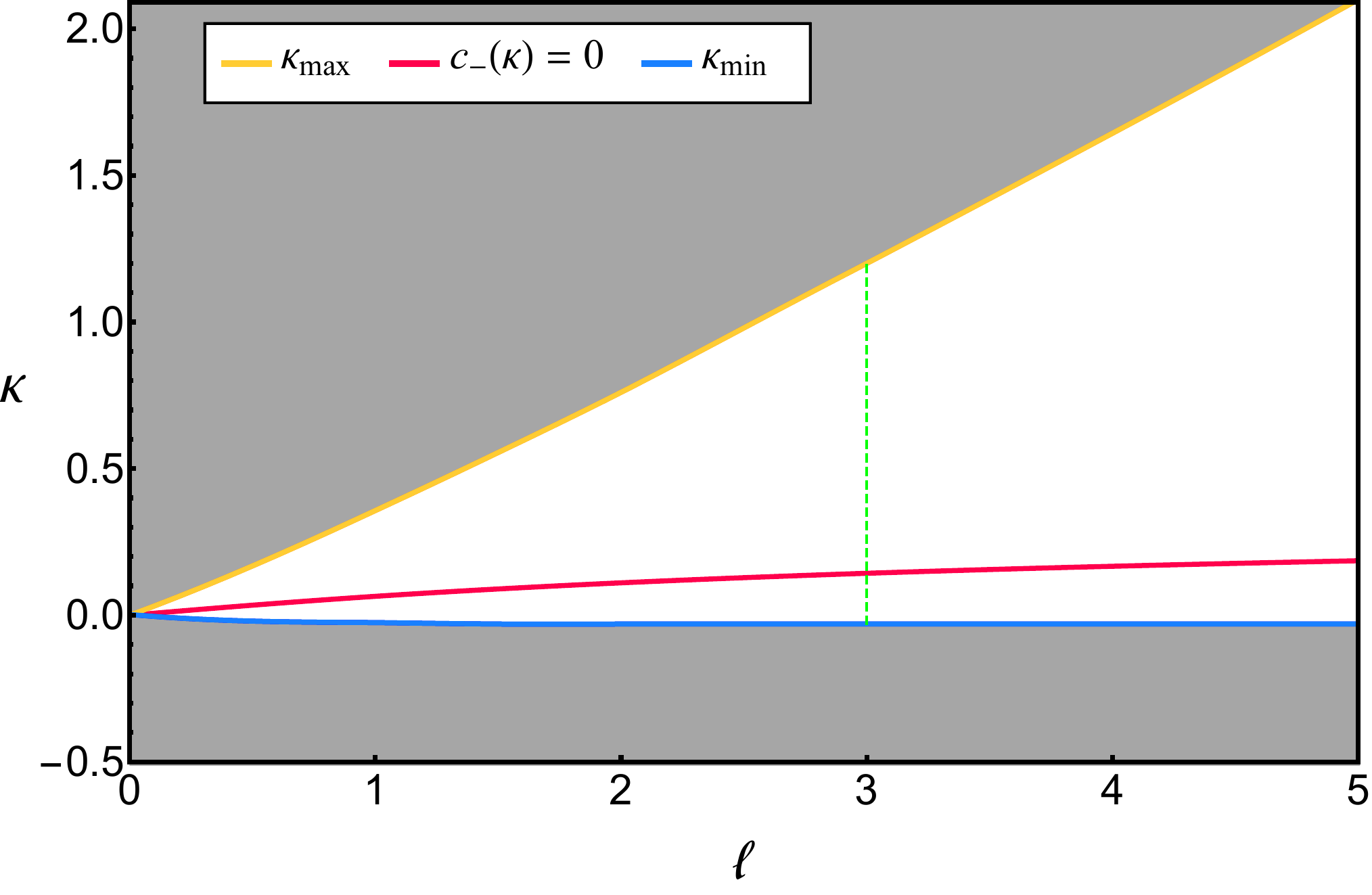}
			\caption{Diagram $(\ell,\kappa)$\\ $(d,r)=(3,r_2)$}
		\end{subfigure}%
		\begin{subfigure}[b]{0.35\textwidth}
			\centering		
			\includegraphics[width=5.5cm]{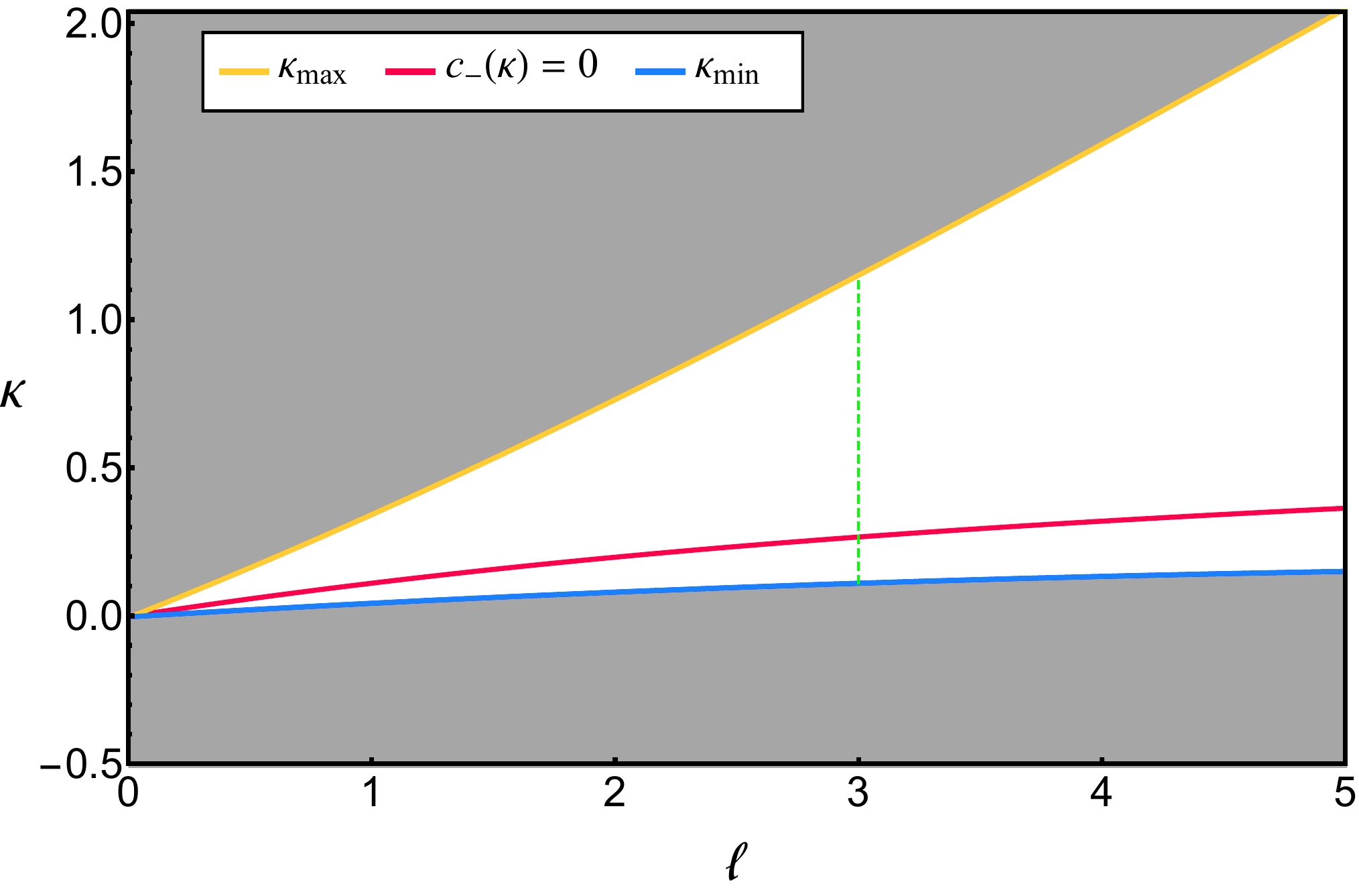}
			\caption{Diagram $(\ell,\kappa)$\\ $(d,r)=(3,r_4)$}
		\end{subfigure}%
		\begin{subfigure}[b]{0.35\textwidth}
			\centering		
			\includegraphics[width=5.5cm]{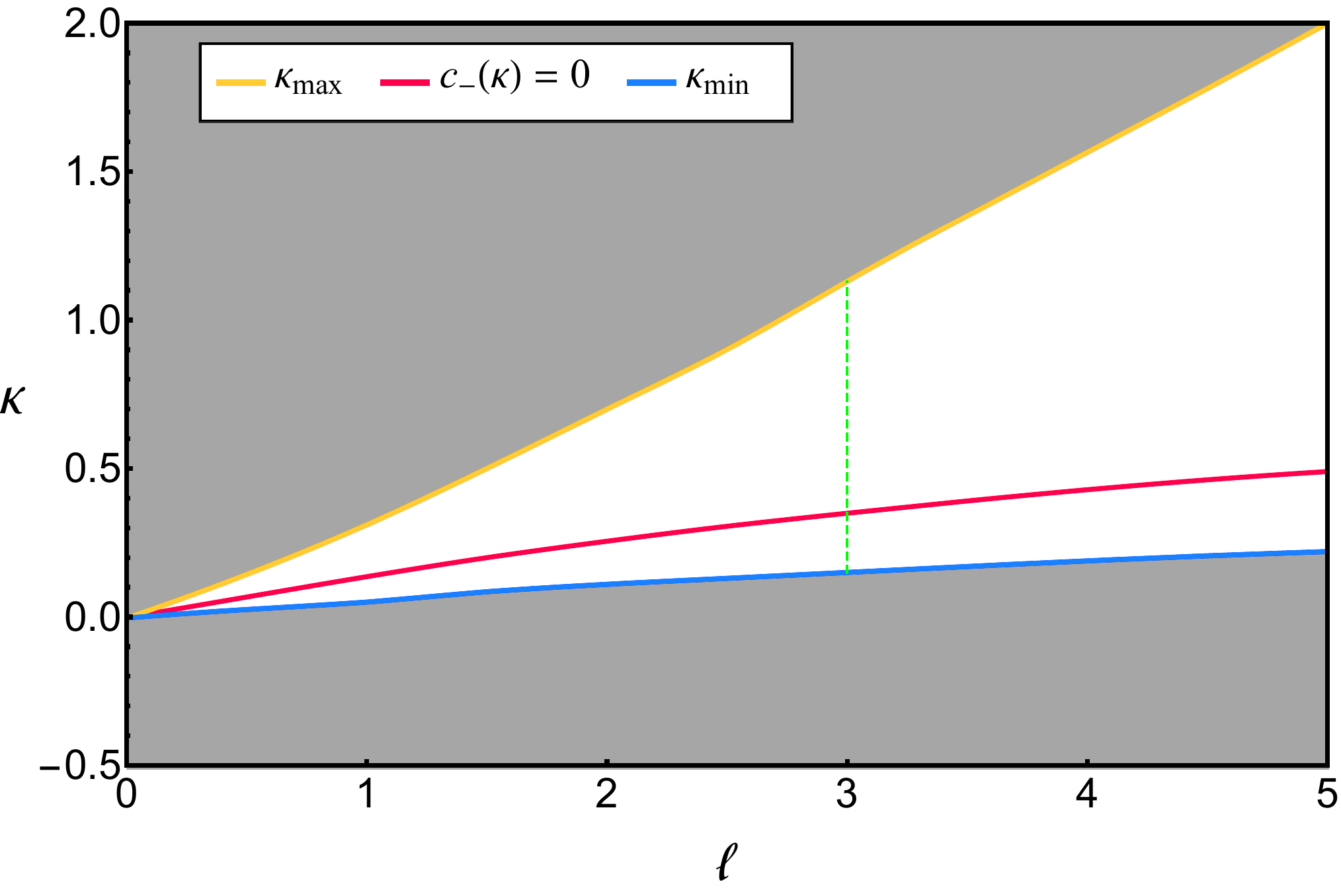}
			\caption{Diagram $(\ell,\kappa)$\\ $(d,r)=(3,r_6)$}
		\end{subfigure}%
		\caption{\small These plots show the search of the trajectories $P_2$-$P_4$ that have a smooth connection with the smooth trajectories $P_2$-$P_6$ ($r=r_n$); namely, the values of $\kappa_{\min}<\kappa<\kappa_{\max}$ such that $c_-(r_n,\kappa)=0$. We show the case of $d=3$ but the structure is representative of other dimensions. Plots from (a) to (f) show how the coefficient $c_-(\kappa)$ changes in $\kappa_{\min}<\kappa<\kappa_{\max}$; namely, along the family of trajectories $P_2$-$P_4$ that can be glued to $P_2$-$P_6$. Blue points indicate that $c_->0$, orange points that $c_-<0$, while the red dashed line marks the value of $\kappa$ where $c_-(\kappa)=0$. We observe that for $r_{1,3,5}$ (odd $n$) there is no $c_-=0$ while for $r_{2,4,6}$ (even $n$) we clearly find $c_-=0$. Plots (g) to (i) show the diagrams $(\kappa,\ell)$ with $d=3$ and $r=r_{2,4,6}$. These diagrams show the evolution of the zeros of $c_-$ (red lines) with $\ell$. We can observe that they lie in $\kappa_{\min}<\kappa<\kappa_{\max}$, even when this window is significantly shrunk. It provides strong confidence in the conclusion that the parity condition described in the main text is associated with $n$ and independent of $\ell$. Finally, the green dashed lines represent the location of plots (b), (d), (f) in these diagrams. \\}
		\label{fig:C_minus_search}
	\end{figure}

		The idea to construct SSs is very simple, we have to scan the space of parameters $(d,\ell,r,\kappa)$ in order to find trajectories $P_2$-$P_6$ with $c_+=0$ and trajectories $P_2$-$P_4$ with $c_-=0$. The most difficult part is at the technical level to reach the appropriate precision that this problem requires. To accomplish our goal we exploit the uniqueness of the curve $P_2$-$P_6$ (for given $(d,\ell,r)$) and the one-parameter family of curves $P_2$-$P_4$ given by $\kappa$. These two properties show that for given $(d,\ell)$ the key coefficients have the following dependence $c_{+}=c_{+}(r)$ and $c_{-}=c_{-}(r,\kappa)$. With it, the search is split into two parts, first we scan values of $r$ looking for $c_{+}=0$ and after that, using this value of $r$, we scan $\kappa$ to get $c_{-}=0$. Then, the full set of parameters $(d,\ell,r,\kappa)$ is determined. It shows that the problem of constructing SSs is reduced to two ``eigenvalue problems". The common strategy to solve this kind of problems is a shooting method where we choose a value of $r$ $(\kappa)$, impose the  condition $c_+=0$ ($c_-=0$) at $P_2$ and perform the integration from this point to see if the trajectory converges to $P_6$ ($P_4$). The iteration of this process for different values of $r$ $(\kappa)$ may determine the desired trajectory $P_2$-$P_6$ ($P_2$-$P_4$) with $c_+=0$ ($c_-=0$). However, note that this strategy moves through the space of inadmissible curves (they begin at $P_2$ but do not end at $P_6$ ($P_4$)). For this reason, the only valuable information that we obtain are the values $(r,\kappa)$ such that $c_+=c_-=0$. In order to get a good understanding of this problem we do not follow the process described above. Instead, our search moves through the space of trajectories $P_2$-$P_6$ ($P_2$-$P_4$). We construct these curves for different values of $r$ $(\kappa)$ to study the structure of $c_+(r)$ $(c_-(\kappa))$ restricted to these families of trajectories. It allows us to perform a more efficient search and to learn important lessons about the structure of this problem. Technical details can be found in appendix~\ref{sec:Appendix_Numerical_Methods}. The main observations that we obtain for $d\geq2$, $\ell>0$ and $\nu <8$ are ($[\nu]$ denotes the integer part of $\nu$)
		\begin{enumerate}
			\item $c_{+}(r)$ has zeros for discrete values of $r$ that we denote by $\{r_n\}$ such that they satisfy $r_1<r_2<r_3<...$ (see fig.~\ref{fig:Search_rk_dimensions}). We have found the first values of this sequence for $d\geq2$, $\ell>0$ and $\nu<8$. Some of them are provided in table~\ref{table:r_1_full}.
			
			\item $\{r_n\}$ are continuous functions of $\ell$ that do not intersect each other (see fig.~\ref{fig:Search_rk_B}); namely, we have that $r_1(\ell)<r_{2}(\ell)<r_3(\ell)<...$
			
			\item $c_-(r_n,\kappa)$ has a parity condition in $n$. Only for $n$ even ($n$ odd for $d=2$) there is a trajectory $P_2$-$P_4$ with $c_-(r_n,\kappa)=0$ (see fig.~\ref{fig:C_minus_search}). Furthermore, we find that this property is independent\footnote{That the parity condition is associated with $n$ and independent of $\ell$ is a strange property because, from a naive point of view, trajectories $P_2$-$P_6$ and $P_2$-$P_4$ are independent. For this reason we expect that there is a deeper explanation for this property. An alternative option is that this is an artifact of the looking at a finite region in the space of parameters. For example, it is possible that the parity condition is no longer associated with $n$ for high values of this number or that for $\ell\ll 1$ the parity condition is not satisfied. However, these scenarios are not convincing due to the robust structure that we observe in fig.~\ref{fig:C_minus_search}; note that even for small $\ell$, when the range of $\kappa$ is significantly shrunk, the parity condition is preserved.} of $\ell$. Note that this parity condition is very important because only for $n$ even (odd in $d=2$) there are SSs. 
			
			\item For a given $\ell$ we usually find a single $r_n$ between consecutive integer values of $[\nu]$. One may be tempted to associate $n$ with $[\nu]$; however, the description in terms of $[\nu]$ is not simple. This is because each curve $r_n(\ell)$ has special properties related to $n$ but some of them can cross lines $\nu=[\nu]$ and transit from a value of $[\nu]$ to another. It can be seen in fig.~\ref{fig:Search_rk_B_d5}.
			
			\item The distance $|r_{n+1}-r_n|$ decreases when $n$ grows. This is consequence that the interval of $r$ where $[\nu]$ is a particular integer shrinks when $[\nu]$ increases, going to zero for $[\nu]\to\infty$; see fig.~\ref{fig:Search_rk_dimensions},\ref{fig:Search_rk_B}.
			
			\item In $d\geq3$, $c_{+}(r)$ goes to zero when $r$ approaches $1$ (see  fig.~\ref{fig:Search_rk_dimensions}). An inspection of this region shows that there is no zero for $r$ close to $1$; however, our numerical methods cannot rule out the existence of a zero for $0<r-1\ll1$. Note that it would be $r_0$ and given the parity condition it should be associated with a SS, which is the perfect candidate to be stable. In $d=2$ we do not report results in the region $r\sim1$ for technical reasons explained in fig.~\ref{fig:Search_rk_dimensions}.
		\end{enumerate}
		This structure, that we find in the regime of low regularity of NSSs, $[\nu]\leq7$ ($r\gtrsim1$), is in excellent agreement with the results obtained in \cite{MerleEuler} in the regime of high regularity of NSSs, $[\nu]\gg 1$ ($r\to \reye$). There, the authors also found that given $(d,\ell)$ SSs exist for a specific parity\footnote{In \cite{MerleEuler} the authors discuss the existence of SSs based on the parity of $[\nu]$ which must be even or odd in relation to the sign of a particular function of $(d,\ell)$. However, the formulation of the parity condition in terms of $[\nu]$ is not appropriate for the full range $1<r<\reye$ because it may change; note that some curves $r_{n}(\ell)$ can cross lines $\nu=[\nu]$ and the parity in $[\nu]$ transits. A visual example can be found in fig.~\ref{fig:Search_rk_B_d5}.} of $[\nu]$ and discrete values of $r$ that accumulate at $\reye$. It provides strong confidence that the structure that we report is the general structure of the problem ($1<r<\reye$).

		\subsection{Linear Perturbations}
		\label{sec:Linear_Perturbations}
		
		In order to explore the vicinity of blow-up profiles we study the spectrum of linear perturbations. It is done from the blow-up profiles, denoted by $(\hat\rho_0,\hat u_0)$, following ansatz
		\begin{align}
			&\hat\rho(\tau, Z) = \hat\rho_{0}(Z) + \epsilon\  \alpha(Z) e^{\LL \tau} + \mathcal{O}(\epsilon^2), \label{eq:ansatz_linear_modes_alpha}\\
			&\hat u(\tau, Z) = \hat u_{0}(Z) + \epsilon\ \beta(Z)e^{\LL \tau}  + \mathcal{O}(\epsilon^2), \label{eq:ansatz_linear_modes_beta}
		\end{align}
		with $|\epsilon| \ll 1$, $\LL \in \mathbb{R}$. The most important element is the sign of the exponent $\LL$, it provides the criterion to determine whether a linear mode is stable ($\LL<0$) or unstable ($\LL>0$). Plugging this ansatz into the time-dependent equations for $(\hat\rho,\hat{u})$, (\ref{eq:Rho_dtau})-(\ref{eq:U_dtau}), at linear order in $\epsilon$ we get the equations for linear perturbations ($\ ':=\partial_Z$)
		\begin{align}
			&\left(\LL + \ell(r-1) + \hat{u}_0'+\frac{d-1}{Z}\hat{u}_0\right)\alpha + (Z+\hat{u}_0)\alpha' + \hat{\rho}_0 \beta' + \left(\hat{\rho}_0' + \frac{d-1}{Z}\hat{\rho}_0\right)\beta=0, \label{eq:alpha_pert}\\
			&\left(\LL + r - 1 + \hat{u}_0'\right)\beta+(Z+\hat{u}_0)\beta' + (\gamma-1)\hat{\rho}_0^{\gamma-2}\alpha'+(\gamma-1)(\gamma-2)\hat{\rho}_0^{\gamma-3}\hat{\rho}_0'\alpha=0. \label{eq:beta_pert}
		\end{align}
		Linear modes $(\alpha,\beta)$ also extend beyond the acoustic cone to connect the origin with the infinity, and again there is no a trivial way to smoothly cross this surface. For this reason $Z_2$ is also a {\em regular singular point} for these equations (appendix~\ref{sec:Appendix_Regularity_Linear_Modes_Z2}). The construction of smooth perturbations also consists of the exploration of the space of parameters to find trajectories from the origin and infinity to the acoustic cone that have a smooth connection at this point.
		 
\begin{figure}[t!]
	 	\centering	
	 	\begin{subfigure}[b]{0.5\textwidth}
	 		\centering
	 		\includegraphics[width=7.5cm]{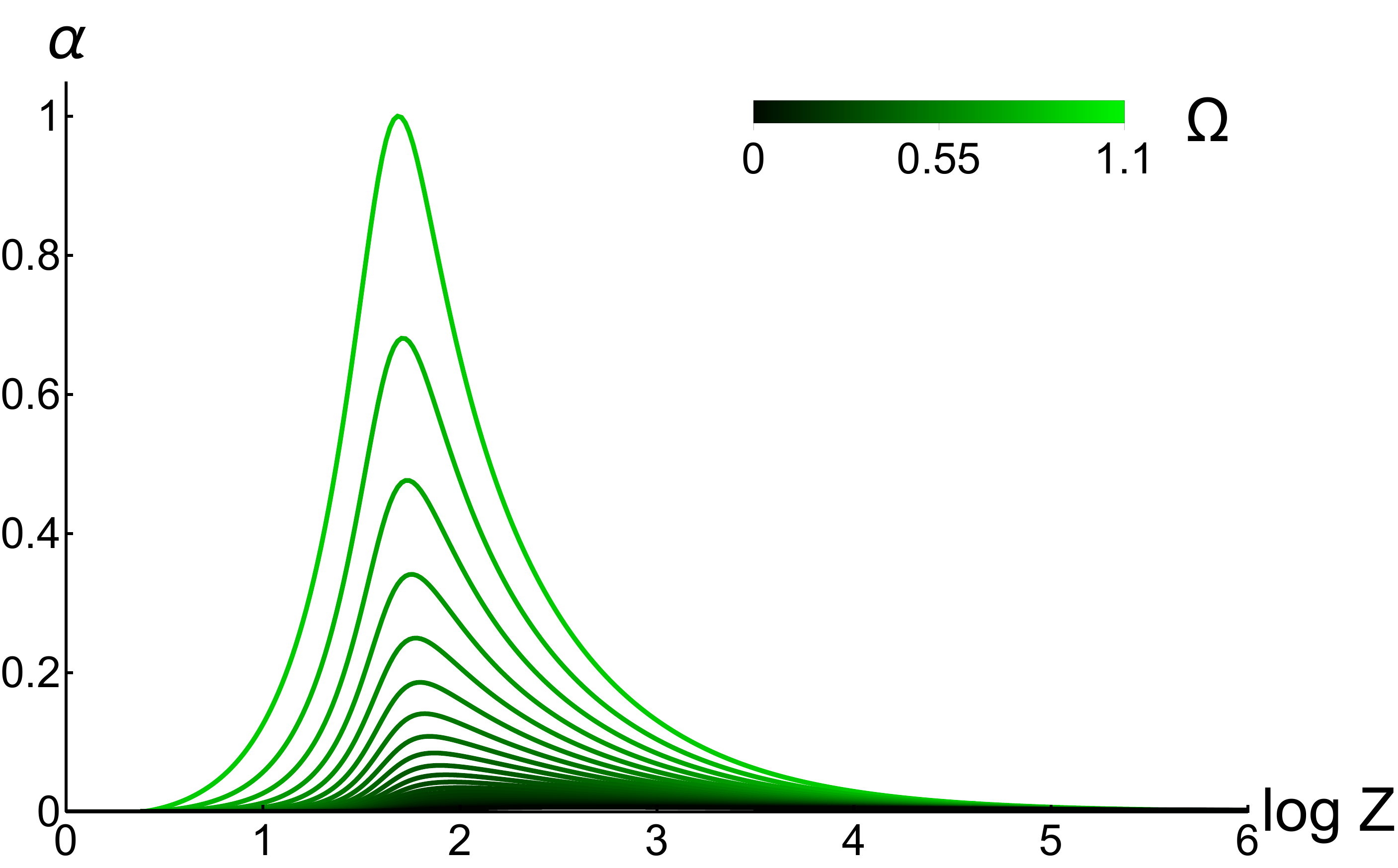}
	 		\caption{0-modes}
	 		\label{fig:linear_modes_0_alpha}
	 	\end{subfigure}%
	 	\begin{subfigure}[b]{0.5\textwidth}
	 		\centering		\hspace{2cm}\includegraphics[width=7.5cm]{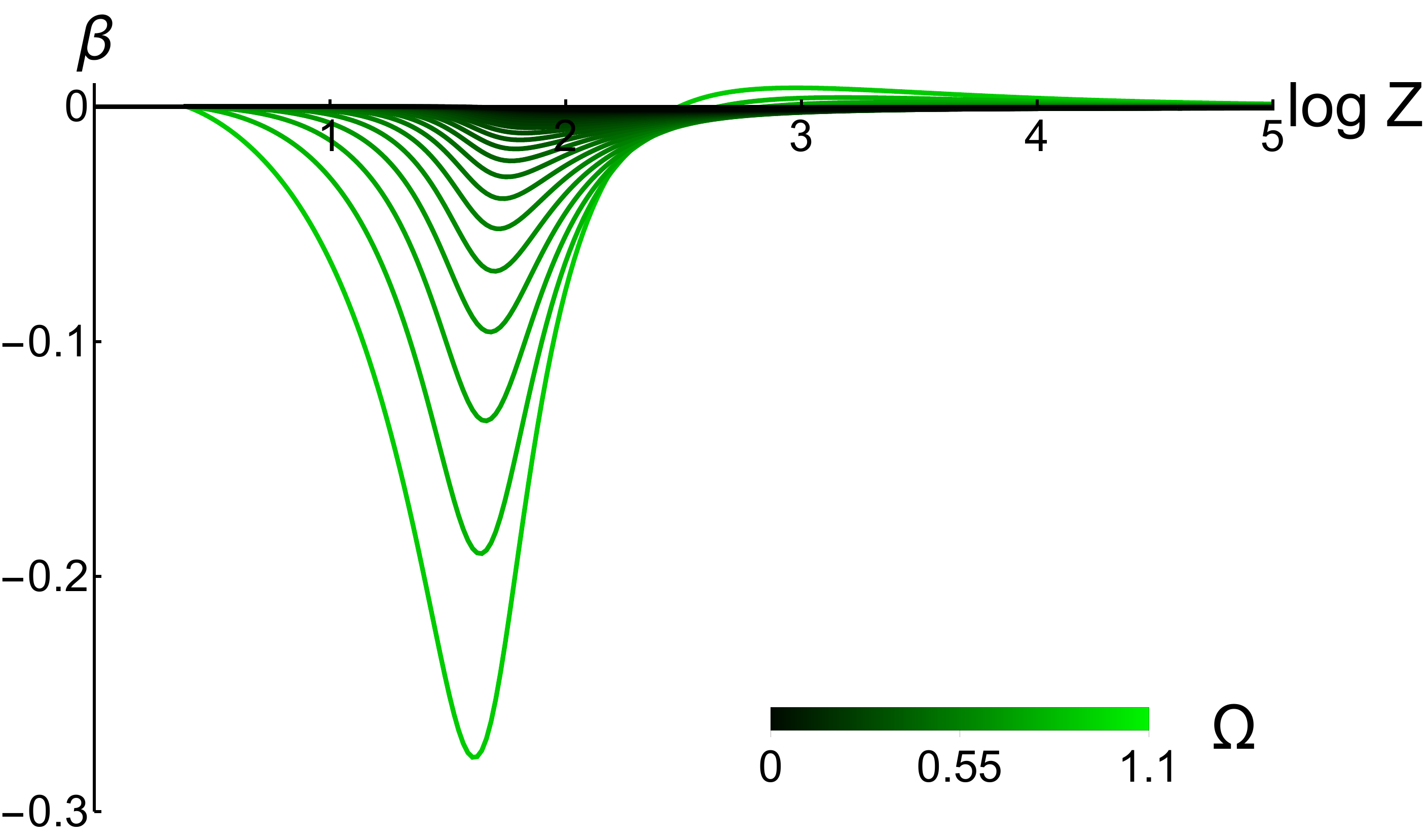}
	 		\caption{0-modes}
	 		\label{fig:linear_modes_0_beta}
	 	\end{subfigure}%
	 	
	 	\vspace{0cm}
	 	\begin{subfigure}[b]{0.5\textwidth}
	 		\centering
	 		\includegraphics[width=7.5cm]{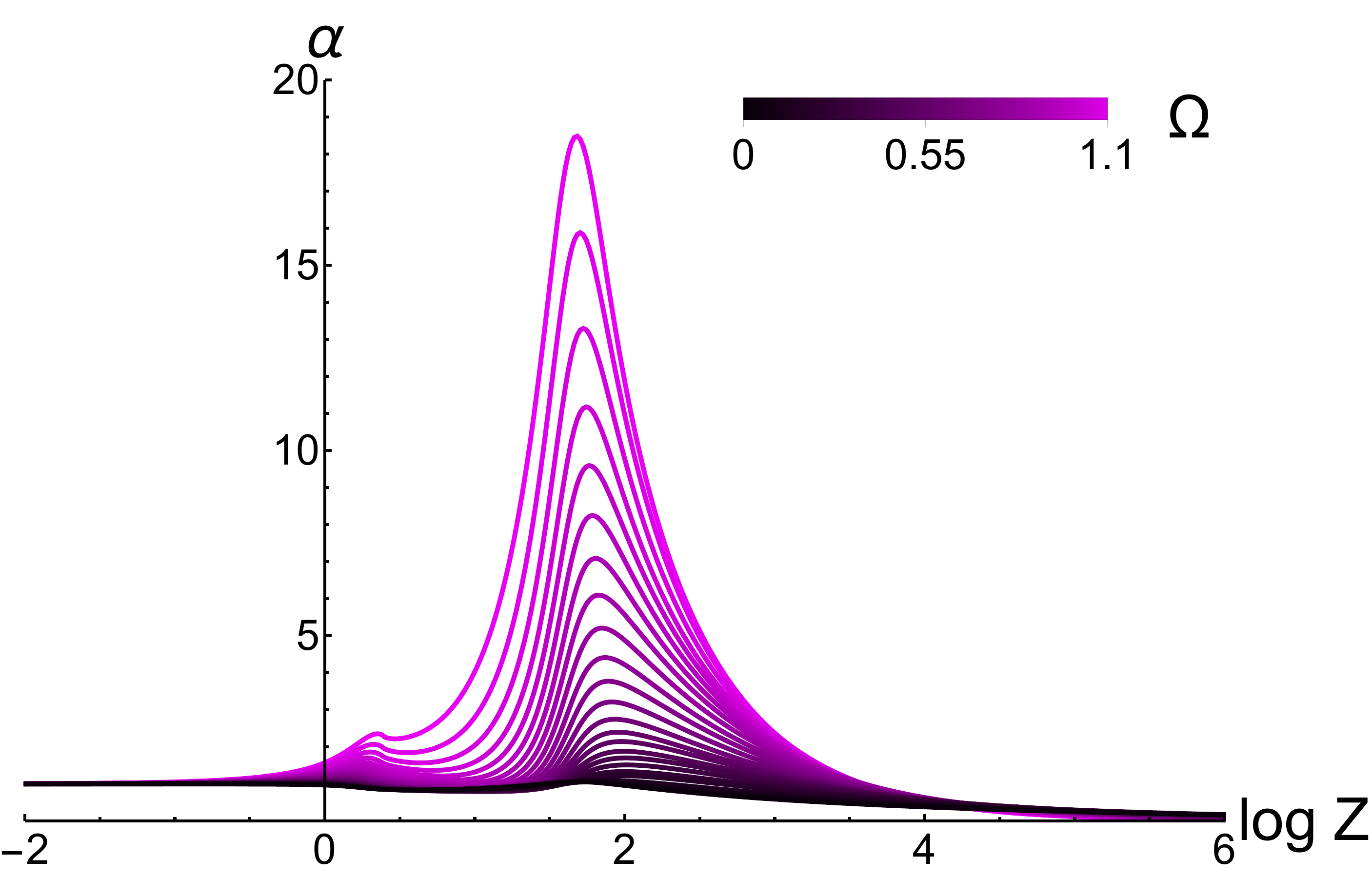}
	 		\caption{1-modes}
	 		\label{fig:linear_modes_1_alpha}
	 	\end{subfigure}%
	 	\begin{subfigure}[b]{0.5\textwidth}
	 		\centering		\hspace{2cm}\includegraphics[width=7.5cm]{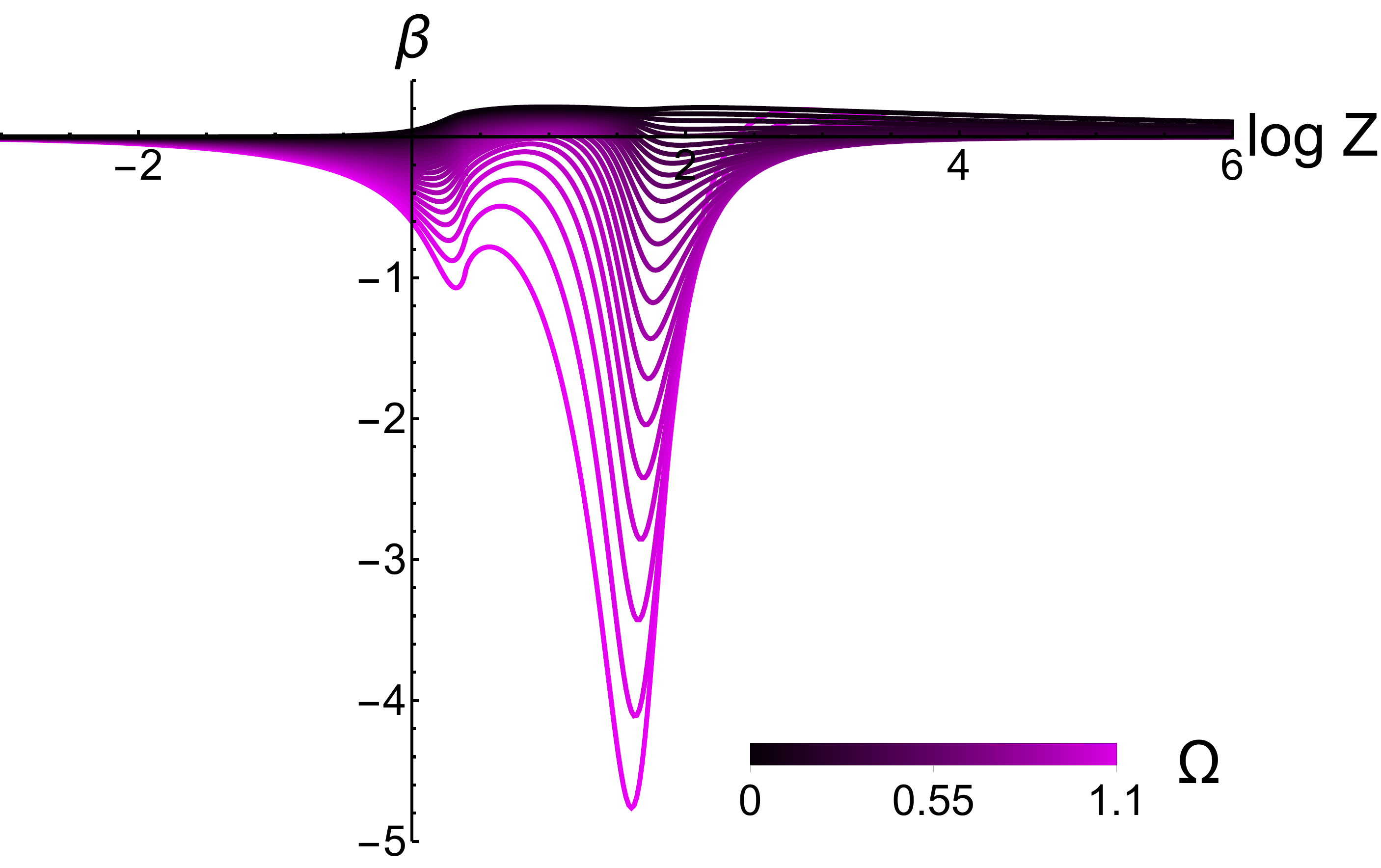}
	 		\caption{1-modes}
	 		\label{fig:linear_modes_1_beta}
	 	\end{subfigure}%
	 	\caption{\small Linear modes associated with the set of parameters $(d,\ell,r,\kappa) = (3,2,1.2,0.6)$. 0-modes ($\alpha(Z_2)=\beta(Z_2) = 0$) in the upper row and 1-modes ($|\alpha(Z_2)|+|\beta(Z_2)| \neq 0$) in the lower row. The color labels the value of $\LL$ following the legends.}
	 	\label{fig:linear_modes}
	 \end{figure}

		 Our analysis of the set of equations (\ref{eq:alpha_pert})-(\ref{eq:beta_pert}) given in appendix~\ref{sec:Appendix_Regularity_Linear_Modes_Z2} shows that there is a continuous range of values $\LL_{\min}<\LL< \LL_{max}$ where linear modes have two linearly independent admissible solutions (they vanish at infinity and have at least a continuous derivative on the cone). We classify them into two types in relation to their values at the origin and $Z_2$:
		\begin{eqnarray}
		\text{0-modes:}  &\quad \alpha(0)=\beta(0)=\alpha(Z_2)=\beta(Z_2) = 0 \label{eq:0-mode_definition}\\
		\text{1-modes:} & \quad |\alpha(0)|+|\beta(0)| \neq 0 \quad  |\alpha(Z_2)|+|\beta(Z_2)| \neq 0 \label{eq:1-mode_definition}
		\end{eqnarray}
		Fig.~\ref{fig:linear_modes} provides a visual representation of generic 0-modes and 1-modes. They are constructed using similar techniques to obtain self-similar profiles (appendix~\ref{sec:Appendix_Numerical_Methods}).

	 
	\subsubsection{0-modes}
\begin{figure}[h!]
	\centering	
	\vspace{1cm}
	\begin{subfigure}[b]{0.5\textwidth}
		\centering
		\includegraphics[width=7.5cm]{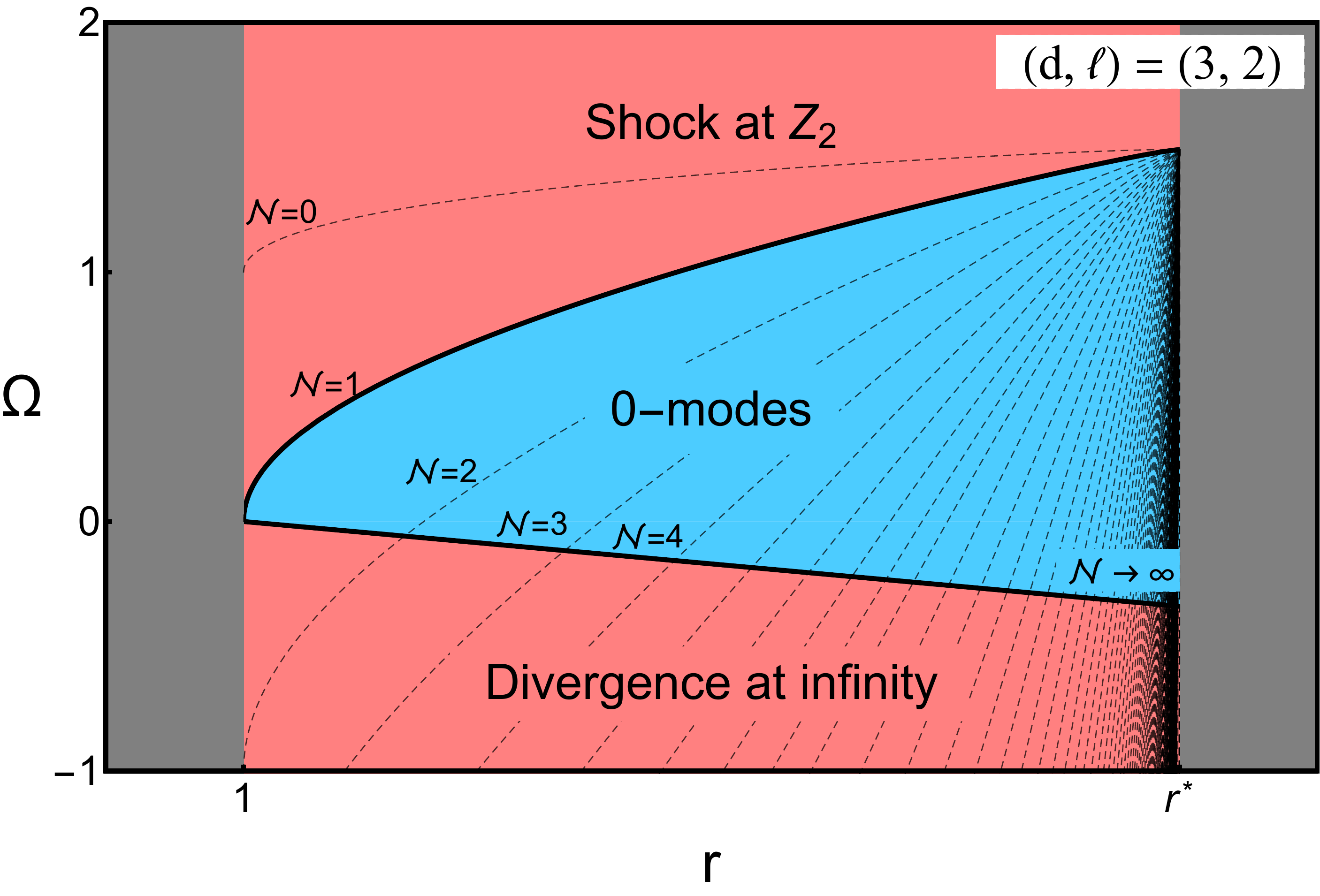}
		\label{fig:Diagram_range_lambda_d3_L2}
	\end{subfigure}%
	\begin{subfigure}[b]{0.5\textwidth}
		\centering		\hspace{2cm}\includegraphics[width=7.5cm]{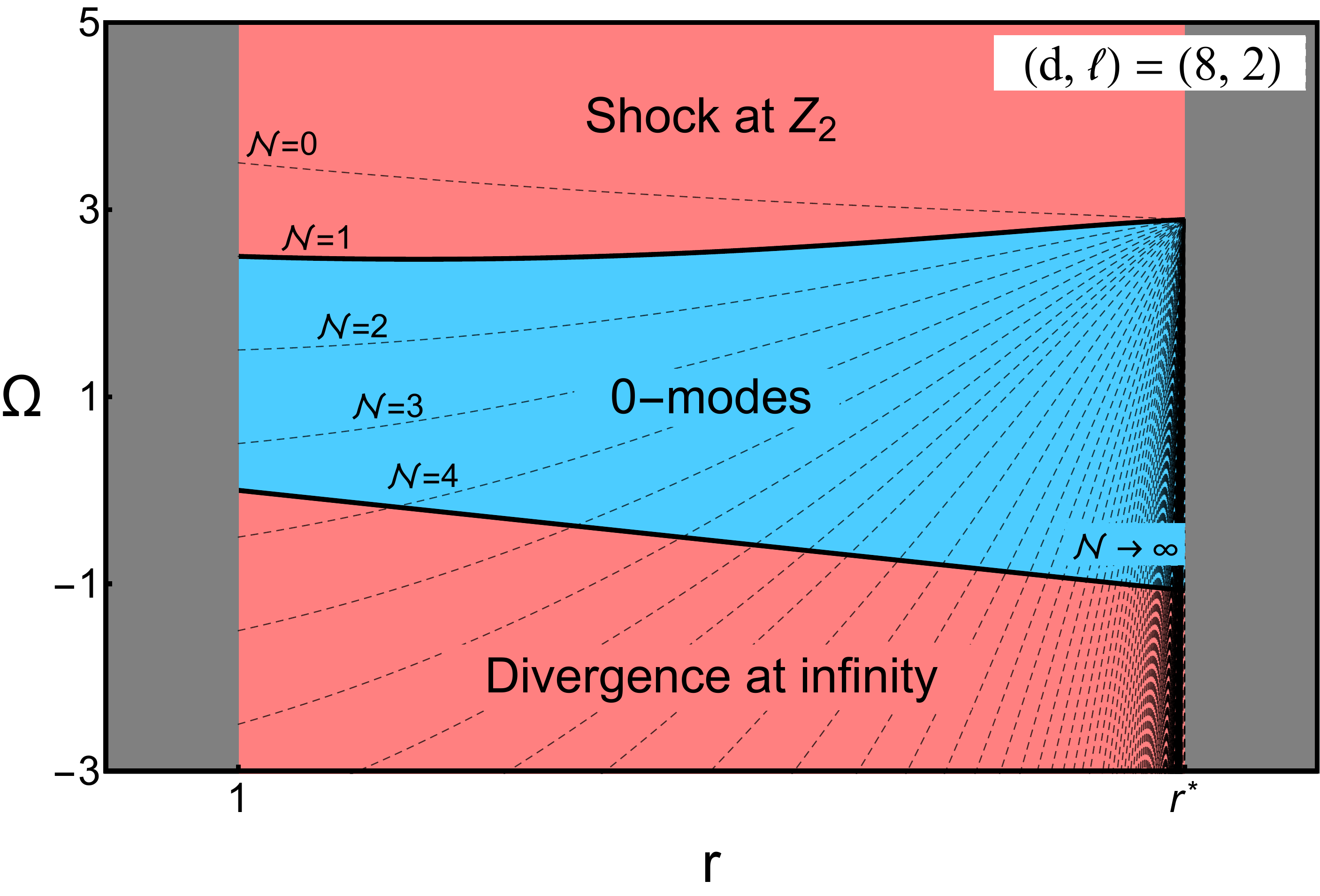}
		\label{fig:Diagram_range_lambda_d8_L2}
	\end{subfigure}%
	\caption{\small Diagrams representing the admissible values of $(r,\LL)$ (blue area) for which 0-modes are at least $\mathcal{C}^1(Z_2)$ and vanish at infinity. Black dashed lines represent curves $\LL(r)$ where $\mathcal{N}$ is integer. This regularity goes to infinity when $r\to\reye$ providing a neighborhood close to this value of admissible 0-modes with arbitrarily high regularity, but finite.}
	\label{fig:Diagram_range_lambda_0-modes_full}
\end{figure}

		 0-modes are engineered by gluing the trivial solution $\alpha=\beta=0$ on the interior of the acoustic cone, $Z\in[0,Z_2]$, to solutions that vanish at $Z_2$ on the exterior of the cone, $Z\in[Z_2,\infty)$. This process results in linear modes in $Z\in[0,\infty)$ with the following properties (see appendix~\ref{sec:Appendix_Regularity_Linear_Modes_Z2} for further details):
		 \begin{itemize}
		 	\item 0-modes have regularity $\mathcal{C}^{\mathcal{N}}(Z_2)$ that depends on $\LL$ in the form ($\nu$ given in (\ref{eq:expansion_w_in_powers_sigma_Z2}))
		 	\beq
		 		\mathcal{N}(\LL) = \nu + \frac{2}{(\ell+1)(r-1)-(d+1)} \left(\nu+1\right) \LL.
		 		\label{eq:regularity_lambda_0_modes}
		 	\eeq
		 	\item $\LL$ has an upper bound $\LL_{max}$ coming from the condition that $\mathcal{N}(\LL_{\max})=1$. $\LL_{\max}>0$ in our range of parameters (\ref{eq:parameters}), meaning that blow-up profiles always have unstable directions triggered by some 0-modes. See fig.~\ref{fig:Diagram_range_lambda_0-modes_full} for a visual representation.
		 	\item $\LL$ has a lower bound $\LL_{min} = (1-r)\min(\ell,1) <0$. This bound comes from the condition that linear modes vanish at infinity. See fig.~\ref{fig:Diagram_range_lambda_0-modes_full} for a visual representation.
		 	\item 0-modes with $\LL = 0$ have regularity $\nu$. Roughly speaking these modes are the difference between two self-similar solutions that are very close in $\kappa$, defined in (\ref{eq:expansion_ws_boundary}); namely, they move one of these solutions along this family of self-similar profiles. They arise as the first non-vanishing contribution in the expansion of our self-similar solutions in powers of $\kappa$. Using the uniqueness of these solutions in the interval $[0,Z_2)$ and continuity at $Z_2$, the resulting mode vanishes in $[0,Z_2]$; namely, this is a 0-mode.  
		 	\item In our range of parameters, (\ref{eq:parameters}), unstable 0-modes ($\LL>0$) have less regularity than the self-similar solution ($\mathcal{N}<\nu$). 
		 \end{itemize}

	\subsubsection{1-modes and Smooth Linear Perturbations}
	
\begin{figure}[h!]
		\centering	
		\vspace{0.5cm}
		\begin{subfigure}[b]{0.5\textwidth}
			\centering		\hspace{2cm}\includegraphics[width=8.cm]{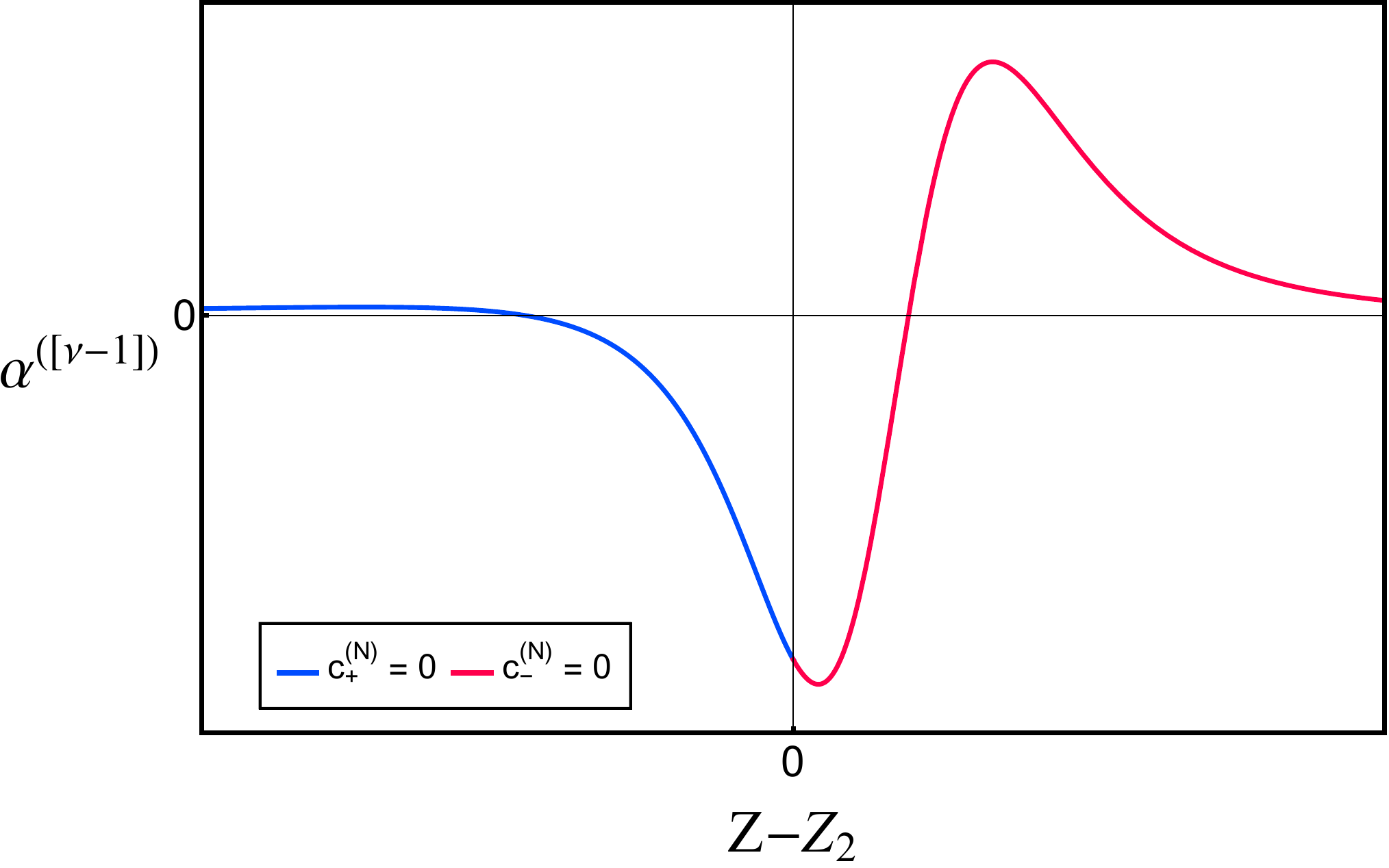}
			\caption{SS}
			\label{fig:Search_Omega_k_1-modes_D}
		\end{subfigure}%
		\begin{subfigure}[b]{0.5\textwidth}
			\centering		\hspace{2cm}\includegraphics[width=8.cm]{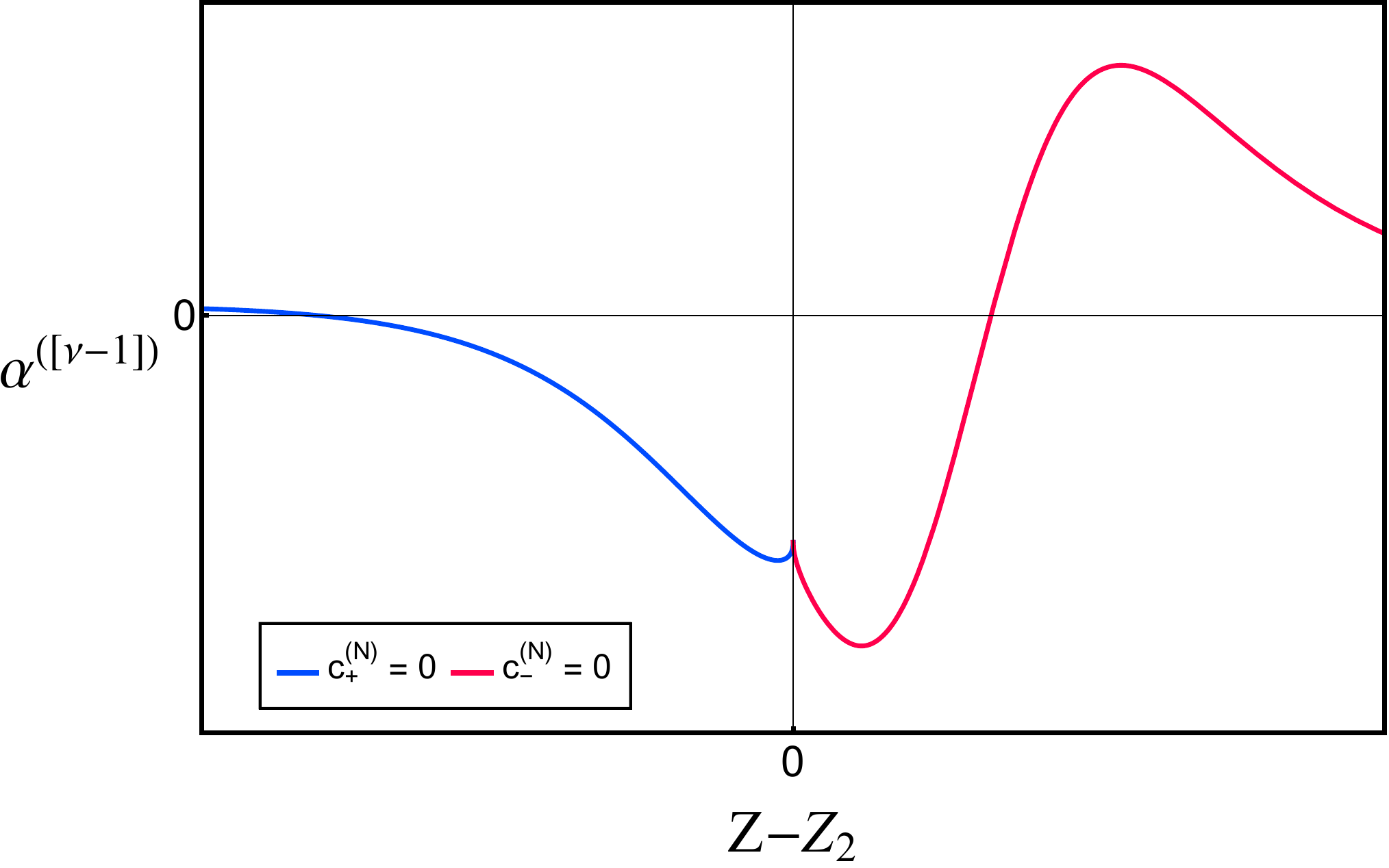}
			\caption{NSS}
			\label{fig:Search_Omega_k_1-modes_E}
		\end{subfigure}%
		\caption{\small  $[\nu-1]$-derivative of $\alpha(Z)$ with $c_{\pm}^{(\mathcal{N})}=0$ and $\mathcal{N}<\nu-1$. We see that for SSs (NSSs) $c_{\pm}^{(\nu)}= 0$ ($c_{\pm}^{(\nu)} \neq 0$) and therefore this 1-mode is smooth (non-smooth).}
		\label{fig:Search_Omega_k_1-modes_nu-der_smooth_non-smooth_profiles}
	\end{figure}

\begin{figure}[h!]
		\vspace{0.5cm}
		\centering	
		\begin{subfigure}[b]{0.5\textwidth}
			\centering
			\includegraphics[width=8.2cm]{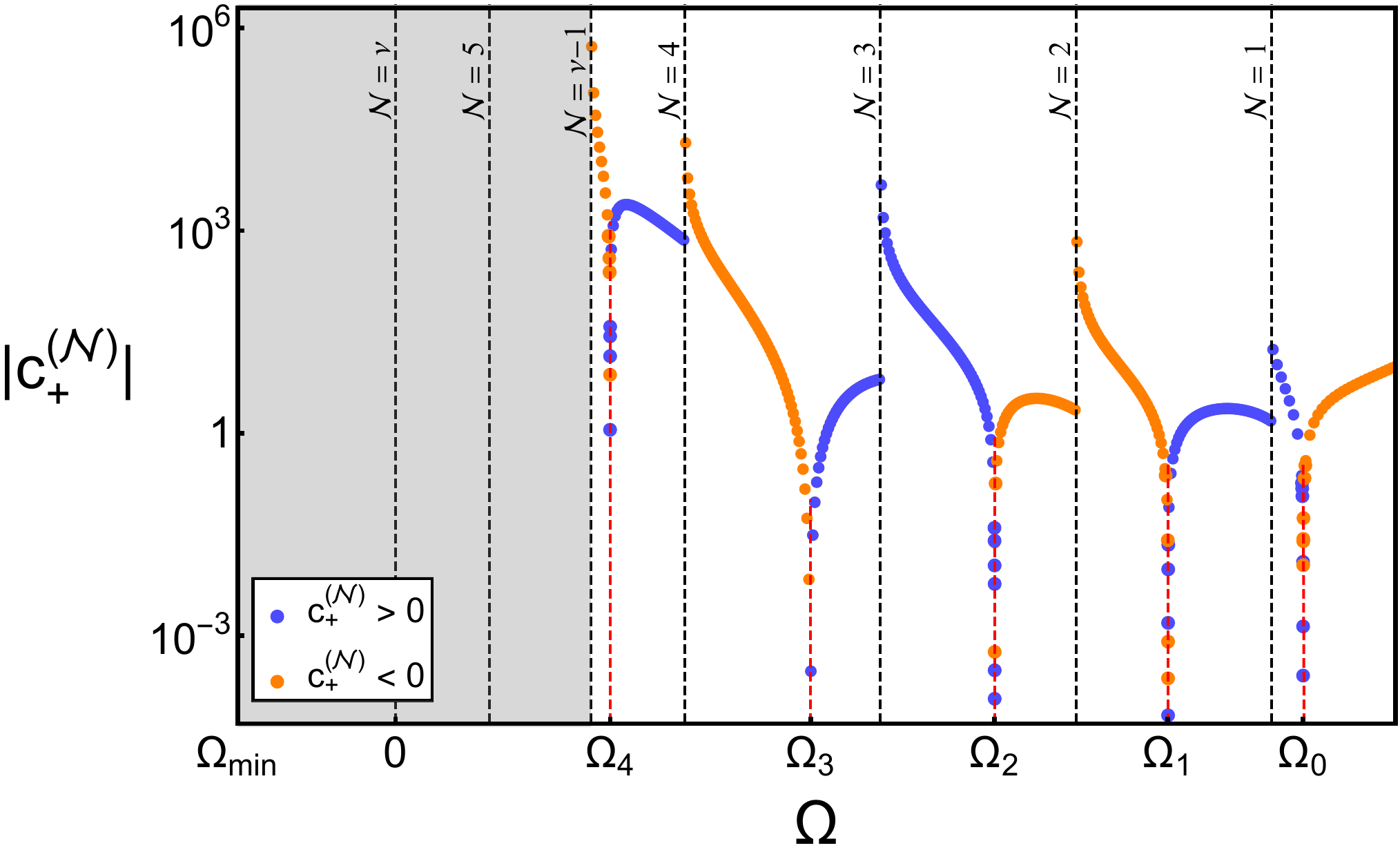}
			\caption{Function $|c_+^{(\mathcal{N})}(\LL)|$}
			\label{fig:Search_Omega_k_1-modes_A}
		\end{subfigure}%
		\begin{subfigure}[b]{0.5\textwidth}
			\centering
			\includegraphics[width=8.3cm]{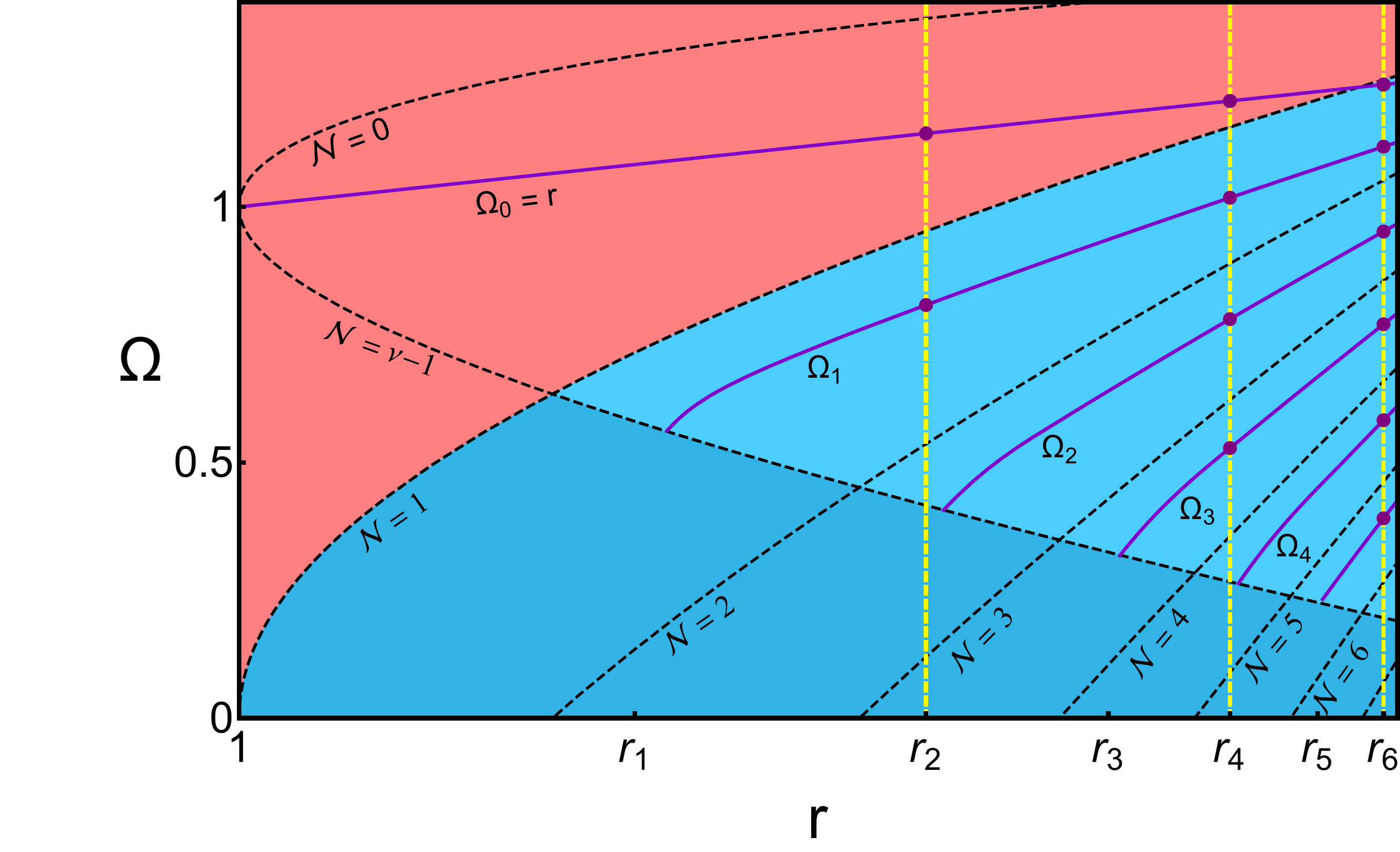}
			\caption{Diagram $(r,\LL)$ of 1-modes in $d=3$}
			\label{fig:Diagram_Smooth_1-modes}
		\end{subfigure}%
								
		\vspace{0.5cm}

		\begin{subfigure}[b]{0.5\textwidth}
		\centering
		\includegraphics[width=6.5cm]{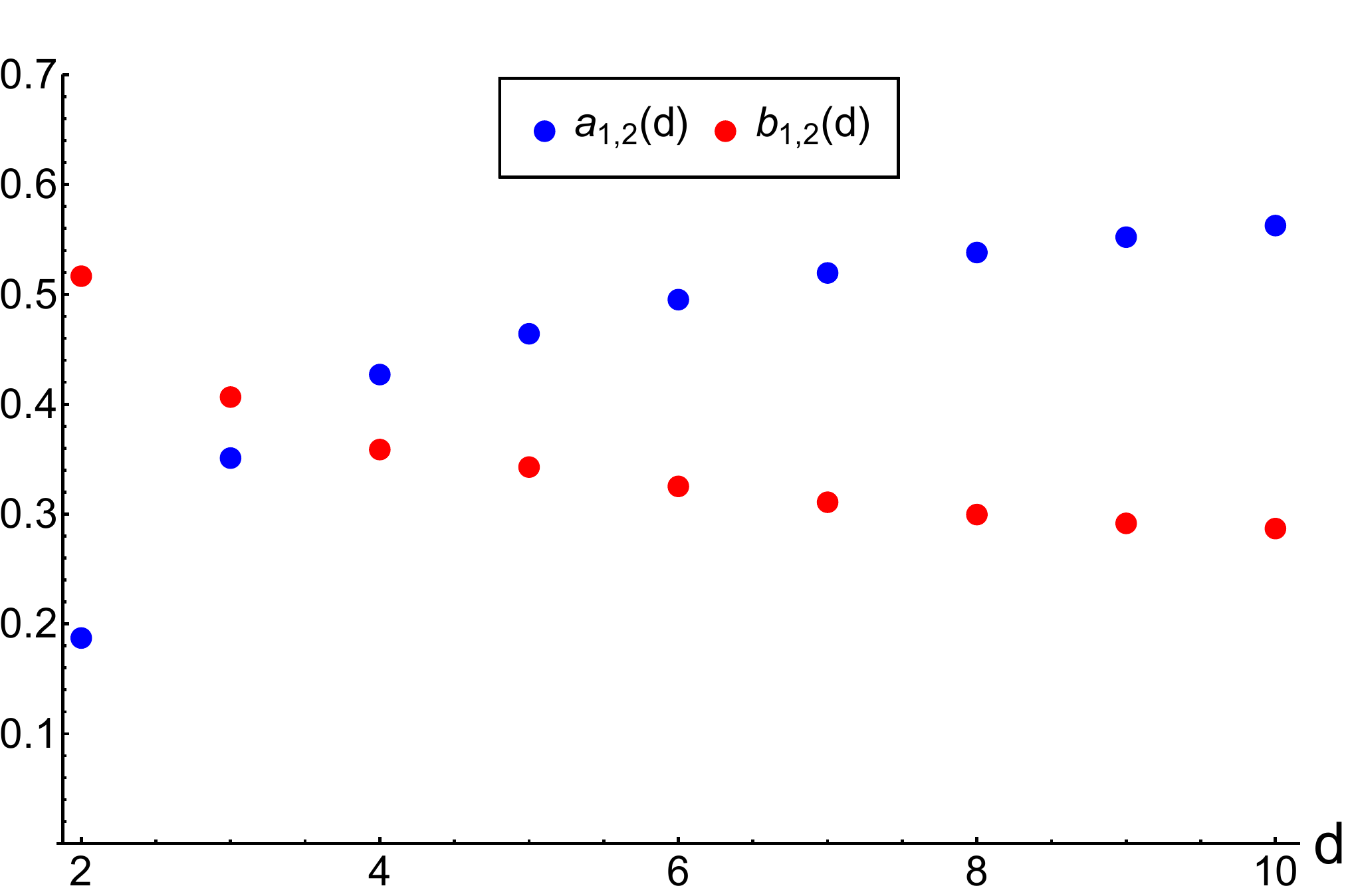}
		\caption{Coefficients $\LL_{1}\simeq a_{1,2} r_2(\ell)+b_{1,2}$}
		\label{fig:Search_Omega_k_AB_coefficients_Omega-1}
		\end{subfigure}%
	\begin{subfigure}[b]{0.5\textwidth}
		\centering
		\includegraphics[width=6.5cm]{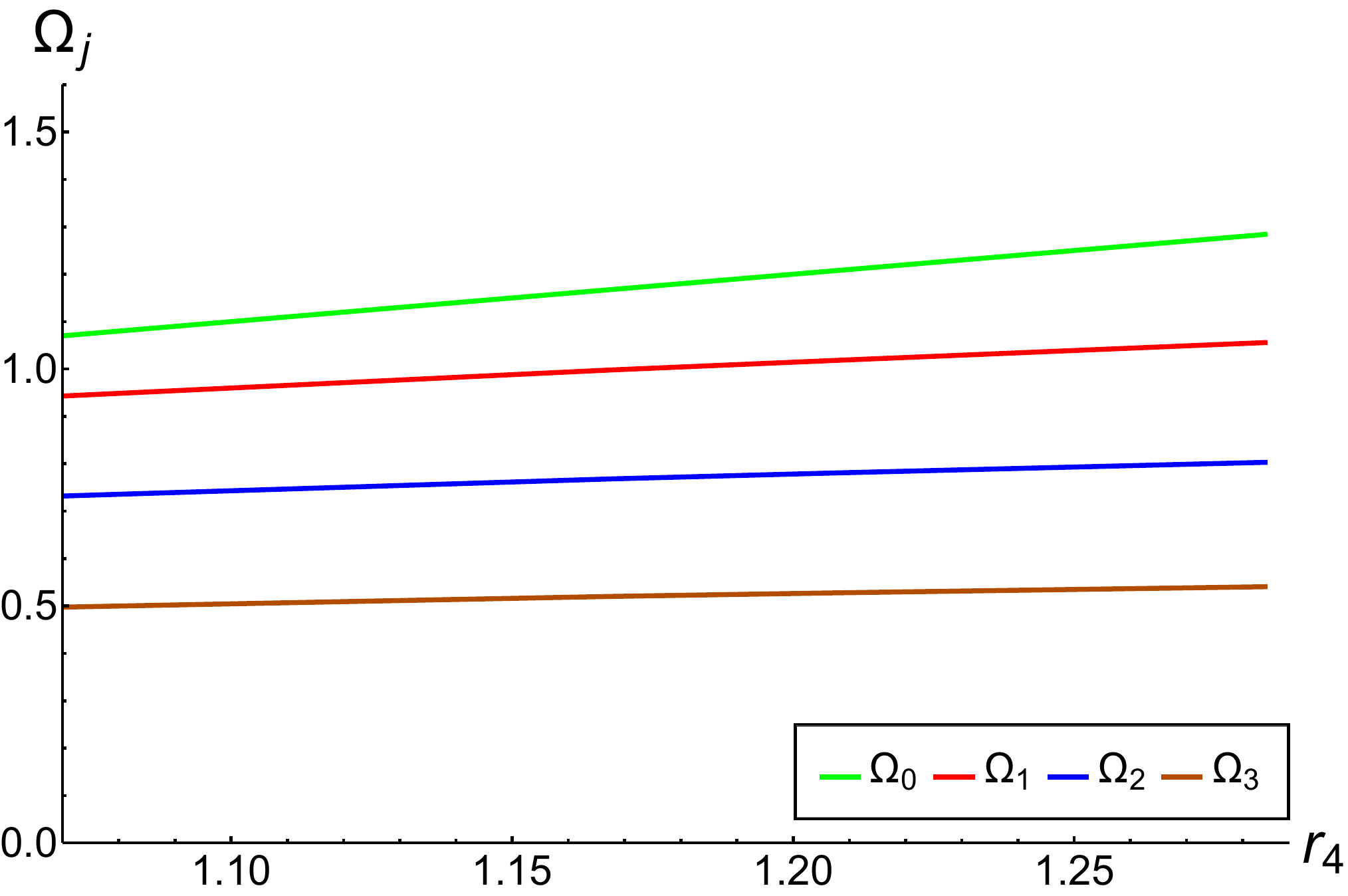}
		\caption{$\LL_j \simeq a_{j,4}r_4+b_{j,4}$}
		\label{fig:Search_Omega_k_modes_0_1_2_3_R6}
	\end{subfigure}%
		\caption{\small Representative examples of our search of unstable SLMs ($\LL>0,\ c_{\pm}^{(\mathcal{N})}=0$). These plots were obtained in $d=3$ but we find the same structure (under minor modifications) in other dimensions. In (a) we see the structure of $|c_{+}^{(\mathcal{N})}|$ which, in this case, has a single zero (red dashed lines) between consecutive integer values of $\mathcal{N}$ (black dashed lines). Recall that these zeros are associated with a SLM in case of SSs. The gray area is the region of $\LL$ where our method does not have accurate access. (b) Diagram ($r,\LL$) of unstable 1-modes. Here we see the region of non-admissible (red area $c_{+}^{(\mathcal{N})}\neq0,\ \mathcal{N}<1$) and admissible (blue areas $c_{+}^{(\mathcal{N})}\neq0,\ \mathcal{N}>1$) generic 1-modes. Among these modes we find some of them with $c_{+}^{(\mathcal{N})}=0$, $\LL_{j}(r)$, represented by purple lines. Yellow dashed lines show values of $r_n$ associated with SSs. Hence, intersections of purple and yellow lines mark points $(r_n,\LL_j)$ (purple points) that contain a SLM. Note that even the lowest SS $(r_2)$ has an unstable SLM (in addition to the gauge SLM $\LL_0$ (\ref{eq:1-mode_lambda_1})). Black dashed lines show the curves $\LL(r)$ where $\mathcal{N}$ is integer or $\nu-1$. In the dark blue sector ($\mathcal{N}>\nu-1$) our method to search 1-modes with $c_+^{(\mathcal{N})}=0$ does not work appropriately and this region was not deeply explored (gray area in plot (a)). Plot (c) shows the values of the coefficients $a_{1,n}(d)$ $b_{1,n}(d)$ of the unstable SLM $\LL_1$ associated with the lowest SS ($r_2$). We observe that these coefficients are greater than zero and seem to saturate when $d$ grows. It indicates that this mode is also unstable ($\LL_1>0$) in higher dimensions. Finally, plot (d) shows that the exponents $\LL_j$ have an almost linear dependence on $r_n(\ell)$ for fixed $n$. In this case we show the exponents $\LL_j$ for $j=0,1,2,3$  associated with the second family of SSs, $r_4$.}
		\label{fig:Search_Omega_k_1-modes_new}
	\end{figure}

	The family of 1-modes does not exclude linear contributions from the 0-modes. Note that the linear combination of the 0-mode and any 1-mode (both associated with the same $\LL$) is a 1-mode. For this reason, the properties of 0-modes described above apply to generic 1-modes as well. Their structure at $Z_2$ is of the form ($\xi := Z-Z_2$)
	\begin{align}
	& \alpha(\xi) = \underbrace{\sum \tilde{\alpha}_n\xi^n }_{\text{integer powers}} + \underbrace{c_{\pm}^{(\mathcal{N})} \ |\xi|^{\mathcal{N}}\sum\left(...\right) + c_{\pm}^{(\nu)} \ |\xi|^{\nu-1}\sum\left(...\right)}_{\text{non-integer powers}}, \label{eq:expansion_1-modes_alpha_non_integer_Z2_plusminus}\\
	& \beta(\xi) = \underbrace{\sum \tilde{\beta}_n\xi^n }_{\text{integer powers}} + \underbrace{c_{\pm}^{(\mathcal{N})} \ |\xi|^{\mathcal{N}}\sum\left(...\right)  + c_{\pm}^{(\nu)} \ |\xi|^{\nu-1} \sum\left(...\right)}_{\text{non-integer powers}},
	\label{eq:expansion_1-modes_beta_non_integer_Z2_plusminus}
	\end{align}	
	where $c_+^{(\mathcal{N})}$ and $c_+^{(\nu)}$ ($c_-^{(\mathcal{N})}$ and $c_-^{(\nu)}$) are associated with trajectories on the interior, $Z<Z_2$, (exterior, $Z>Z_2$) of the acoustic cone, $\mathcal{N}(\LL)$ is the regularity of 0-modes and $\nu$ the regularity of NSSs. From these expressions we see that there are two sources of non-smoothness $|\xi|^{\mathcal{N}}$ and $|\xi|^{\nu-1}$.  Therefore, the regularity of generic 1-modes is the minimum between $\mathcal{N}$ and $\nu-1$. However, we have found that for fine-tuned values of the parameters, the coefficients $c_{\pm}^{(\mathcal{N})}$ vanish. Then, these modes have regularity $\nu-1$ governed by NSSs and in case of SSs they are smooth linear modes (SLMs) (see fig.~\ref{fig:Search_Omega_k_1-modes_nu-der_smooth_non-smooth_profiles}). We will say that 1-modes with $c_{\pm}^{(\mathcal{N})}=0$ do not have contributions from 0-modes. To construct them we exploit the following structure (appedix~\ref{sec:Appendix_Regularity_Linear_Modes_Z2}):
	\begin{itemize}
		\item Given $(d,\ell,r,\kappa,\Omega)$ the nontrivial (regular) trajectory $(\alpha,\beta)$ that connects the origin with $Z_2$ (interior of the cone) is unique up to scaling.
		
		\item Given $(d,\ell,r,\kappa,\Omega)$ trajectories $(\alpha,\beta)$ that connect $Z_2$ with the infinity (exterior of the cone) belong to an one-parameter family. This freedom is materialized at infinity
		\beq
		\alpha(Z)\sim \theta\frac{\beta_0}{Z^{(r-1)\ell+\LL}}+..., \qquad \beta(Z)\sim \frac{\beta_0}{Z^{r-1+\LL}}+...
		\eeq
		where $\beta_0$ represents the scaling symmetry and $\theta$ labels trajectories.
	\end{itemize}
	Note that this is analogous to the structure of self-similar solutions; therefore, the strategy to construct modes with $c_{\pm}^{(\mathcal{N})}=0$ is also the same. First, exploiting the uniqueness of ($\alpha,\beta$) on the interior of the acoustic cone we know that $c_{+}^{(\mathcal{N})}=c_{+}^{(\mathcal{N})}(\LL)$ while on the exterior $c_{-}^{(\mathcal{N})}=c_{-}^{(\mathcal{N})}(\LL,\theta)$. Hence, this problem is also reduced to two ``eigenvalue problems" where we determine $\LL$ such that $c_{+}^{(\mathcal{N})}(\LL)=0$ first and with this value we determine $\theta$ such that $c_{-}^{(\mathcal{N})}(\LL,\theta)=0$. From this construction we learn the following lessons about the structure of 1-modes with $c_{\pm}^{(\mathcal{N})}=0$ (for $\mathcal{N}<8$ and avoiding situations where $\mathcal{N}$ and $\mathcal{N}-\nu$ are integers)
	\begin{enumerate}
			\item $c_{\pm}^{(\mathcal{N})} = 0$ is not equivalent to $c_{\pm}^{(\nu)}=0$; see fig.~\ref{fig:Search_Omega_k_1-modes_nu-der_smooth_non-smooth_profiles}.
		
		\item $c_{\pm}^{(\nu)}$ vanishes for SSs.
		
		\item  $c_{+}^{(\mathcal{N})}(\LL)$ has zeros for discrete values of $\LL$ that we denote by $\{\LL_j\}$ such that they satisfy $\LL_0>\LL_1>\LL_2>...$ (see fig.~\ref{fig:Search_Omega_k_1-modes_A}). We have found the first values of this sequence for $\mathcal{N}<\nu-1$, some of them are provided in table~\ref{table:r_1_full}.
		
		\item $\{\LL_j\}$ are continuous functions of $r$ (at least where $\mathcal{N}<\nu-1$); see fig.~\ref{fig:Diagram_Smooth_1-modes} for a visual representation in the diagram $(r,\LL)$. Our strategy to construct modes with $c_{\pm}^{(\mathcal{N})}=0$ does not give accurate access to the whole space of modes, just a particular region $(\mathcal{N}<\nu-1)$. This is the reason that we cannot study $\LL_j(r)$ for all $r$.
		
		\item For each $\LL_j$ there is $\theta$ such that $c_{-}^{(\mathcal{N})}(\LL_j,\theta)=0$. This is a consequence that our equations are linear and we have found two independent solutions; the 0-mode and any 1-mode. Then, for SSs zeros of $c_{+}^{(\mathcal{N})}(\LL)$ are associated with SLMs.
		
		\item The greatest $\LL_j$ is $\LL_0=r$, this is associated with an artificial mode that we provide later in (\ref{eq:1-mode_lambda_1}). In case of SSs this is a SLM but note that this is not a real instability.
		
		\item For a given $r$, the coefficient $c_{+}^{(\mathcal{N})}$ can have more than a zero or none between consecutive integer values of $[\mathcal{N}]$. 
	\end{enumerate}
	Now, focusing on SSs ($d\geq2$, $\ell>0$ and $\nu<8$) we find that these solutions are unstable under smooth perturbations; they always have at least the SLM $\LL_1>0$. Even the SS associated with the lowest $r_n$, has this unstable mode. Furthermore, we find that, restricted to a family $r_n(\ell)$ of SSs, this exponent follows an almost linear expression
	\beq
		\LL_1\simeq a_{1,n}(d)r_n(\ell)+b_{1,n}(d),
		\label{eq:Omega_linear_ab}
	\eeq 
	where the coefficients depend on $d$ and $a_{1,k}\neq a_{1,n}$, for $k\neq n$. Calculating these coefficients in different dimensions we see that they are always positive and saturate when the dimension grows (see fig.~\ref{fig:Search_Omega_k_AB_coefficients_Omega-1}). Therefore, it indicates that this instability, $\LL_1>0$, is also present in dimensions higher than the ones that we have explored.
	
	 Our study also provides information about the number of unstable directions of SSs. We warn the reader that this number excludes mode $\LL_0$ because this is not a real instability as we will explain below. From our exploration we find that the number of unstable SLMs of SSs associated with $r_n(\ell)$ is equal or greater than $n-1$ for $d\geq3$ ($n$ for $d=2$). Recall that for $d\geq3$ ($d=2$) there are SSs only for $n=2,4,6,...$ ($n=1,3,...$). We expect that there are some extra SLMs given the structure of the space of modes in fig.~\ref{fig:Diagram_Smooth_1-modes}. Furthermore, given that when $\nu$ grows the region $\mathcal{N}<\nu-1$ grows and the region $\nu-1<\mathcal{N}<\nu$ shrinks, we may expect that for $\nu$ large enough the number of unstable SLMs goes almost like $n$. In practice this estimate is not useful when $\nu\gg1$, the calculation of $n$ requires to know a large number of zeros of $c_{+}(r)$. Then, we should relate $n$ with some quantity that we can calculate easily. In this case, given that we observe that for $\nu<8$ there is an almost linear relation $n \sim [\nu] + c$ and in \cite{MerleEuler} for $\nu\gg1$ as well, we may expect that for $\nu\gg1$ the number of unstable SLMs is dominated by $[\nu]$. Therefore, the number of unstable directions of SSs grows with $\nu$, going to infinity when $r\to\reye$. Note that this is a naive estimate because this is based on an optimistic extrapolation of the number of SLMs that we observe for $\nu<8$ to the case $\nu\gg1$. We do not know if there is some saturation in the number of modes or if some of them merge together in case that two exponents $\LL_j(r)$ collide. These processes may drastically reduce the number of SLMs.
	
Finally, explicit expressions for two 1-modes without contributions from 0-modes (SLMs in case of SSs) are easily constructed:
\begin{itemize}
	\item Scaling symmetry: the scaling symmetry
	\beq
	\hat\rho(\tau,Z)\to \eta^{\ell}\hat\rho(\tau,Z/\eta) \qquad \hat u(\tau,Z) \to \eta \hat u(\tau,Z/\eta)
	\eeq
	of equations (\ref{eq:Rho_dtau})-(\ref{eq:U_dtau}) has the following mode associated ($\eta = 1 + \epsilon$ with $|\epsilon|\ll 1$)
	\beq
	\LL = 0, \qquad \begin{cases}
		\alpha(Z) = -\ell\hat{\rho}_0(Z) + Z\partial_{Z}\hat{\rho}_0(Z)\\
		\beta(Z) = -\hat{u}_0(Z) + Z\partial_Z \hat{u}_0(Z)
	\end{cases}
	\label{eq:1-mode_lambda_0}
	\eeq
	
	\item Gauge instability: self-similar variables $(\tau,Z)$ given in (\ref{eq:Self-Similar_Variables}) depend on $T$ and any small deviation $T\to T+\epsilon$ ($|\epsilon|\ll 1$) makes that instead of observing a static profile $(\hat\rho_0,\hat u_0)$ we observe an exponential deviation driven by the following mode
	\beq
	\LL = r, \qquad  \begin{cases}
		\alpha(Z) = \ell(r-1)\hat{\rho}_0(Z) + Z\partial_{Z}\hat{\rho}_0(Z)\\
		\beta(Z) = (r-1)\hat{u}_0(Z) + Z\partial_Z \hat{u}_0(Z)
	\end{cases}
	\label{eq:1-mode_lambda_1}
	\eeq
\end{itemize}
	Note that these two modes cannot be considered as instabilities of self-similar solutions because they are the materialization of a symmetry and a choice of coordinates. However, they are very useful for code verification and provide explicit information about the structure of the problem. Moreover, they are two explicit examples of SLMs when they are associated with a SS. In case of NSSs, despite these modes come from continuous transformations they have less regularity than the NSS; they exactly miss one derivative, in agreement with (\ref{eq:expansion_1-modes_alpha_non_integer_Z2_plusminus}-\ref{eq:expansion_1-modes_beta_non_integer_Z2_plusminus}).

	\begin{table}[h!]
	\centering
	\begin{tabular}{|c|c|c|c|c|c|c|c|} 
		\hline
		$d$ &	$\ell$ & $r_2$ $\left(\pm 3\cdot{}10^{-6}\right)$ & $\kappa$ $\left(\pm 10^{-5}\right)$ & $\LL_1$ $\left(\pm 10^{-5}\right)$& $\theta$ $\left(\pm 10^{-5}\right)$& $\Delta\tilde{\omega}_{[\nu]}$  &  $\Delta\tilde{\omega}_{[\nu]+1}$ \\ [0.5ex] 
		\hline
		3 & 0.5 & 1.248159 & 0.03429 & 0.84448 & -0.91152 & $8 \cdot{} 10^{-9}$ & $2 \cdot{} 10^{-5}$  \\  [0.5ex]
		\hline
		3 & 1 & 1.198877 & 0.06456 & 0.82719& -1.29945 & $8 \cdot{} 10^{-9}$ & $4 \cdot{} 10^{-5}$ \\  [0.5ex]
		\hline
		3 & 1.5 & 1.166517 & 0.08968 & 0.81596 &-1.55846 & $4 \cdot{} 10^{-9}$ & $5 \cdot{} 10^{-5}$  \\  [0.5ex]
		\hline
		3 & 2 & 1.143517 & 0.11056 & 0.80796 & -1.74581 & $5\cdot{} 10^{-11}$ & $6 \cdot{} 10^{-5}$  \\  [0.5ex]
		\hline
		3 & 2.5 & 1.126267 & 0.12813 & 0.80190 & -1.88780  & $3\cdot{} 10^{-9}$ & $6 \cdot{} 10^{-5}$  \\  [0.5ex]
		\hline
		3 & 3 & 1.112816 & 0.14311 & 0.79711 & -1.99949  & $2\cdot{} 10^{-9}$ & $6 \cdot{} 10^{-5}$  \\  [0.5ex]
		\hline
		3 & 5 & 1.079404 & 0.18591 & 0.78491 & -2.28200  & $9\cdot{} 10^{-9}$ & $4 \cdot{} 10^{-5}$  \\  [0.5ex]
		\hline
		3 & 10 & 1.045911 & 0.23869 & 0.77200 & -2.60000  & $5\cdot{} 10^{-9}$ & $3 \cdot{} 10^{-5}$  \\  [0.5ex]
		\hline
		4 & 0.5 & 1.390371 & 0.01171 & 0.95142 & -1.85753 & $10^{-10}$ & $9 \cdot{} 10^{-6}$ \\
		\hline
		4 & 1 & 1.321444 & 0.03264 & 0.92218 & -2.35813 & $6 \cdot{} 10^{-10}$ & $2 \cdot{} 10^{-5}$ \\
		\hline
		4 & 1.5 & 1.273995 & 0.05457 & 0.90240 & -2.76720 & $10^{-9}$ & $3 \cdot{} 10^{-5}$ \\
		\hline
		4 & 2 & 1.239224 & 0.07536 & 0.88795 & -3.08776 & $2 \cdot{} 10^{-9}$ & $4 \cdot{} 10^{-5}$ \\
		\hline
		4 & 2.5 & 1.212575 & 0.09447 & 0.87683 & -3.33737 & $3 \cdot{} 10^{-9}$  & $4 \cdot{} 10^{-5}$ \\
		\hline
		4 & 3 & 1.191452 & 0.11189 & 0.86793 & -3.53260 & $3 \cdot{} 10^{-9}$  & $5 \cdot{} 10^{-5}$ \\
		\hline
		4 & 3.5 & 1.174271 & 0.12771 & 0.86062 & -3.68554 & $3 \cdot{} 10^{-9}$  & $5 \cdot{} 10^{-5}$ \\
		\hline
	\end{tabular}
	\caption{\small Some SSs associated with the lowest $r_n(\ell)$  and their first unstable SLM given by $(\Omega_1,\theta)$. For these parameters $[\nu]=3$ and $[\mathcal{N}]=1$. $\Delta\tilde{\omega}_{j}:=|\tilde{\omega}_{j}^{(N)}-\tilde{\omega}_{j}^{(A)}|/|\tilde{\omega}_{j}^{(A)}|$, where $\tilde{\omega}_{j}^{(N)}$ represents the numerical value that we obtain for the coefficient $\tilde{\omega}_{j}$ in (\ref{eq:expansion_w_in_powers_sigma_Z2_plusminus}) and $\tilde{\omega}_{j}^{(A)}$ its analytic value obtained from the expansion at $P_2$ imposing $c_{\pm}=0$. These quantities are used to quantify the accuracy of our numerical results. We warn the reader that these values of $(r,\kappa,\Omega,\theta)$ do not provide the exact SSs and SLMs. The levels of uncertainty take into account that these parameters may slightly differ between different numerical implementations. The reader has to scan NSSs around these values to add extra digits.}
	\label{table:r_1_full}
\end{table}

	\newpage
	\section{The Nonlinear Time-Evolution: Endpoints of Instabilities \vspace{3mm}}\label{sec:Quasi_Shock_Formation}

	After constructing self-similar solutions and their linear modes we wonder what are the endpoints of unstable directions triggered by these perturbations. For example, is the singularity formation process unstable? To provide information about this question a numerical time-evolution is carried out with the initial data prepared as follows 
	\beq
	\hat\rho(0,Z) = \hat\rho_{0}(Z) + \epsilon \alpha(Z), \qquad \hat{u}(0,Z) = \hat{u}_0(Z) + \epsilon \beta(Z),
	\label{eq:initial_data_blow-up_plus_perturbations}
	\eeq
	where $|\epsilon|\ll 1$, $(\hat\rho_0,\hat{u}_0)$ is the blow-up profile, $(\alpha,\beta)$ a linear mode associated with $\LL$ and we fix the scale such that $\max_z\alpha = 1 > |\min_z\alpha|$. Given the two kinds of regularity in this problem we discuss two scenarios in different subsections: the generic behavior of instabilities and the behavior of SSs under unstable SLMs.

	\subsection{Unstable Directions of Non-Smooth Perturbations: A Shock-Formation}

	 We will focus this discussion on 0-modes because this is the cleanest case to extract conclusions. Given that generic 1-modes have contributions from 0-mode, these conclusions are easily adapted to them. Furthermore, our results are valid for both SSs and NSSs under generic perturbations. The main statement extracted from our numerical experiments is the following
	
	\vspace{0.4cm}
	\textbf{Numerical Claim:} {\em Given an initial data at $\tau=0$ prepared as in (\ref{eq:initial_data_blow-up_plus_perturbations}) with a generic 0-mode with $\LL\in(0,\LL_{max})$ and $0<\epsilon\ll 1$; then, its initial regularity is $\mathcal{N}(\LL)$ given in (\ref{eq:regularity_lambda_0_modes}), and there exists a time $0<\tau^*<\infty$ such that $(\hat{\rho},\hat{u})$ remain finite but their gradients go to infinity at a single point out of the origin; namely, a shock formation.}

	\vspace{0.4cm}
	\textbf{Remarks:}
	\vspace{-0.5cm}
	\begin{enumerate}
		\item The shock forms at $Z>Z_2$ before the original singularity associated with the self-similar solutions.
		\item The original singularity formation associated with self-similar solutions is unstable under perturbations with less regularity than $\nu$, but arbitrarily close.
		
		\item There is no evidence of loss of regularity in finite $\tau$ when the perturbations described in the {\em numerical claim} have the opposite sign ($\epsilon < 0$). This remark is subject to the maximum time that we can simulate (see the numerical ramp in appendix~\ref{sec:Appendix_Numerical_Methods}).
	
		\item \label{item:remarks_1-modes} The fact that 1-modes have linear contributions from 0-modes makes that the {\em numerical claim} usually holds as our numerical experiments show. 
	\end{enumerate}

	\vspace{0.4cm}
	\textbf{Justification for the numerical claim:}
	\vspace{-0.3cm}
	
\begin{figure}[h!]
	\centering	
	\begin{subfigure}[b]{0.35\textwidth}
		\centering
		\includegraphics[width=5.5cm]{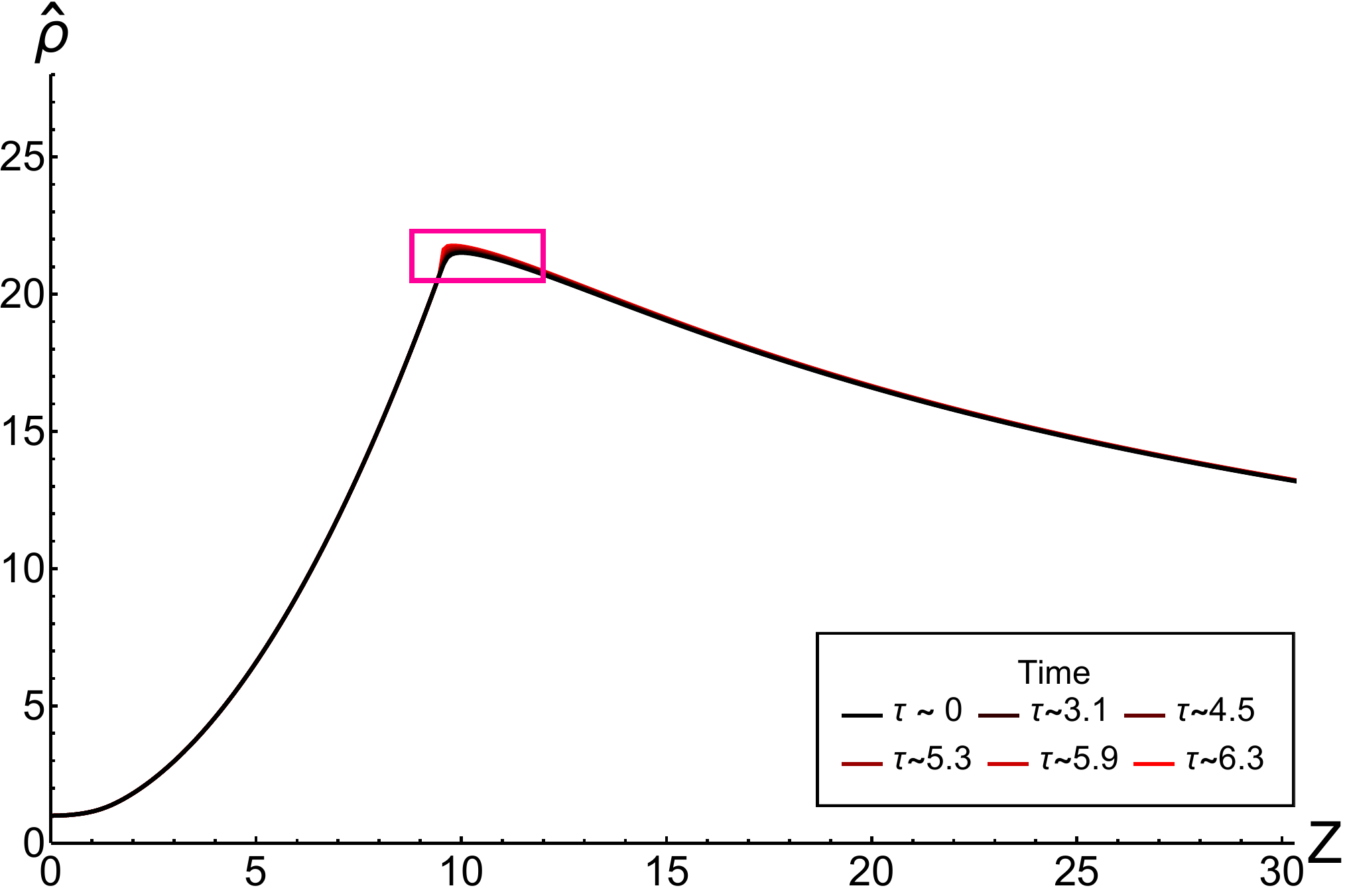}
		\caption{}
	\end{subfigure}%
	\begin{subfigure}[b]{0.35\textwidth}
		\centering
		\includegraphics[width=5.5cm]{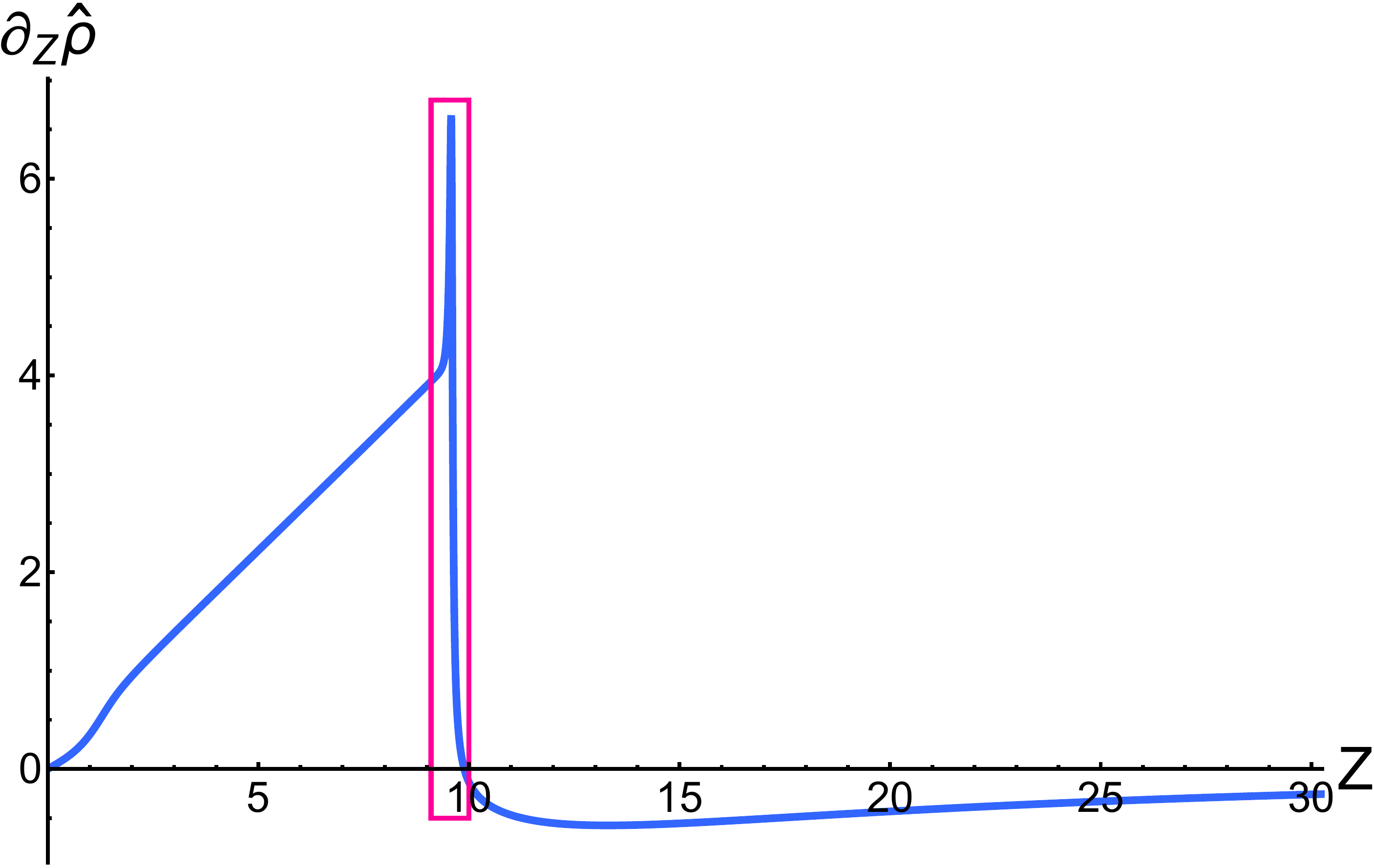}
		\caption{}
	\end{subfigure}%
	\begin{subfigure}[b]{0.35\textwidth}
		\centering		\hspace{2cm}\includegraphics[width=5.5cm]{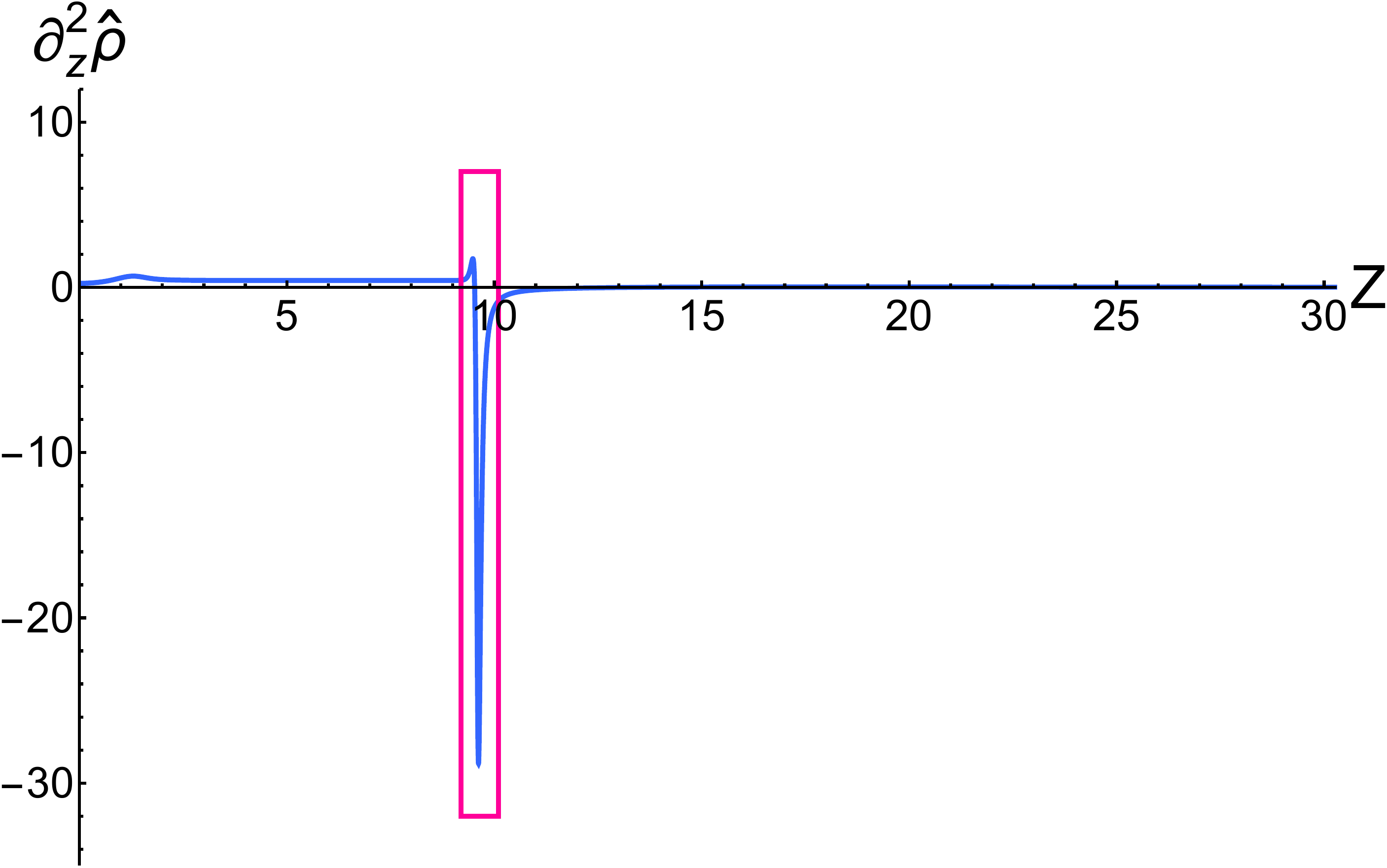}
		\caption{}
	\end{subfigure}%

\vspace{0.2cm}
	
	\begin{subfigure}[b]{0.35\textwidth}
	\centering		\hspace{2cm}\includegraphics[width=5.5cm]{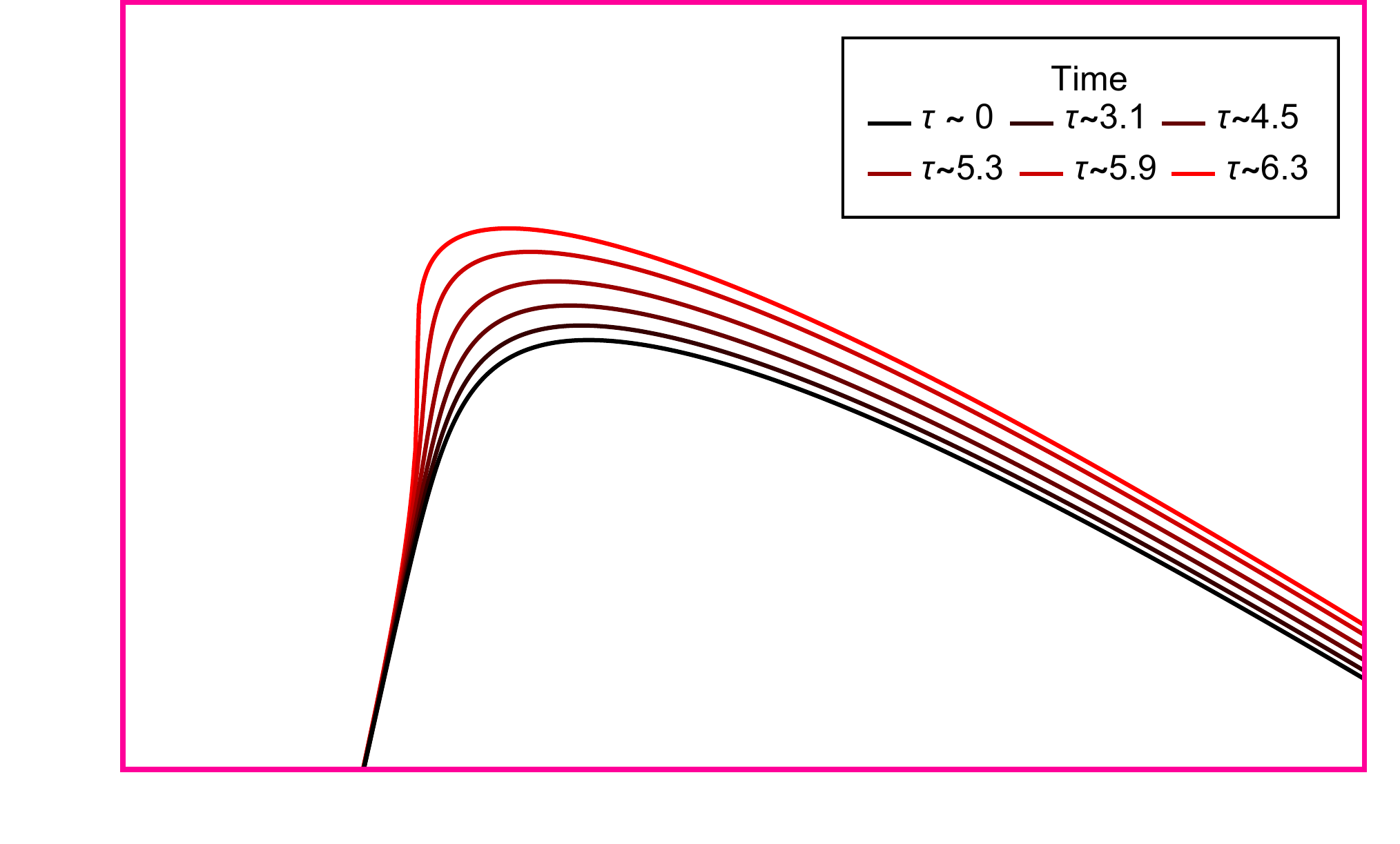}
	\caption{Pink frame in (a)}
\end{subfigure}%
	\begin{subfigure}[b]{0.35\textwidth}
	\centering		\hspace{2cm}\includegraphics[width=5.5cm]{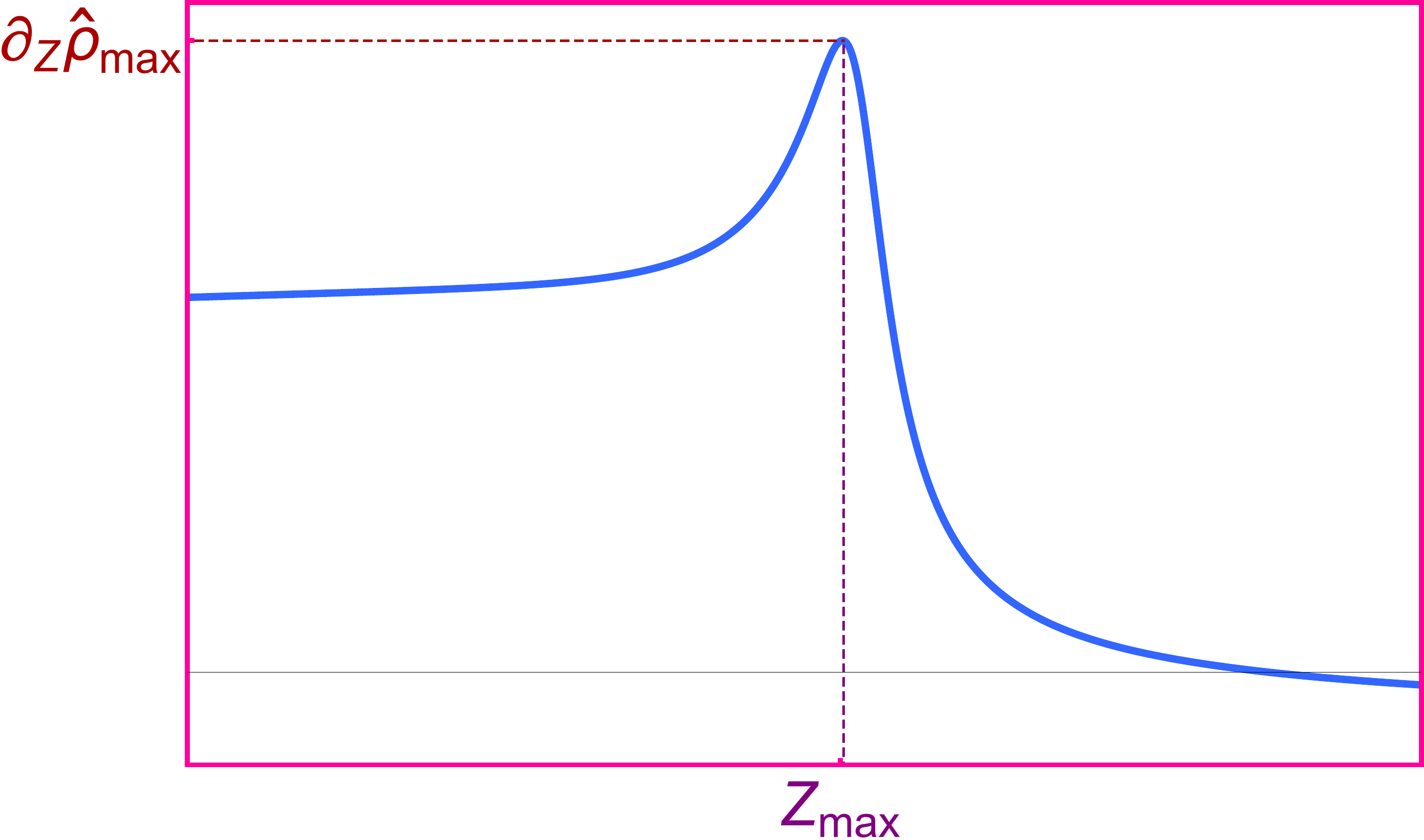}
	\caption{Pink frame in (b)}
\end{subfigure}%
	\begin{subfigure}[b]{0.35\textwidth}
	\centering		\hspace{2cm}\includegraphics[width=5.5cm]{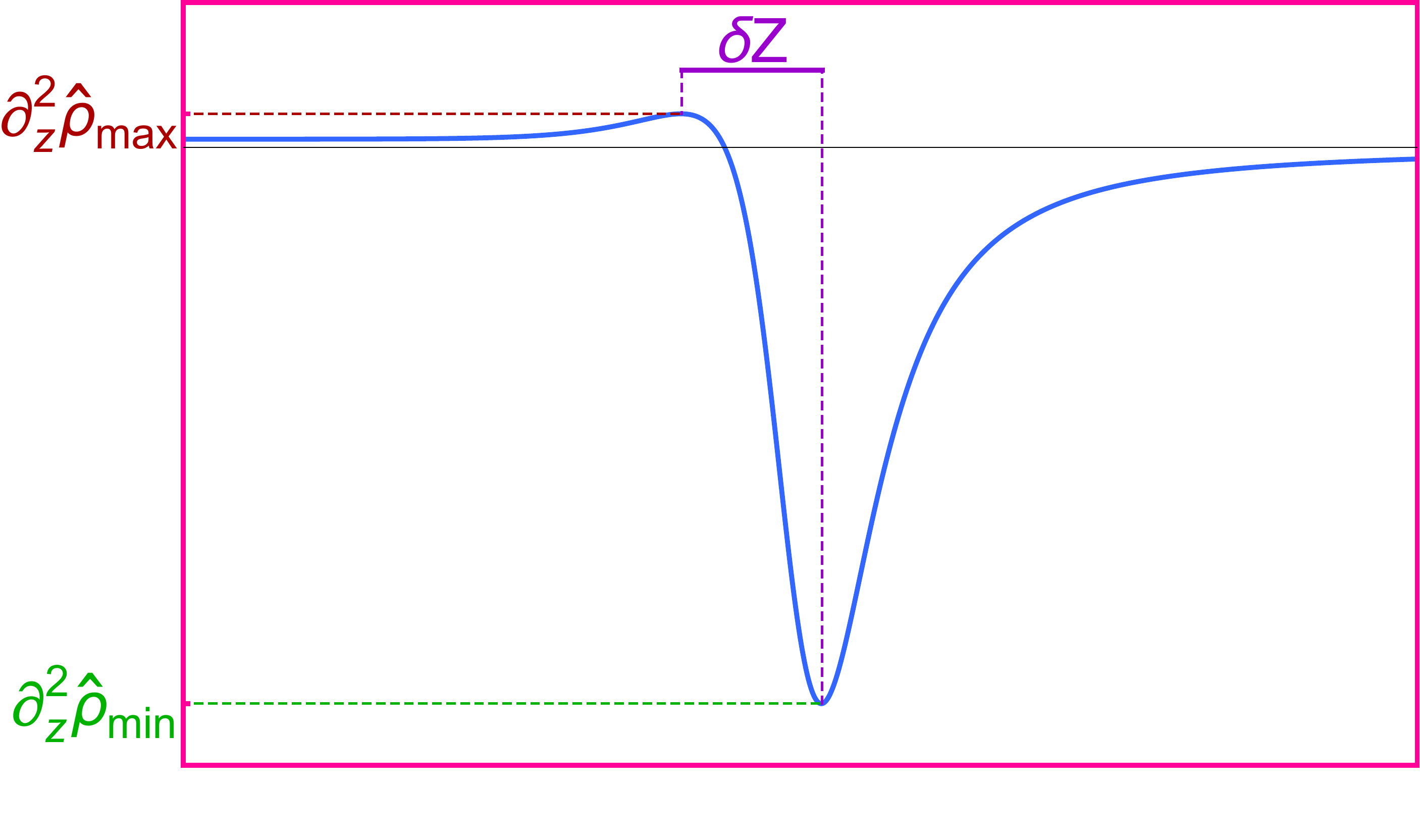}
	\caption{Pink frame in (c)}
\end{subfigure}%
	\caption{\small (a, d): time-evolution, in $(\tau,Z)$, of the NSS $(d,\ell,r,\kappa)=(3,2,1.33,0.6)$ and unstable 0-mode $\LL=1/2$ with $\epsilon = 10^{-2}$ following (\ref{eq:initial_data_blow-up_plus_perturbations}).  (b, e): generic shape of $\partial_z\hat\rho$. (c, f): generic shape of $\partial_z^2\hat\rho$. (e, f) show a visual representation of important quantities in the determination of the shock formation; their time-evolution can be found in fig.~\ref{fig:tau-d2Z}.}
	\label{fig:Rho_Shock}
\end{figure}

\begin{figure}[h!]
	\centering	
	\begin{subfigure}[b]{0.5\textwidth}
		\centering		\hspace{2cm}\includegraphics[width=8cm]{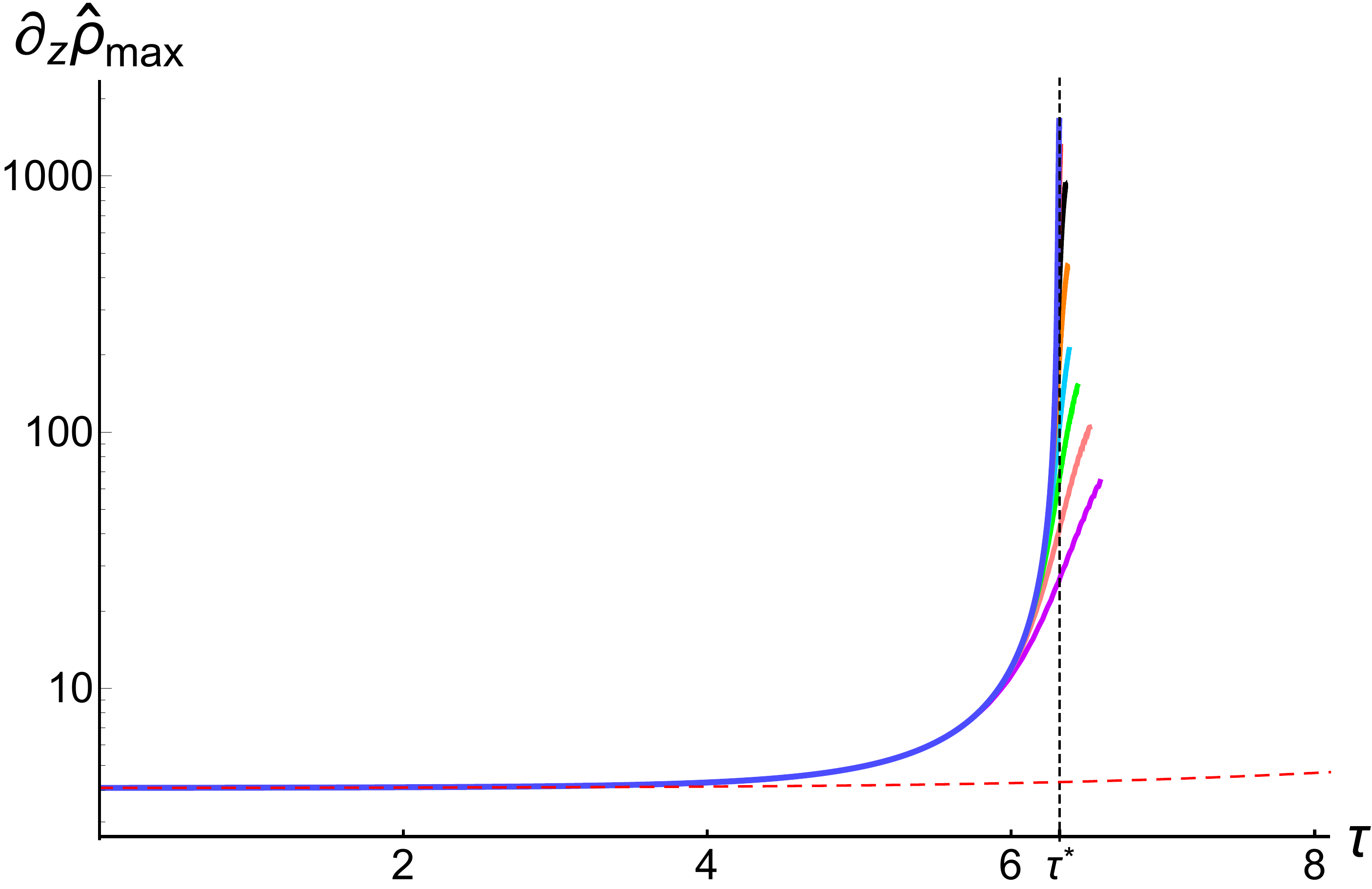}
		\caption{$\max_z\left(\partial_z\hat\rho\right)$}
	\end{subfigure}%
	\begin{subfigure}[b]{0.5\textwidth}
		\centering
		\includegraphics[width=8cm]{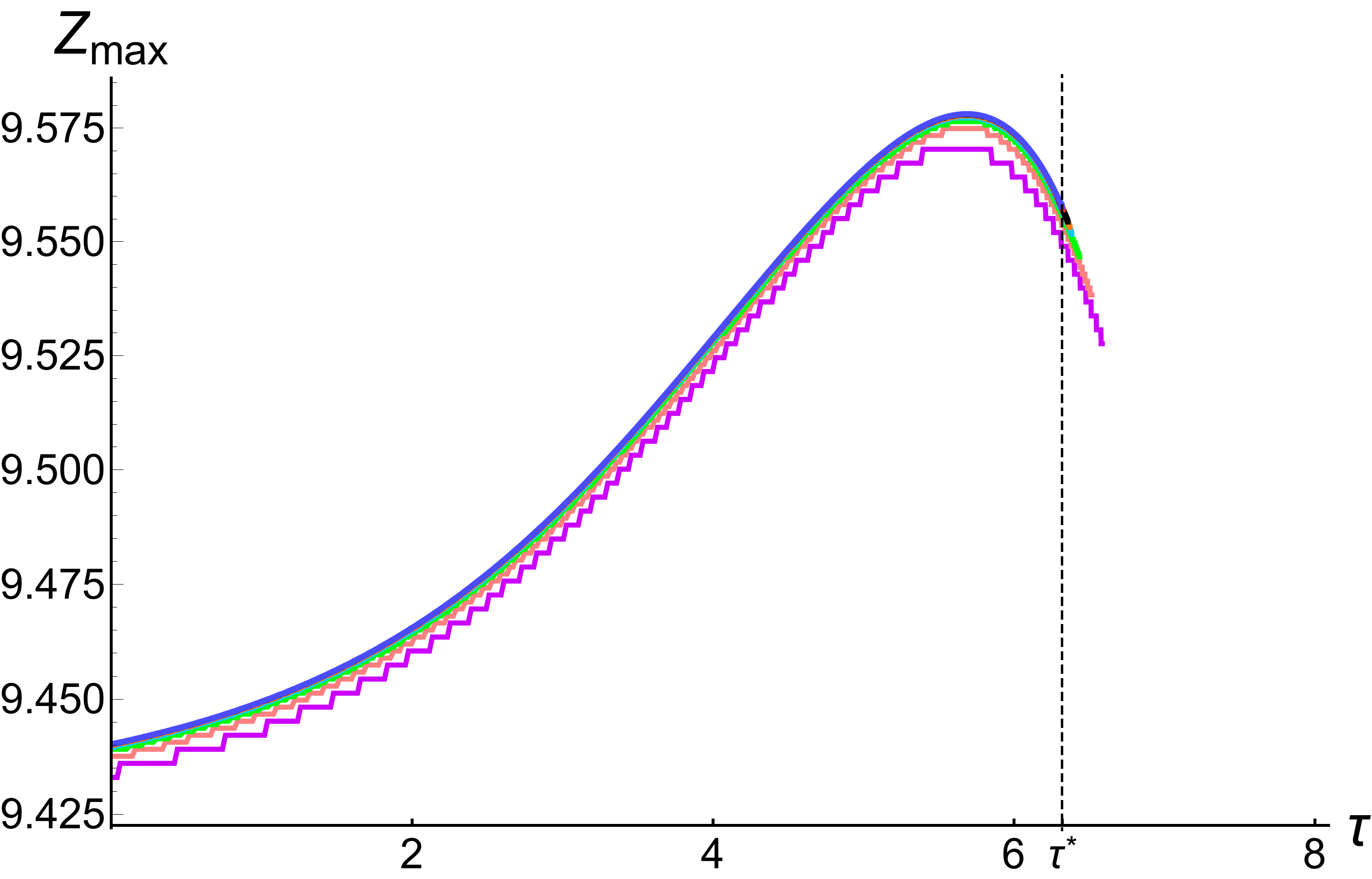}
		\caption{$\partial_z\hat\rho(Z_{max})=\partial_z\hat\rho_{max}$}
	\end{subfigure}%
	
	\vspace{0.cm}
	\begin{subfigure}[b]{0.5\textwidth}
		\centering
		\includegraphics[width=8cm]{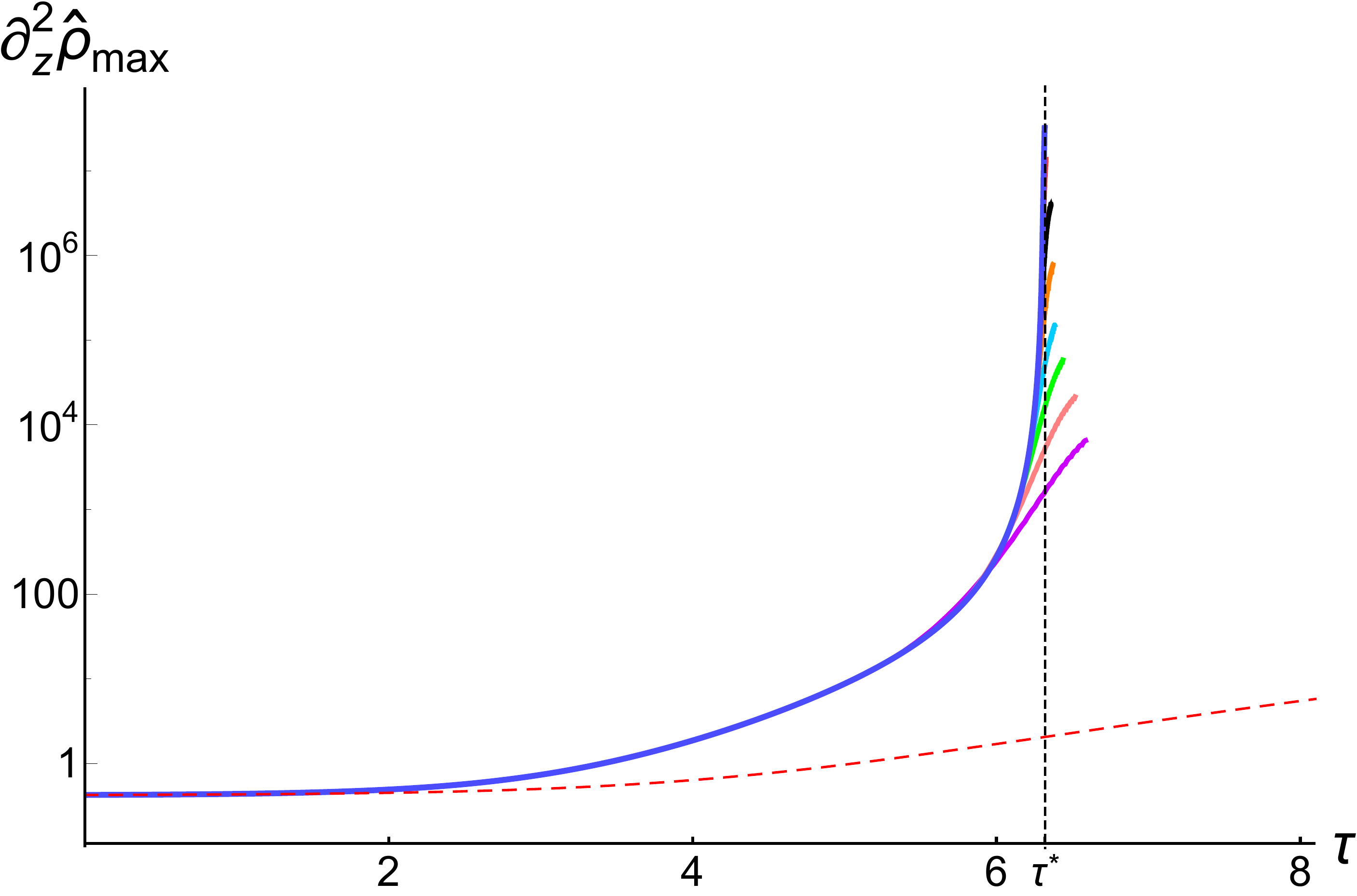}
		\caption{$\max_z\left(\partial_z^2\hat\rho\right)$}
	\end{subfigure}%
	\begin{subfigure}[b]{0.5\textwidth}
		\centering		\hspace{2cm}\includegraphics[width=8cm]{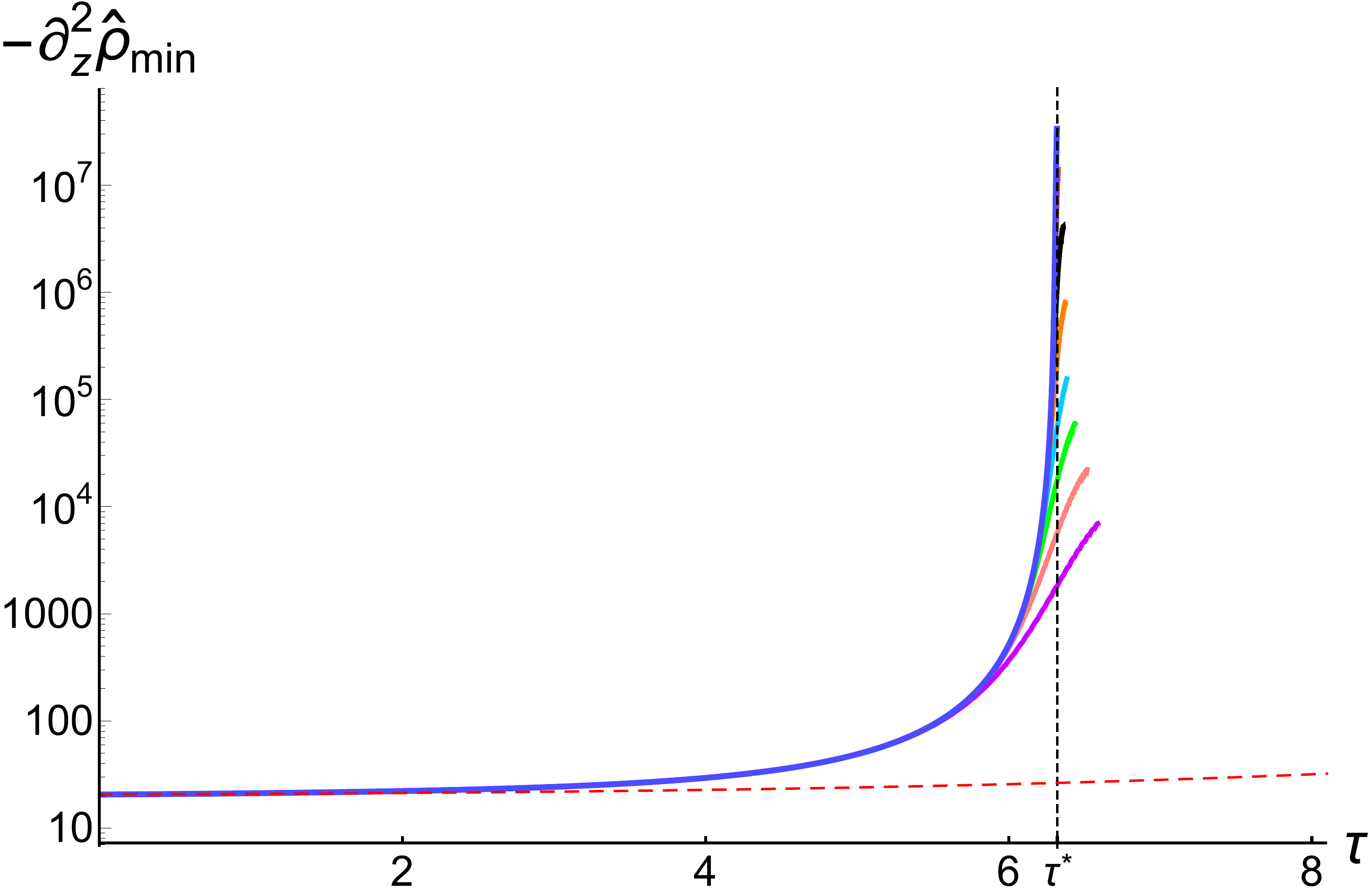}
		\caption{$-\text{min}_z\left(\partial_z^2\hat\rho\right)$}
	\end{subfigure}%
	
	\vspace{0.cm}
	\begin{subfigure}[b]{0.5\textwidth}
		\centering
		\includegraphics[width=8cm]{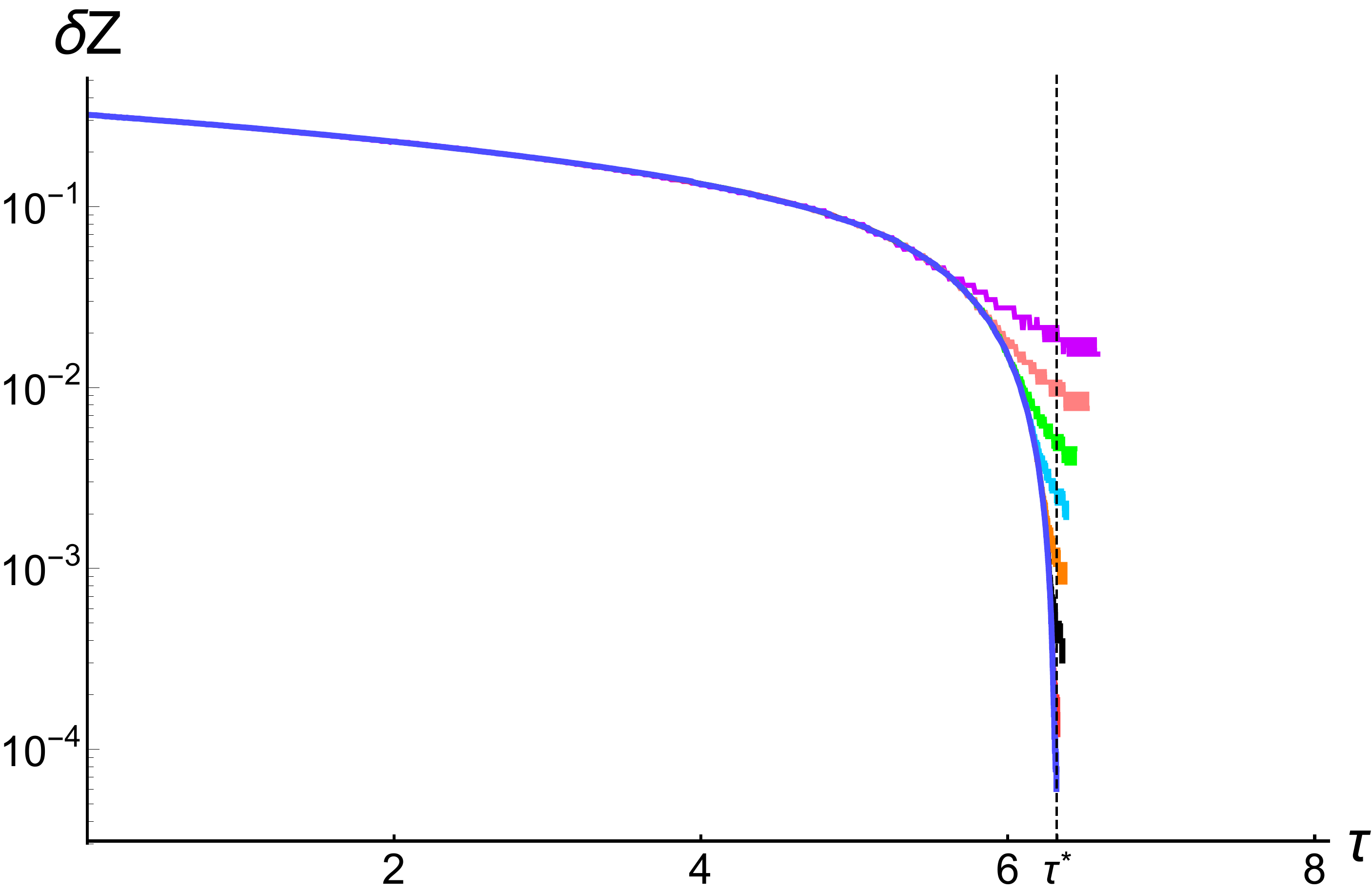}
		\caption{$\delta Z := |Z_{\partial_z^2\hat\rho_{max}} - Z_{\partial_z^2\hat\rho_{min}}|$}
	\end{subfigure}%
	\begin{subfigure}[b]{0.5\textwidth}
		\centering		\hspace{2cm}\includegraphics[width=8cm]{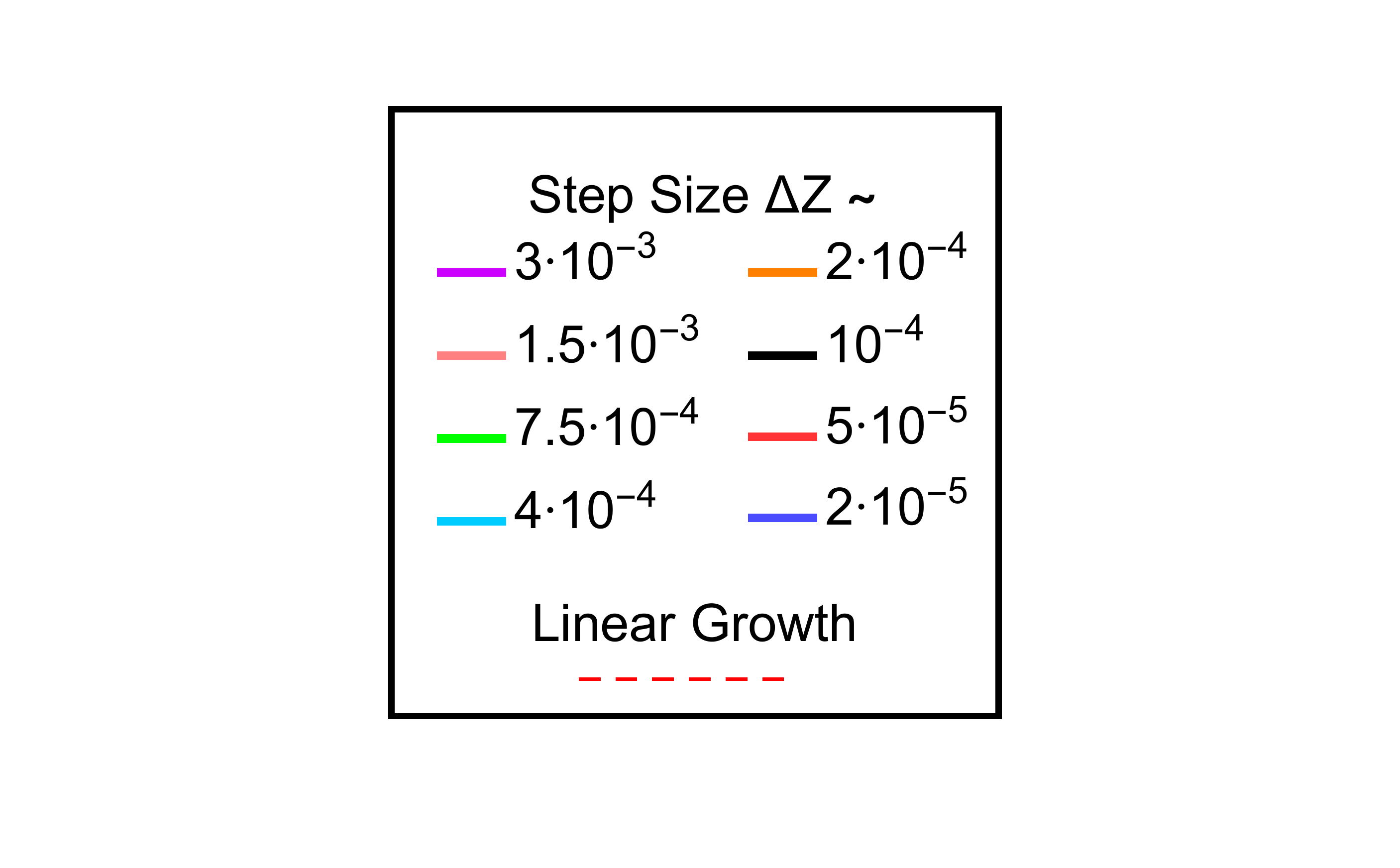}
		\caption{ }
	\end{subfigure}%
	\caption{\small Time-evolution of the maximum of $\partial_z\hat\rho$ (a) the position of this maximum (b), the maximum (c) and minimum (d) of $\partial_z^2\hat\rho$, the distance in $Z$ between them (e) and the legend in these plots (f). From these plots we see that the solution develops a sharp structure (spatial derivatives rapidly grow) when $\tau$ approaches $\tau^*$. Simulations of the same initial data using different resolution ($\Delta Z$) suggest that the spatial derivatives go to infinity for $\tau=\tau^*$ as we argue in fig.~\ref{fig:loglog_shock_formation_plots_Rho_U}. Plots shown in this figure are a representative example with the initial data given in fig.~\ref{fig:Rho_Shock} and $\tau^*\sim 6.3$, $Z_2\sim 5$ but the conclusions are the same in our range of parameters (\ref{eq:parameters}).}
	\label{fig:tau-d2Z}
\end{figure}

\begin{figure}[h!]
	\centering	
	\begin{subfigure}[b]{0.5\textwidth}
		\centering		\hspace{2cm}\includegraphics[width=8cm]{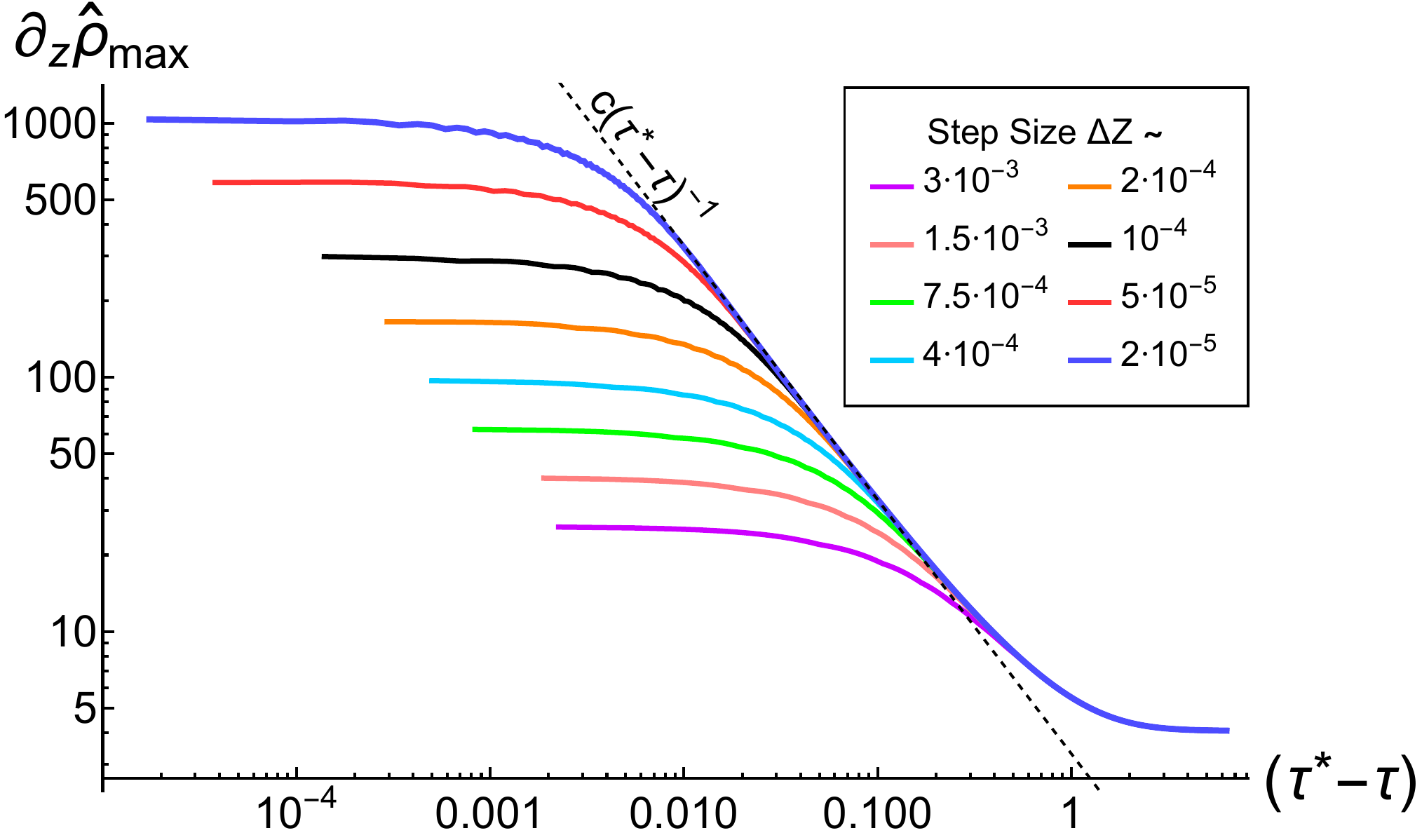}
		\caption{$\max_z\left(\partial_z\hat\rho\right)$}
	\end{subfigure}%
	\begin{subfigure}[b]{0.5\textwidth}
		\centering		\hspace{2cm}\includegraphics[width=8cm]{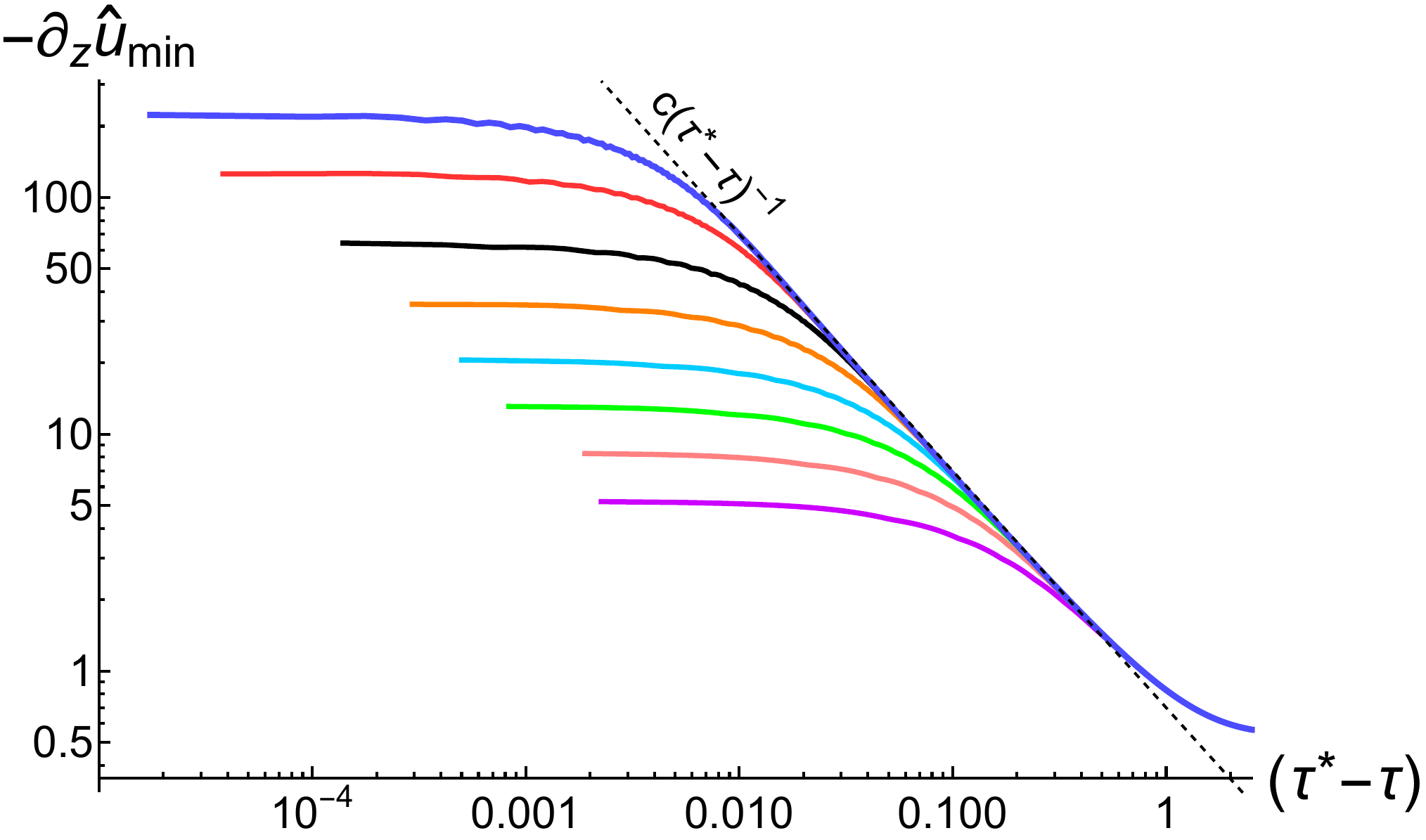}
		\caption{$-\min_z\left(\partial_z\hat{u}\right)$}
	\end{subfigure}%

	\vspace{0.5cm}
	\begin{subfigure}[b]{0.5\textwidth}
		\centering		\hspace{2cm}\includegraphics[width=8cm]{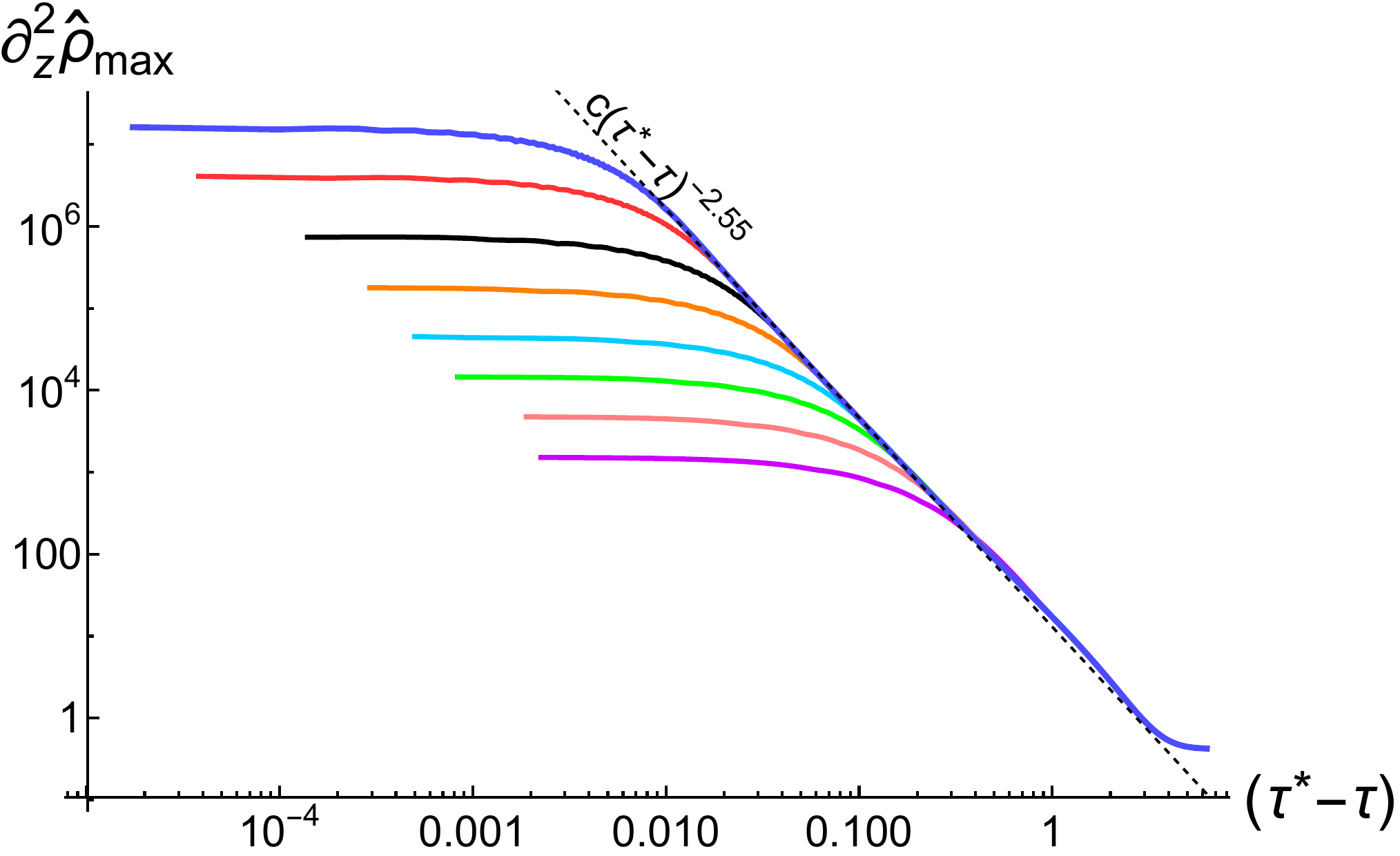}
		\caption{$\max_z\left(\partial_z^2\hat\rho\right)$}
	\end{subfigure}%
	\begin{subfigure}[b]{0.5\textwidth}
		\centering		\hspace{2cm}\includegraphics[width=8cm]{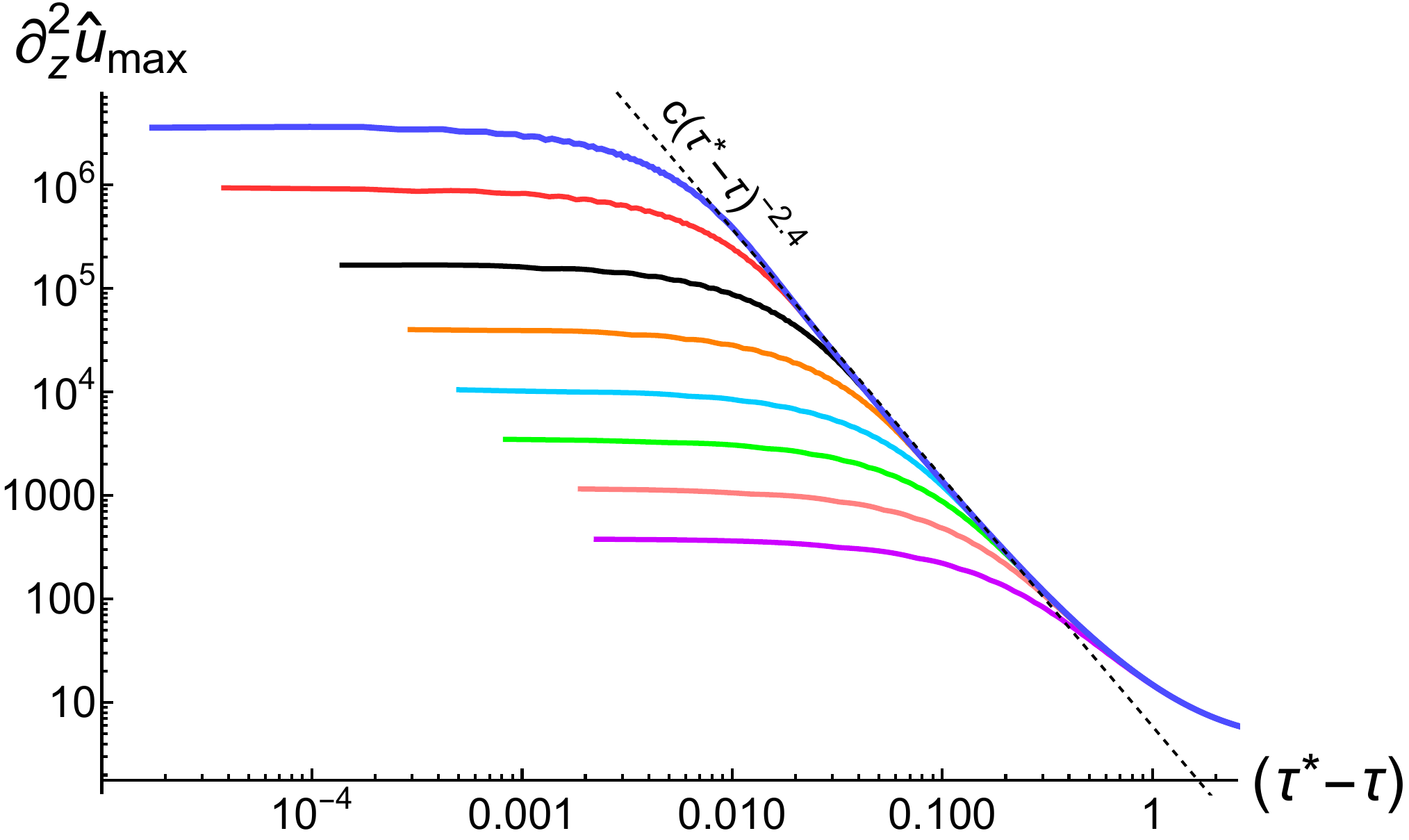}
		\caption{$\max_z\left(\partial_z^2\hat{u}\right)$}
	\end{subfigure}%
	
	\vspace{0.5cm}
	\begin{subfigure}[b]{0.5\textwidth}
		\centering
		\includegraphics[width=8cm]{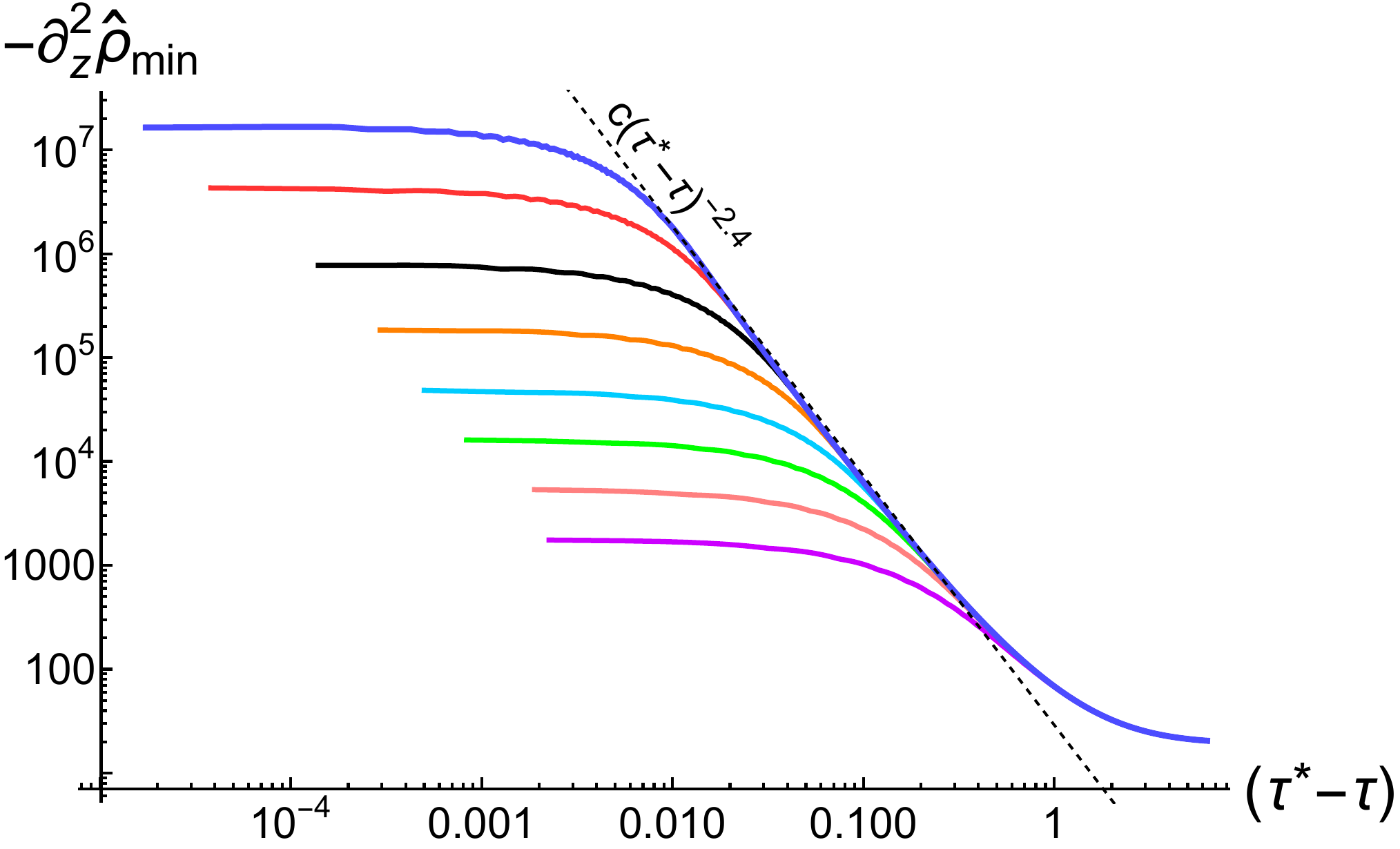}
		\caption{$-\text{min}_z\left(\partial_z^2\hat\rho\right)$}
	\end{subfigure}%
	\begin{subfigure}[b]{0.5\textwidth}
		\centering
		\includegraphics[width=8cm]{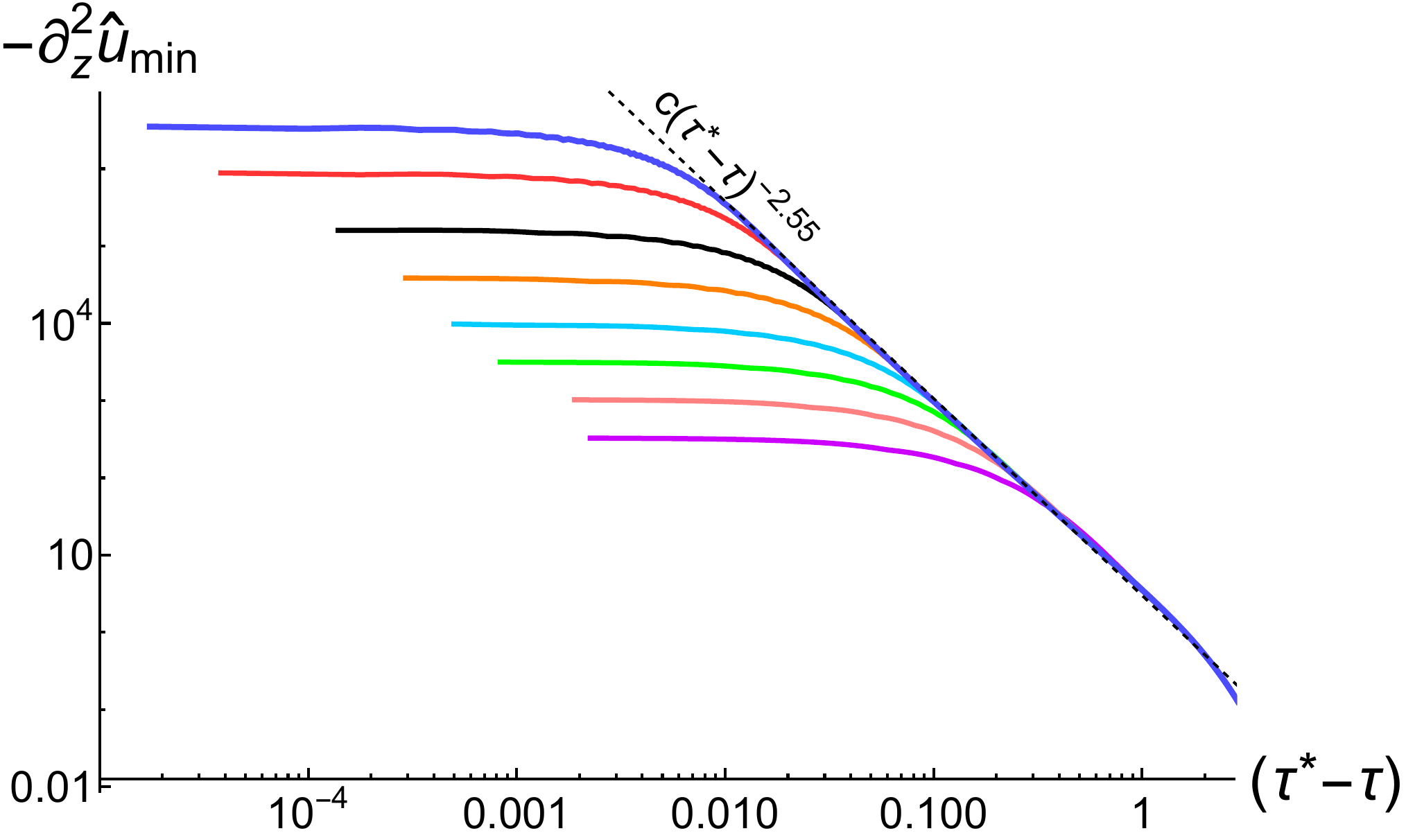}
		\caption{$-\text{min}_z\left(\partial_z^2\hat{u}\right)$}
	\end{subfigure}%
	
	\caption{\small These plots show, in a log-log scale and following the legend in (a), how the first and second derivatives of $\hat{\rho}$ and $\hat{u}$ (maximum and/or minimum) evolve close to $\tau^*$. We observe that while the simulation converges these quantities evolve as $c(\tau^*-\tau)^{-s}$ with $s> 0$. At some point the resolution (quantified by the step-size $\Delta Z$) is not enough to describe the sharp structure developed in the density and velocity and we observe a deviation from $c(\tau^*-\tau)^{-s}$. However, simulations of higher resolution (smaller $\Delta Z$) show that this deviation does not actually happen.\\ \\}
	\label{fig:loglog_shock_formation_plots_Rho_U}
\end{figure}

	The {\em numerical claim} is formulated after we have performed a large number of numerical experiments. The main obstacle that one has to face when tries to numerically describe the formation of a singularity is the finite resolution (quantified by the step-size $\Delta Z$). At some point, close to the singularity, the numerical evolution is not able to accurately describe the process; for this reason it is important to understand the behavior of such simulations for different $\Delta Z$ to extract conclusions that provide intuition for $\Delta Z \to 0$. From our simulations we extract the following points (a representative example can be found in fig.~\ref{fig:Rho_Shock},~\ref{fig:tau-d2Z},~\ref{fig:loglog_shock_formation_plots_Rho_U}):
	\begin{itemize}		
		\item The region where our simulations converge shows that close to $\tau^*$ the maximum of the first spatial derivative of the density ($\partial_Z\hat{\rho}_{\max}$), the minimum of the first spatial derivative of the velocity ($\partial_Z\hat{u}_{\min}$) and the maximum and minimum of their second spatial derivatives ($\partial_Z^2\hat{\rho}_{\max}$, $\partial_Z^2\hat{\rho}_{\min}$, $\partial_Z^2\hat{u}_{\max}$, $\partial_Z^2\hat{u}_{\min}$) exhibit a growth of the form ${c(\tau^*-\tau)^{-s}}$ with $s>0$ close to $\tau^*$. Specifically, for  $\partial_Z\hat{\rho}_{\max}$ and $\partial_Z\hat{u}_{\min}$ we find $s\sim 1$, for $\partial_Z^2\hat{\rho}_{\max}$ and $\partial_Z^2\hat{u}_{\min}$ that $s \sim 2.55$ while for $\partial_Z^2\hat{\rho}_{\min}$ and $\partial_Z^2\hat{u}_{\max}$ that $s\sim 2.4$; however, these values for $s$ should be understood as a rough estimate for the reasons explained below.
		
		\item The quantities mentioned in the previous point deviate from a growth of the form $c(\tau^*-\tau)^{-s}$ because the resolution of the spatial grid is not enough to describe the process. The point of deviation grows when $\Delta Z$ decreases.
	\end{itemize}
	These points suggest that for the exact system of equations there exists a time $\tau^*$ when the spatial derivatives of $\hat{\rho}$ and $\hat{u}$ go to infinity (at a single point) while $\hat{\rho}$ and $\hat{u}$ remain finite. Therefore, when $\tau$ approaches $\tau^*$ these functions develop a region where their values have an abrupt change, ending in the formation of a shock in finite $\tau$. 
	
	At this point it is worth to mention that the values for the exponent $s$ of $c(\tau^*-\tau)^{-s}$ provided above should be understood by the reader as rough estimates because they are subject to significant uncertainties. To reduce these uncertainties one has to perform simulations with higher resolution (much smaller $\Delta Z$) in order to get accurate access to regions closer to $\tau^*$. With it the leading terms become much more dominant. However, smaller $\Delta Z$ is a challenge due to the needs of this problem (see appendix~\ref{sec:Appendix_Numerical_Methods}). An hydro-code specifically developed to simulate shocks could provide a significant boost in this direction. For example, one expects that both the the maximum and minimum of the second derivatives grow with the same exponent. An improvement in the resolution could show that $2.4$ and $2.55$ converge to the same value.
	
	Finally, we want to remark that in order to get a robust confirmation of the shock formation for $\epsilon \ll 1$, we have made use of sophisticated numerical techniques. Among others, the numerical construction of highly accurate blow-up profiles and the use of extended precision operations. These two points are crucial when $\epsilon$ is very small. If the arithmetic precision is restricted to standard ``double-numbers" we observe that for $\epsilon$ small enough unstable directions of a self-similar solution intersect stable directions of other self-similar solutions. This stability is {\em spurious} because when the arithmetic precision is extended the shock formation happens again. Further details about our numerical methods can be found in appendix~\ref{sec:Appendix_Numerical_Methods}.


	\subsection{Smooth Solutions under Smooth Perturbations}

\begin{figure}[h!]
		\centering	
		\begin{subfigure}[b]{0.5\textwidth}
			\centering		\hspace{2cm}\includegraphics[width=7.5cm]{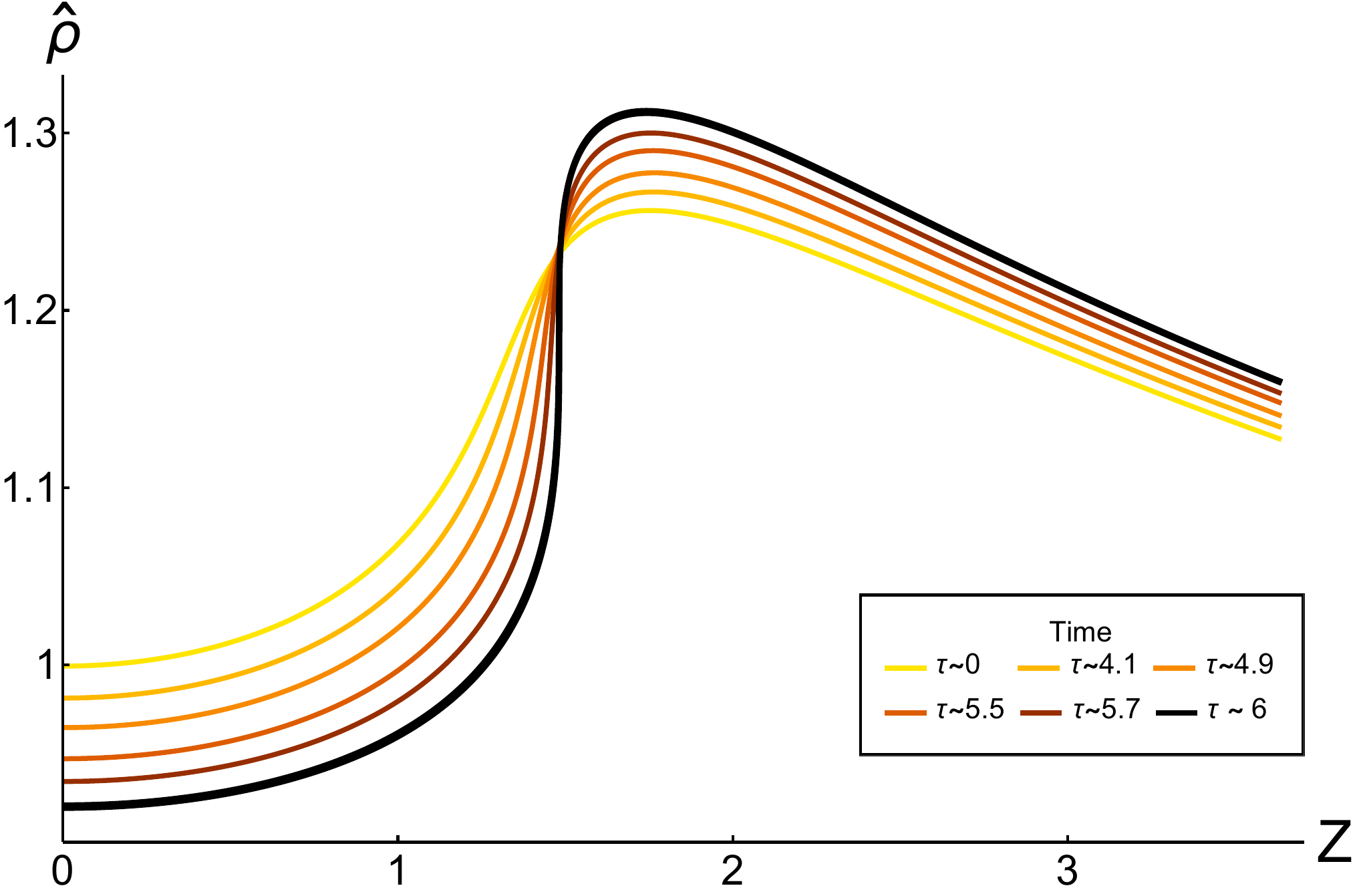}
			\caption{Time evolution of the density}
			\label{fig:Shock_SS_SLM_A}
		\end{subfigure}%
		\begin{subfigure}[b]{0.5\textwidth}
			\centering		\hspace{2cm}\includegraphics[width=7.5cm]{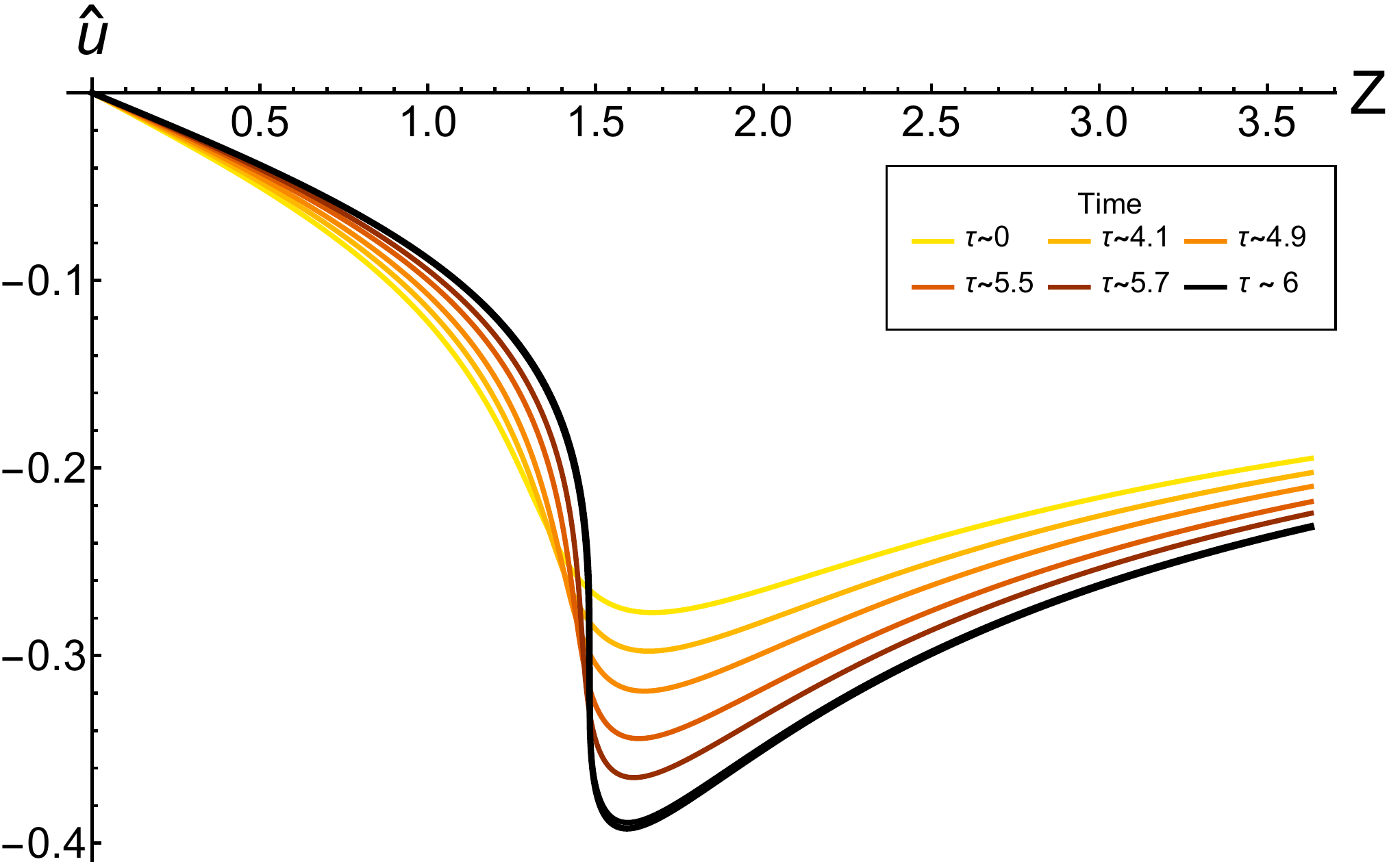}
			\caption{Time evolution of the velocity}
			\label{fig:Shock_SS_SLM_A_U}
		\end{subfigure}%
	
		\begin{subfigure}[b]{0.5\textwidth}
			\centering		\hspace{2cm}\includegraphics[width=7.5cm]{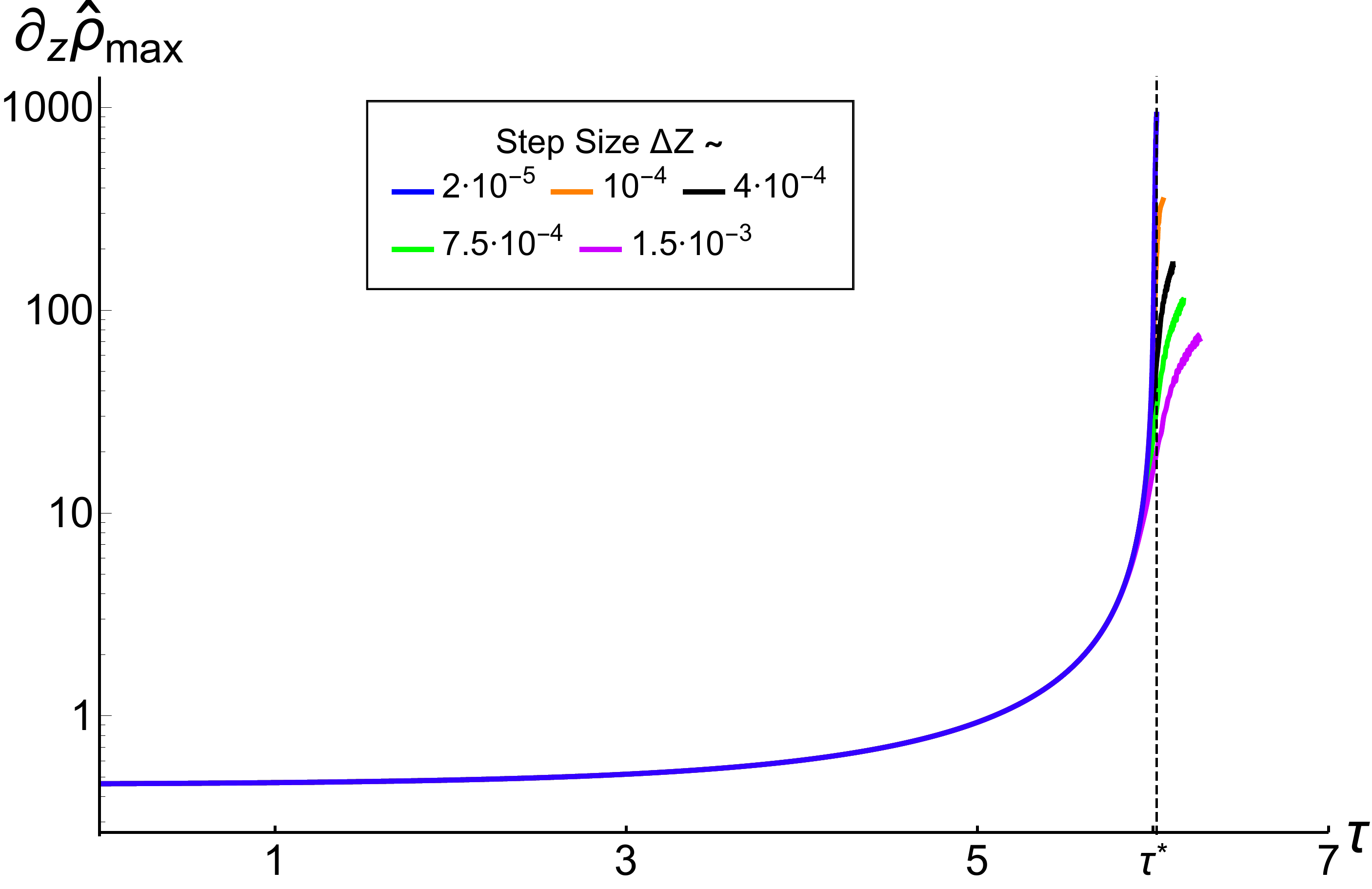}
			\caption{$\max_z\left(\partial_z\hat\rho\right)$}
			\label{fig:Shock_SS_SLM_B}
		\end{subfigure}%
		\begin{subfigure}[b]{0.5\textwidth}
			\centering		\hspace{2cm}\includegraphics[width=7.5cm]{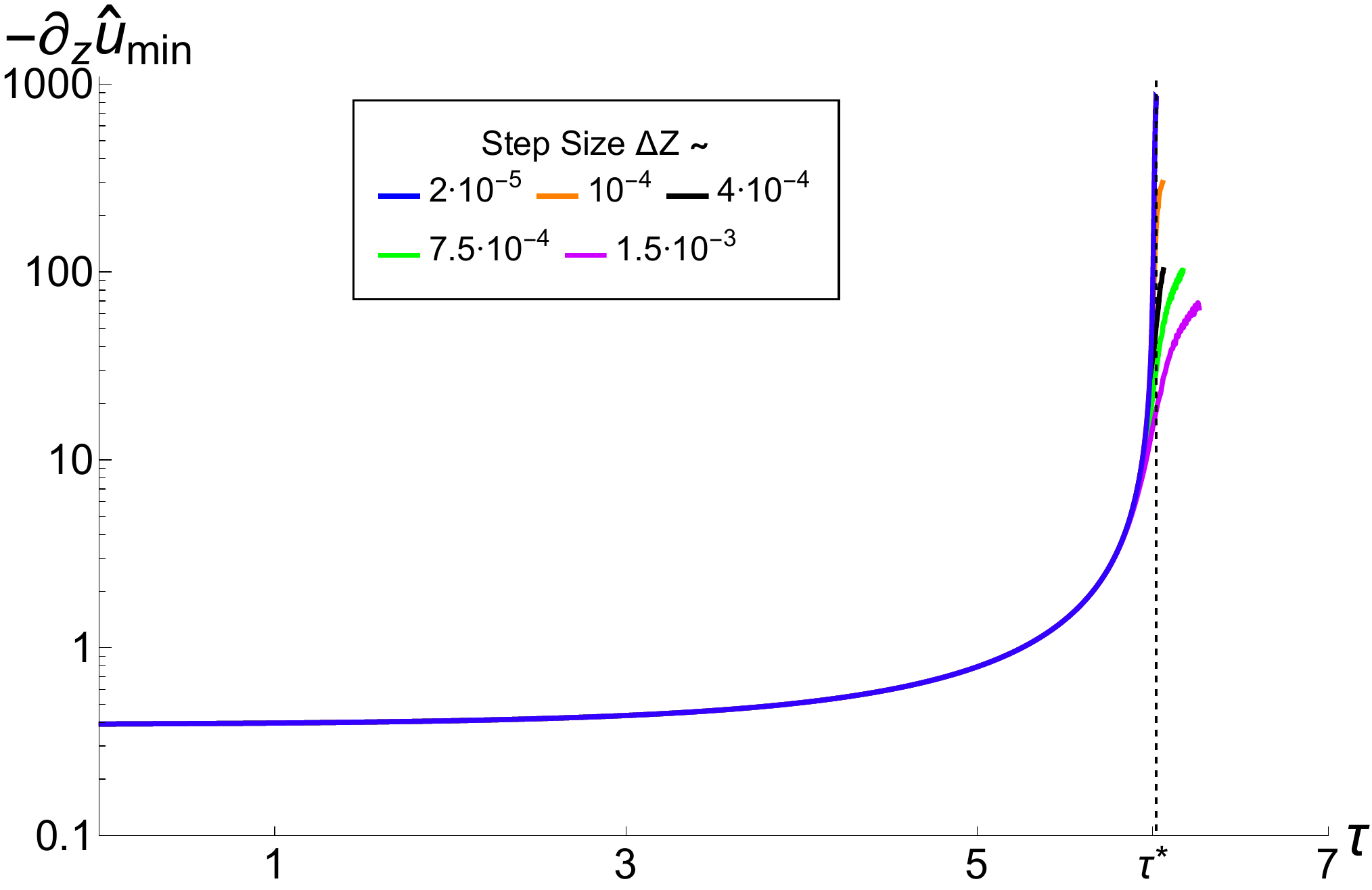}
			\caption{$-\min_z\left(\partial_z\hat{u}\right)$}
			\label{fig:Shock_SS_SLM_B_U}
		\end{subfigure}%
	
		\caption{\small  These plots show the time evolution of a SS perturbed by the first unstable SLM. In (a) and (b) we observe that when $\tau$ approaches $\tau^*\sim 6$, the density $\hat{\rho}$ and the velocity $\hat{u}$ remain finite but develop a shock. In (c) and (d) we see that the maximum of $\partial_Z\hat{\rho}$ and the minimum of $\partial_Z\hat{u}$ rapidly grow when $\tau$ approaches $\tau^*$. Our analysis of this simulation shows that the structure of the formation of the shock is the same as the one presented in the previous section. The specific initial data in these plots is prepared following (\ref{eq:initial_data_blow-up_plus_perturbations}) with  $\epsilon=10^{-3}$, $(d,\ell)=(3,2)$ and $(r_2,\kappa,\LL_1,\theta)$ from table~\ref{table:r_1_full}.}
		\label{fig:Shock_SS_SLM}
	\end{figure}

	After determining common endpoints of generic perturbations, what remains is the exploration of instabilities of SSs triggered by SLMs.	When there is a special element in a continuous family of perturbations, like in this case, one may expect that the properties of this element gradually manifest when we approach to it. For example, if we assume that these perturbations develop a shock in finite time but this element does not, we expect that the time of the shock formation grows when we are close to the special perturbation. For this reason, even if we are not able to exactly construct this perturbation, the exploration of its neighborhood provides useful information. This is our situation because we actually construct NSSs (modes) very close to SSs (SLMs). Recall that once $(d,\ell)$ are fixed SSs (SLMs) are isolated points in a 2-dimensional (4-dimensional) space. In practice, from a numerical point of view, our profiles have more continuous derivatives at $Z_2$ than generic NSSs, but a finite number. For this reason we also explore the endpoints of modes in the neighborhood of SLMs. This process does not show evidence of any transition, concluding that the first unstable SLM ($\LL_1$) leads the system to the formation of a shock (see fig.~\ref{fig:Shock_SS_SLM}) with the same structure as the one presented in the previous section. Amusingly, fine-tuned perturbations of fine-tuned solutions have an ordinary endpoint. 
	
	Our study does not leave significant room from surprises, but we still see some options:
	\begin{itemize}
		\item The existence of a SSs for $r\to1$. Given our analysis of linear modes, if these solutions exist they are the best candidates to be stable.
		
		\item The space of parameters in this problem is huge, five continuous parameters and the dimension. We have explored reasonable values and extracted generic behaviors of generic and fine-tuned solutions. Nevertheless, we cannot say that for a problem that has fine-tuned perturbations for fine-tuned solutions there is no another level of fine-tuning.
		\item It is possible (but unlikely) that the neighborhood of SLMs where some transition is manifestly visible is so narrow that we did not get access to it. This scenario is unlikely because our construction places us very close to the exact solutions relatively to generic profiles and nothing special was observed. However, a more sophisticated method should be developed to reduce errors and get access to regions closer to the exact SS. 
	\end{itemize}
	
	\newpage
	\section{Discussion\vspace{3mm}}\label{sec:Discussions}
	
	This paper shows that the combination of analytical and numerical methods is a very competitive strategy to study the structure of PDEs. We made use of this duo to get a better understanding about the structure of the isotropic implosion of a gas. Almost eighty years after the publication of the pioneering Guderley problem \cite{Guderley}, the authors of \cite{MerleEuler} proved the existence of smooth self-similar solutions to the isentropic compressible Euler equations and, in the present paper, we have developed the first method to construct such smooth profiles. Furthermore, we have designed the search of smooth profiles such that the process by itself provides extra information of the problem. It allowed us to understand how these special solutions arise among non-smooth profiles. We must highlight that our results (in our regime $\nu<8$) are in excellent agreement with \cite{MerleEuler} (in their regime $\nu \to \infty$). It provides strong confidence that the structure that we have reported is actually the general structure of the problem $(1<\nu<\infty)$.

	Regarding the stability of smooth solutions we have studied the spectrum of linear radial perturbations of smooth and non-smooth profiles. Among other results it was shown that smooth solutions are unstable. We adapted our search for smooth profiles to provide the first construction of unstable smooth modes. To conclude our study we explored the endpoints of some instabilities. It was done performing a numerical time evolution of self-similar solutions perturbed by unstable modes. We found that a common endpoint is the classical singularity in the context of the Euler equations, a shock formation; even for smooth solutions under smooth perturbations. The new singularity, the shock out of the origin, has a structure different from the original one, a blow-up of density and velocity at the origin. Therefore, the singularity formation process associated with smooth solutions is unstable. This is of special interest in the context of the energy-supercritical NLS and the compressible NS because in suitable regimes our profiles dominate the dynamics of these systems.
	
	Finally, our work leaves some problems that remain to be addressed:
	\vspace{-0.4cm}
	\begin{itemize}
		\item Our results rely on a semi-numerical construction (in $\nu<8$) where we have obtained intuition about the general structure  of the problem ($1<\nu<\infty$) but it requires confirmation. We hope that the information provided in this paper will be the starting point of future approaches to this problem. Moreover, some of the properties that we have found seem to be accessible analytically.
		\item We find special difficulties in $d=2$. Zeros of $c_+(r)$ are dangerously located in regions where $\nu$ is very close to an integer. It increases the errors in our construction of SSs and SLMs. It motivates the development of a method that provides a better description of these regions.
		\item Now that we know a method to construct smooth solutions and their smooth perturbations, it is time to work on refinements and improvements to get access to higher values of $\nu$.
	\end{itemize}

	{\large\bf Acknowledgements}
	
	I am indebted to Piotr Bizo\'n  for suggesting this problem, as well as for numerous discussions and comments on the manuscript. I am also grateful to Pierre Raphael for discussions, as well as to Oleg Evnin, Brad Cownden and Javier Mas for helpful comments on the manuscript. This work has been supported by the Polish National Science Centre grant number 2017/26/A/ST2/00530.

	
	\appendix
	
	
	\section{Appendix: Regularity of Linear Modes\vspace{3mm}}\label{sec:Appendix_Regularity_Linear_Modes_Z2}
	
	This appendix studies the behavior of linear modes at the origin, $Z_2$ and infinity, determining their level of regularity and bounds $\LL_{\min}<\LL\leq\LL_{\max}$. To analyze these special points of the equations for linear perturbations (\ref{eq:alpha_pert})-(\ref{eq:beta_pert}) we rewrite them using
	\beq
	\alpha(Z) = \frac{\hat{\rho}(Z)}{Z\s(Z)} \tilde{\alpha}(Z).
	\label{eq:alpha_redefinition}
	\eeq
	After several manipulations they take the form of the non-autonomous linear system (viewing $(\sigma,\omega)$ as functions of $Z$)
	\beq
	Z\frac{dv}{dZ} = M v \qquad v:= (\tilde{\alpha},\beta)^{T}
	\label{eq:System_Linear_Modes_Appendix}
	\eeq
	with
	\beq
	 M = \left(\frac{\LL}{\Delta}M_1 + \frac{1}{\Delta^2} M_2 \right),
	\eeq
	\beq
	M_1 = \begin{pmatrix}
		(\omega-1) & \sigma\\
		\sigma & (\omega-1)
	\end{pmatrix},
	\eeq
	\beq
	M_2= \begin{pmatrix}
		\left(\Delta+2\sigma^2\right)\left(\Delta-\frac{\Delta_2}{\sigma}\right) & \sigma\left(-(d+\ell-r(\ell+1)+2\omega)\Delta + 2\Delta_1\right)\\
		-2(1-\omega)\sigma \left(\Delta-\frac{\Delta_2}{\sigma}\right) & -2\ell\sigma\Delta_2 - \left(r-1-(d-1+2\ell)\sigma^2\right)\Delta
	\end{pmatrix}.
	\eeq 
	
	
	\subsection{Expansion of $(\alpha,\beta)$ at $Z_2$}
	
	Linear modes have a singular point on the acoustic cone, it can be seen from system (\ref{eq:System_Linear_Modes_Appendix}) using that $M_2=\Delta_1=\Delta_2=\Delta=0$ at $Z_2$. To study the structure of linear modes in the neighborhood of this point we expand matrix $M$  using that ($\xi := (Z-Z_2)$)
	\beq
	\sigma(\xi)\sim \sigma_2 + \frac{s_1}{Z_2} \xi + ..., \qquad  \omega(\xi)\sim (1-\sigma_2) + \frac{w_1}{Z_2} \xi + ...
	\label{eq:expansion_sigma_w_to_first_order_appendix}
	\eeq
	Then, the leading system is
	\beq
	\xi\frac{dv}{d\xi} = \begin{pmatrix}
		m_1 & -m_2\\
		-m_1 & m_2
	\end{pmatrix}v
	\label{eq:matrix_Z2}
	\eeq
	with
	\begin{align}
	& m_1 = \frac{\LL}{2 (s_1+\omega_1)} + \frac{(2 - d + \ell (r-1) - r) + 2 (-1 + d) \sigma_2}{2\ell(s_1+\omega_1)^2}\omega_1 \label{eq:m1_Z2}\\
	& m_2= \frac{\LL}{2 (s_1+\omega_1)} - \frac{(2 - d + \ell (r-1) - r) + 2 (-1 + d) \sigma_2}{2\ell(s_1+\omega_1)^2}\ell s_1.  \label{eq:m2_Z2}
	\end{align}
	With it we see that $Z_2$ ($\xi = 0$) is a {\em regular singular point}. The eigenvectors of this system of equations are 
	\beq
	\tilde{e}_0 = \left(\frac{m_2}{m_1},1\right)^{T} \qquad \text{and} \qquad \tilde{e}_1 = (-1,1)^{T},
	\eeq
	with their respective eigenvalues $0$ and $\mathcal{N}$ given by
	\beq
	\mathcal{N}(\LL) = \nu + \frac{2}{(\ell+1)(r-1)-(d+1)} \left(\nu+1\right) \LL,
	\label{eq:N_regularity}
	\eeq
	where $\nu$ is the regularity of NSSs introduced in (\ref{eq:expansion_w_in_powers_sigma_Z2}). Details about the derivation of (\ref{eq:N_regularity}) are provided in appendix~\ref{sec:Appendix_Derivations}.

	
	In light of these eigenvalues the structure of linear modes around $Z_2$ is the following. If we consider that $(\s,\w)$ is a SS, matrix M in (\ref{eq:System_Linear_Modes_Appendix}) can be decomposed in integer powers of $\xi$. Then, avoiding values of $\LL$ where $\mathcal{N}$ is integer, to avoid logarithms, the structure around a regular singular point is
	\begin{align}
		& \alpha(\xi) = \underbrace{\sum \tilde{\alpha}_n\xi^n }_{\text{integer powers}} + \underbrace{ c_{\pm}^{(\mathcal{N})} \ |\xi|^{\mathcal{N}}\sum \hat{\alpha}_n \xi^n}_{\text{non-integer powers}}, \label{eq:expansion_alpha_Z2_APPENDIX_v0}\\
		& \beta(\xi) = \underbrace{\sum \tilde{\beta}_n\xi^n }_{\text{integer powers}} + \underbrace{c_{\pm}^{(\mathcal{N})} \ |\xi|^{\mathcal{N}}\sum \hat{\beta}_n \xi^n}_{\text{non-integer powers}}.
		\label{eq:expansion_beta_Z2_APPENDIX_v0}
	\end{align}	
	In case of NSSs, as we explained in the main text, we know that their structure at $Z_2$ is of the form \cite{MerleEuler} 
	\beq
	\sigma(\xi) = \underbrace{\sum \frac{s_n}{Z_2^n}\xi^n }_{\text{integer powers}} + \underbrace{c_{\pm} \ |\xi|^{\mathcal{\nu}}\sum\left(...\right)}_{\text{non-integer powers}} \qquad  \omega(\xi) = \underbrace{\sum \frac{\omega_n}{Z_2^n}\xi^n }_{\text{integer powers}} + \underbrace{c_{\pm} \ |\xi|^{\nu}\sum\left(...\right)}_{\text{non-integer powers}}
	\eeq
	The presence of non-integer terms related to $\xi^{\nu}$ makes that the structure of linear modes given in (\ref{eq:expansion_alpha_Z2_APPENDIX_v0})-(\ref{eq:expansion_beta_Z2_APPENDIX_v0}) is not enough to solve the equations. Plugging the expansion of $(\s,\w)$ in (\ref{eq:System_Linear_Modes_Appendix}) we see that extra powers must be added to compensate terms coming from the non-integer powers initiated by $\xi^{\nu}$ resulting in an structure of the form 
	\begin{align}
		& \alpha(\xi) = \underbrace{\sum \tilde{\alpha}_n\xi^n }_{\text{integer powers}} + \underbrace{c_{\pm}^{(\nu)} \ |\xi|^{\nu-1}\sum\left(...\right) + c_{\pm}^{(\mathcal{N})} \ |\xi|^{\mathcal{N}}\sum\left(...\right)}_{\text{non-integer powers}}, \label{eq:expansion_alpha_Z2_APPENDIX}\\
		& \beta(\xi) = \underbrace{\sum \tilde{\beta}_n\xi^n }_{\text{integer powers}} + \underbrace{c_{\pm}^{(\nu)} \ |\xi|^{\nu-1} \sum\left(...\right) + c_{\pm}^{(\mathcal{N})} \ |\xi|^{\mathcal{N}}\sum\left(...\right)}_{\text{non-integer powers}},
		\label{eq:expansion_beta_Z2_APPENDIX}
	\end{align}
	where $\xi^{\nu-1}$ is present because for linear modes depend on $(\rho',u')$ (as we show below), and $(...)$ gather towers of powers initiated by $\xi^{\nu-1}$ and $\xi^{\mathcal{N}}$ as well as additional powers coming from their combination as $\xi^{\mathcal{N}+\nu-1}$. The important lesson that we learn is that the structure of linear modes have two sources of non-integer powers $c_{\pm}^{(\nu)} \ |\xi|^{\nu-1}$ and $c_{\pm}^{(\mathcal{N})} \ |\xi|^{\mathcal{N}}$. When we work with a SS ($c_{\pm}=0$) only non-integer powers related to $\xi^{\mathcal{N}}$ remain; while for NSS even when we have $c_{\pm}^{(\mathcal{N})}=0$ there are non-integer powers coming from $\xi^{\nu-1}$.
	
	This derivation can be also done from equations (\ref{eq:alpha_pert})-(\ref{eq:beta_pert}) working with $(\alpha,\beta)$. In order to simplify the derivations and expressions we are going to assume that $(\sigma,\omega)$  is a SS ($c_{\pm}=0$). With this simplification integer and non-integer powers can be analyzed separately. First we are going to determine $\mathcal{N}$. In this case we use (\ref{eq:expansion_sigma_w_to_first_order_appendix}) and
	\beq
	\alpha(\xi) \sim  |\xi|^{\mathcal{N}} \sum_{n=0}\hat\alpha_{n}\xi^n, \qquad \beta(\xi) \sim  |\xi|^{\mathcal{N}} \sum_{n=0}\hat\beta_{n}\xi^n.
	\label{eq:expansion_Z2}
	\eeq
	where $\mathcal{N}$ is unknown. Plugging this ansatz into equations (\ref{eq:alpha_pert})-(\ref{eq:beta_pert})  for linear perturbations, removing the factors $\xi^{\mathcal{N}}$ and gathering terms with power $k$, yields
	\begin{align}
	&(1+k+\mathcal{N})\left(\hat\alpha_{k+1}+ \left(\frac{\ell}{2}\right)^\frac{\ell}{2}\left(Z_2 \sigma_2\right)^{\ell-1}\hat\beta_{k+1}\right)+c_1\hat\alpha_{k} + c_2 \hat\beta_{k} + P_k(\LL,\mathcal{N}, \hat\alpha_{i<k-1},\hat\beta_{i<k-1}) = 0, \label{eq:alpha_expans_Z2}\\
	&(1+k+\mathcal{N})\left(\hat\alpha_{k+1}+ \left(\frac{\ell}{2}\right)^\frac{\ell}{2}\left(Z_2 \sigma_2\right)^{\ell-1}\hat\beta_{k+1}\right)+\tilde{c}_1\hat\alpha_{k} + \tilde{c}_2 \hat\beta_{k} + \tilde{P}_k(\LL,\mathcal{N}, \hat\alpha_{i<k-1},\hat\beta_{i<k-1}) = 0, \label{eq:beta_expans_Z2}
	\end{align}
	where $P_k$ and $\tilde{P}_k$ are linear in $\hat\alpha_{i}$ and $\hat\beta_{i}$, with regular coefficients and no independent term. Coefficients in (\ref{eq:alpha_expans_Z2})-(\ref{eq:beta_expans_Z2}) are
	\begin{align}
	& c_1 = (r-1)\ell+\LL-d+(d-1)\sigma_2+(k+\mathcal{N}+1)(\sigma_2-\w_1),\\ 
	&c_2 = \left(\frac{\ell}{2}\right)^{\frac{\ell}{2}}\left(Z_2\sigma_2\right)^{l-1}\left((d-1)\sigma_2+ (1+k+\mathcal{N})\ell (\sigma_2 + s_1)\right),\\
	& \tilde{c}_1 =(1+k+\mathcal{N})\left(2-\ell\right)(\sigma_2+s_1),\\
	&\tilde{c}_2 = \left(\frac{\ell}{2}\right)^{\frac{\ell}{2}}\left(Z_2\sigma_2\right)^{l-1}\left(r-2+\LL + (1+k+\mathcal{N}) (\sigma_2 - \w_1) \right).
	\end{align}
	From (\ref{eq:alpha_expans_Z2})-(\ref{eq:beta_expans_Z2}) we obtain an equation for $(\hat{\alpha}_k,\hat{\beta}_k)$ and an equation for $(\hat{\alpha}_{k+1},\hat{\beta}_{k+1})$. We see that for $\mathcal{N}$ given in (\ref{eq:N_regularity}) and $k=0$ the equation for $(\hat{\alpha}_0,\hat{\beta}_0)$ vanishes; then, after relating $\hat{\beta}_0$ to $\hat{\alpha}_0$ using the equation for $k=-1$ and the scaling symmetry, we see that there is a single freedom $\hat{\alpha}_0=0,1$. With it, from $(\LL,\hat{\alpha}_0)$ we can perform  an iterative resolution of $(\hat\alpha_k,\hat\beta_k)$. Furthermore, we can see that if $\hat{\alpha}_0 = 0$, all the coefficients vanish and then the expansions in (\ref{eq:beta_expans_Z2}) are trivially zero. This is a key point in our construction of SLMs, if we find a trajectory where the coefficient of $\xi^{\mathcal{N}}$ vanishes, then, the tower of non-integer powers associated with $\mathcal{N}$ is not present. This result works independently for trajectories coming from the left and from the right of $Z_2$; therefore, a smooth connection at $Z_2$ requires that the coefficient of $\xi^{\mathcal{N}}$ vanishes for both. Additionally, generic trajectories have a finite regularity given by $\mathcal{N}(\LL)$. It imposes a bound for $\LL$ given by the minimum level of regularity $\mathcal{N}(\LL_{max})=1$,
	\beq
	\LL_{max} = - \frac{(\ell+1)(r-1)-(d+1)}{2}\left(\frac{\nu - 1}{\nu + 1}\right).
	\eeq
	In our range of parameters, (\ref{eq:parameters}), $\LL_{max}>0$, it can be shown using that $\nu=\lambda_-/\lambda_+$,\\ $\lambda_-<\lambda_+<0$, and 
	\beq
	\frac{(\ell+1)(r-1)-(d+1)}{2} < 0.
	\eeq

	The analysis of the tower of integer powers is simple, we only need to fix $\mathcal{N}=0$ in (\ref{eq:alpha_expans_Z2})-(\ref{eq:beta_expans_Z2}). In this case $k=-1$ is the equation that vanishes providing the freedom $\hat{\alpha}_0=0,1$ (after using the scaling symmetry). Then, an iterative construction can be performed to obtain $(\hat{\alpha}_k,\hat{\beta}_k)$ (recall that we avoid $\Omega$ such that $\mathcal{N}$ is integer). With this process we see that $\hat{\alpha}_0=0$ leads to expansions that are trivially zero at $Z_2$. Hence, a nontrivial construction relies on $\hat{\alpha}_0=1$. 
	
	Additionally, extracting contributions $(s_{k+1},\omega_{k+1})$ from (\ref{eq:alpha_expans_Z2})-(\ref{eq:beta_expans_Z2}) (with $\mathcal{N}=0$) we get
	\begin{align}
	&(1+k)\left(\hat\alpha_{k+1}+ \left(\frac{\ell}{2}\right)^\frac{\ell}{2}\left(Z_2 \sigma_2\right)^{\ell-1}\hat\beta_{k+1}\right)+c_1\hat\alpha_{k} + c_2 \hat\beta_{k} + Q_k \nonumber\\
	&+\frac{(k+1)}{Z_2\sigma_2}\left(\left(\frac{\ell}{2}\right)^{\frac{\ell}{2}}\left(Z_2\sigma_2\right)^{\ell-1}\ell \hat{\beta}_0 s_{k+1} - \hat{\alpha}_0 \omega_{k+1}\right) = 0, \label{eq:alpha_expans_Z2_1-modes}\\
	&(1+k)\left(\hat\alpha_{k+1}+ \left(\frac{\ell}{2}\right)^\frac{\ell}{2}\left(Z_2 \sigma_2\right)^{\ell-1}\hat\beta_{k+1}\right)+\tilde{c}_1\hat\alpha_{k} + \tilde{c}_2 \hat\beta_{k} + \tilde{Q}_k \nonumber\\
	& +\frac{(k+1)}{Z_2\sigma_2}\left((2-\ell)\hat{\alpha}_0 s_{k+1}-\left(\frac{\ell}{2}\right)^{\frac{\ell}{2}}\left(Z_2\sigma_2\right)^{\ell-1} \hat{\beta}_0 \omega_{k+1}\right) = 0, \label{eq:beta_expans_Z2_1-modes}
	\end{align}
	where $Q_k$ and $\tilde{Q}_k$ are functions $P_k$ and $\tilde{P}_k$ in (\ref{eq:alpha_expans_Z2})-(\ref{eq:beta_expans_Z2}) after extracting all terms with $s_{k+1}$ and $\omega_{k+1}$. From these expansions we can see that coefficients $\hat{\alpha}_{k}$ and $\hat{\beta}_{k}$ have contributions from $s_{k+1}$ and $\omega_{k+1}$ because the subtraction of (\ref{eq:alpha_expans_Z2_1-modes})-(\ref{eq:beta_expans_Z2_1-modes}) that cancels $\hat{\alpha}_{k+1}$ and $\hat{\beta}_{k+1}$ does not remove  $s_{k+1}$ and $\omega_{k+1}$. Therefore, when we consider NSSs these contributions are the cause that $\xi^{\nu-1}$ is the lowest non-integer power (in this tower) instead of $\xi^{\nu}$.
	
	\vspace{0.2cm}
	\textbf{Remarks}
	\vspace{-0.2cm}
	\begin{enumerate}
		\item There are two independent trajectories emerging from $Z_2$. One of them vanishes at this point and goes like $v\sim\xi^\mathcal{N}$. The second trajectory does not vanish at $Z_2$ and does not have terms related to $\xi^{\mathcal{N}}$. Generic trajectories are the linear combination of both of them and have the structure given in (\ref{eq:expansion_alpha_Z2_APPENDIX})-(\ref{eq:expansion_beta_Z2_APPENDIX}).
		\item 0-modes vanish at $Z_2$, then, they correspond to the trajectory $v\sim\xi^{\mathcal{N}}$ and their regularity is determined by $\mathcal{N}(\LL)$.
		\item $\mathcal{N}(\LL)$ imposes a positive upper bound $\LL_{max}>0$. Meaning that blow-up profiles are always linearly unstable under some 0-modes.
		\item 0-modes with $\LL=0$ have the regularity of blow-up profiles, $\mathcal{N}(0) = \nu$.
		\item Unstable 0-modes, $\LL>0$, have less regularity than blow-up profiles.
		\item Generic 1-modes are associated with generic trajectories (\ref{eq:expansion_alpha_Z2_APPENDIX})-(\ref{eq:expansion_beta_Z2_APPENDIX}) and have regularity given by the minimum of $\mathcal{N}(\LL)$ and $\nu-1$. However, for fine-tuned values of the parameters there are 1-modes associated with the trajectory (on both sides) that does not have contributions from $\xi^{\mathcal{N}}$; they have $c_{\pm}^{(\mathcal{N})}=0$. Then, their regularity is $\nu-1$.
		\item For 1-modes associated with SSs the term $c_{\pm}^{(\nu)}\xi^{\nu-1}$ is not present and their regularity is given by $\mathcal{N}$. For fine-tuned situations $c_{\pm}^{(\mathcal{N})}=0$ and then these 1-modes are smooth.
		
	\end{enumerate}
	
	\subsection{Expansion of $(\alpha,\beta)$ at the origin}
	
	To extract the structure of linear modes at the origin we use the expansion of blow-up profiles at this point \cite{MerleEuler}
	\beq
	\sigma(Z)\sim \frac{s_{-1}}{Z} + \sum_{n=0}s_{2n+1}Z^{2n+1} \qquad \omega(Z) \sim \sum_{n=0}\omega_{2n}Z^{2n}.
	\eeq
	Expanding matrix $M$ at the origin we see that the leading term of system (\ref{eq:System_Linear_Modes_Appendix}) takes the form
	\beq
		Z\frac{dv}{dZ} = \begin{pmatrix}
			0 & 0\\
			0 & (1-d)
		\end{pmatrix} v.
		\label{eq:Linear_system_leading_term_Origin_APPENDIX}
	\eeq
	Then, regular trajectories emerging from the origin are initiated by $\beta_0=0$ and $\tilde{\alpha}_0=0,1$ (using the scaling symmetry). 
	This conclusion can be also obtained from $(\alpha,\beta)$ and (\ref{eq:alpha_pert})-(\ref{eq:beta_pert}), using
	\beq
	\hat\rho(Z)\sim \sum_{n=0}\hat{\rho}_{2n}Z^{2n} \qquad \hat{u}(Z) \sim \sum_{n=0}\hat{u}_{2n+1}Z^{2n+1},
	\eeq
	and an ansatz
	\beq
	\alpha(Z)\sim \sum\alpha_{n}Z^{n} \qquad \beta(Z) \sim \sum\beta_{n}Z^{n}.
	\label{eq:expansion_alpha_beta_Origin_APPENDIX}
	\eeq
	Plugging them into equations (\ref{eq:alpha_pert})-(\ref{eq:beta_pert}) for linear perturbations and gathering terms with power $k$ we get
	\beq
	\hat\rho_0 (d+k-1)\beta_{k} = P_k(\LL,\alpha_{i< k},\beta_{i< k}) \qquad (\gamma-1)\hat{\rho}_0^{\gamma-2}k\alpha_{k}= \tilde{P}_k(\LL,\alpha_{i< k},\beta_{i< k})
	\eeq
	where $P_k$ and $\tilde{P}_k$ are linear in $\alpha_i$ and $\beta_i$, their constant coefficients remain regular for any $\LL$ and $k$ and there are no independent coefficients. Hence, if we impose $\alpha_{n<0}=\beta_{n<0}=0$ and after obtaining $\beta_0=0$ when $k=0$, we can follow the constructive sequence
	\beq
	(\alpha_0,\LL) \to (\alpha_1,\beta_1) \to (\alpha_2,\beta_2) \to ... \to (\alpha_k, \beta_{k}) \to ...
	\label{eq:constructive_sequence_Z_0}
	\eeq
	to calculate any coefficient; showing that ansatz (\ref{eq:expansion_alpha_beta_Origin_APPENDIX}) is consistent. We can easily show that $\beta_{2k}=\alpha_{2k+1}=0$ in agreement with the parity of $\rho$ and $u$. Once $\LL$ is fixed the solution that emerges from the origin only depends on $\alpha_0$ and after rescaling just two possibilities remain $\alpha_0 = 1$ and $\alpha_0 = 0$. The constructive sequence (\ref{eq:constructive_sequence_Z_0}) leads to a nontrivial profile in the first case and the trivial one ($\alpha=\beta=0$) in the second case. If we repeat the process for trajectories $v\sim Z^{1-d}$, we see that ansatz (\ref{eq:expansion_alpha_beta_Origin_APPENDIX}) does not solve the equation for $k=0$ and logarithms must be included, in agreement with the fact that eigenvalues of (\ref{eq:Linear_system_leading_term_Origin_APPENDIX}) differ by an integer. Despite this fact, these trajectories are not relevant in our problem because they diverge at $Z=0$
	
	From this process we have learned that two regular trajectories emerge from the origin, the trivial one, $\alpha=\beta=0$, and a nontrivial trajectory with $\alpha(0)\neq0$ that follows ansatz (\ref{eq:expansion_alpha_beta_Origin_APPENDIX}). Then, given that there are no singular points in $(0,Z_2)$, $\alpha(0)=0$ leads to $\alpha=\beta=0$ in $[0,Z_2)$ while the trajectory with $\alpha(0)\neq0$ generically reaches $Z_2$ with the structure (\ref{eq:expansion_alpha_Z2_APPENDIX})-(\ref{eq:expansion_beta_Z2_APPENDIX}). However, for fine-tuned values of the parameters $\alpha(0)\neq 0$ may be the curve with $c_+^{(\mathcal{N})}=0$.


		\subsection{Expansion of $(\alpha,\beta)$ at infinity}
	
	Now we study the structure of linear modes at infinity. For this purpose we use the expansion of blow-up profiles \cite{MerleEuler}
	\beq
	\s(Z\sim\infty) \sim \sum_{n=1}^{\infty}\frac{s_n}{Z^{rn}} \qquad \w(Z\sim\infty)\sim \sum_{n=1}^{\infty}\frac{\w_n}{Z^{rn}}
	\label{eq:expansion_boundary_SW}
	\eeq
	to get the leading term of the linear system (\ref{eq:System_Linear_Modes_Appendix}) ($y:=1/Z^r$)
	\beq
	y\frac{dv}{dy} = \begin{pmatrix}
		\frac{\LL+r-1}{r} & 0\\
		0 & \frac{\LL+r-1}{r}
	\end{pmatrix} v.
	\eeq
	The eigenvalues of the leading matrix are equal but the space of eigenvectors is not degenerate. Then, this structure provides a freedom materialized in the expansion at infinity
	\beq
	\alpha(Z)\sim \frac{\theta\beta_0}{Z^{(r-1)\ell+\LL}}+..., \qquad \beta(Z)\sim \frac{\beta_0}{Z^{r-1+\LL}}+...
	\eeq
	Here we have introduced $\alpha_0 = \theta \beta_0$ to have a quantity, $\theta$, invariant under the scaling symmetry represented by $\beta_0$.  (Note that from (\ref{eq:alpha_redefinition}) one has $\theta = \left(\ell/2\right)^{\ell/2}\eta^{\ell-1} \hat{\alpha}_0/\hat{\beta}_0$). In order to have linear modes that remain regular at infinity we obtain the bound
	\beq
	\LL \geq \LL_{\min} = (1-r)\min{(\ell,1)}.
	\eeq
	When $\LL>\LL_{\min}$ linear modes vanish at infinity, preserving the behavior of the self-similar solutions.
	
	From this process we have learned that a one-parameter family of trajectories emerges from the infinity. There are no singular points in $(Z_2,\infty)$, therefore, all these trajectories must end at $Z_2$. Among them we have one that reaches $Z_2$ with $c_{-}^{(\mathcal{N})}=0$ and one that vanishes at this point. Any other curve has the generic structure (\ref{eq:expansion_alpha_Z2_APPENDIX})-(\ref{eq:expansion_beta_Z2_APPENDIX}).
	
	\subsection{Derivation of $\mathcal{N}$}\label{sec:Appendix_Derivations}
	This appendix provides details about the derivation of expression (\ref{eq:N_regularity}), that we reproduce here
		\beq
	\mathcal{N}(\LL) = \nu + \frac{2}{(\ell+1)(r-1)-(d+1)} \left(\nu+1\right) \LL.
	\label{eq:N_regularity_second_appendix}
	\eeq
	This eigenvalue comes from the diagonalization of the matrix given in (\ref{eq:matrix_Z2}) and its components $m_1$ and $m_2$ in (\ref{eq:m1_Z2})-(\ref{eq:m2_Z2}) 
	\beq
		\mathcal{N} = m_1 + m_2.
	\eeq
	We reorganize terms in $m_1$ and $m_2$ to split the eigenvalue into two new terms $\mathcal{N} = T_1 + T_2$ as follows 
	\beqa
		&T_1 = \frac{\LL}{(s_1+\omega_1)}, \label{eq:T1}\\
		&T_2 =  \frac{(2 - d + \ell (r-1) - r) + 2 (-1 + d) \sigma_2}{2\ell(s_1+\omega_1)^2}(\ell s_1 + \omega_1). \label{eq:T2}
	\eeqa
	By direct comparison with (\ref{eq:N_regularity_second_appendix}) we see that we have to obtain the following expressions
	\beq
		\left(s_1+\w_1\right) = \frac{(\ell+1)(r-1)-(d+1)}{2}\left(\nu+1\right)^{-1}, \qquad \text{and} \qquad T_2 = \nu.
		\label{eq:expressions}
	\eeq
	First we are going to show the derivation of the first expression. To this end we use the coefficients defined in (\ref{eq:c_i_coefficients}) and reproduced here
	\beq
	 \begin{cases}
		c_1 = \partial_\w \Delta_1\big{|}_{P_2}, &	c_3 = \partial_\w \Delta_1\big{|}_{P_2},\\ 
		c_2 = \partial_\s \Delta_2\big{|}_{P_2}, & c_4 = \partial_\s \Delta_2\big{|}_{P_2}.
	\end{cases}
	\eeq
	Now we expand $(\s,\w)$ at first order around $Z_2$ ($\xi := Z-Z_2$)
	\beq
		\s(\xi)\sim \s_2 + \frac{s_1}{Z_2} \xi + ...
		\qquad
		 \w(\xi)\sim (1-\s_2) + \frac{w_1}{Z_2} \xi + ... 
	\eeq
	and plug these expressions in the system of equations for $(\s(Z),\w(Z))$, (\ref{eq:System_W_S}). Gathering terms with the same power of $\xi$ we see that the coefficient for $\xi^{-1}$ is the quadratic equation that determines $\s_2$
	\beq
		(d-1)\s^2 - (d-1-(\ell-1)(r-1))\s + (r-1) = 0
		\label{eq:quadratic_eq_S2}
	\eeq
	Power $\xi^0$ provides the relations of our interest
	\beq
	\omega_1 + \frac{c_1 w_1 + c_3 s_1}{-2\s_2 (w_1 + s_1)} = 0, \qquad s_1 + \frac{c_2 w_1 + c_4 s_1}{-2\s_2 (w_1 + s_1)} = 0.
	\label{eq:s1_w1_equations}
	\eeq
	Manipulating the first equation we obtain the following expressions
	\beq
	s_1 = - w_1 \frac{(c_1 - 2 \s_2 w_1)}{c_3 - 2 \s_2 w_1}, \qquad	s_1+w_1 = w_1\frac{c_3-c_1}{(c_3-2\s_2 w_1)}
	\label{eq:s1_s1pw1}
	\eeq
	and substituting them into the second equation the result is a quadratic equation for $w_1$ with solutions
	\beq
	w_1^{(\pm)} = \frac{(c_1^2-c_1c_3+2c_2c_3-c_1c_4-c_3c_4)\pm(c_1-c_3)\sqrt{(c_1-c_4)^2+4c_2 c_3}}{4\s_2\left(c_1+c_2-c_3-c_4\right)}.
	\label{eq:w1_solution}
	\eeq
	Therefore, substituting (\ref{eq:s1_s1pw1})-(\ref{eq:w1_solution}) in the second equation of (\ref{eq:s1_w1_equations}) we obtain 
	\beq
	s_1+w_1^{(\pm)}=\frac{c_1+c_4 \pm \sqrt{(c_1-c_4)^2 + 4 c_2 c_3}}{4\s_2} = \frac{\lambda_{\pm}}{2\s_2}
	\label{eq:s1_p_w1_v2}
	\eeq
	where $\lambda_{\pm}$ are the eigenvalues of the Jacobian matrix given in (\ref{eq:lambda_eigenvalues}) and replicated here
	\beq
		\lambda_{\pm} = \frac{c_1+c_4\pm\sqrt{(c_1-c_4)^2+4c_2c_3}}{2}.
	\eeq
	Now we are going to write $\s_2$ in terms of $\lambda_\pm$ using that $\lambda_++\lambda_- = c_1+c_4$ and that
	\beq
	c_1 + c_4 - ((\ell+1)(r-1)-(d+1))\s_2 = - \frac{(\ell-1)}{\ell}\left((d-1)\s_2^2 - (d-1-(\ell-1)(r-1))\s_2+(r-1)\right).
	\eeq
	The RHS of this expression is proportional to the quadratic equation (\ref{eq:quadratic_eq_S2}) solved by $\s_2$ yielding 
	\beq
		\lambda_+ + \lambda_- = ((\ell+1)(r-1)-(d+1))\s_2.
	\eeq
	From this last step (\ref{eq:s1_p_w1_v2}) takes the form
	\beq
		s_1+w_1^{(\pm)} = \frac{((\ell+1)(r-1)-(d+1))}{2}\frac{\lambda_{\pm}}{\lambda_+ + \lambda_-}
		\label{eq:s1_p_w1_v3}
	\eeq
	Finally, each sign ($\w_1^{(\pm)}$) is associated with one of the two trajectories that end at $(\s_2,1-\s_2)$ in the phase portrait. We have explained that these trajectories ($\w(\s)$) have slopes $\w_{\pm}$ given in (\ref{eq:slope_w_pm}). Using that $w_\pm = w_1/s_1$ and (\ref{eq:s1_s1pw1}) we can see that $w_-$ is associated with $w_1^{(+)}$; therefore the first expression of our interest in (\ref{eq:expressions}) is obtained (recall that $\nu := \lambda_-/\lambda_+$)
	\beq
		\left(s_1+w_1\right) = \frac{((\ell+1)(r-1)-(d+1))}{2}\left(\nu + 1\right)^{-1}, \quad  T_1 = \frac{2}{(\ell+1)(r-1)-(d+1)} \left(\nu+1\right) \LL.
		\label{eq:s1_p_w1_v4}
	\eeq
	
	To derive the second expression in (\ref{eq:expressions}), $T_2 = \nu$, we do not know a clever strategy; instead of that we made use of a brute force calculation of $T_2-\nu$. First we manipulated $T_2-\nu$ a bit using (\ref{eq:s1_p_w1_v4}) and substituted $\s_2$, $\w_1$ and $\nu$ by its expressions in terms of the parameters $(d,\ell,r)$ of the problem that can be obtained from (\ref{eq:quadratic_eq_S2}), (\ref{eq:w1_solution}), (\ref{eq:lambda_eigenvalues}). Then, we have to perform a large number of operations in an equation that only depends on $(d,\ell,r)$. We encourage any reader that wants to face this task to make use of a program that performs symbolic calculation (this is the simplest and fastest way that we have found to show that $T_2-\nu = 0$). At this point it is worth to mention our heuristic motivation to originally perform such a brute force calculation hoping that $T_2-\nu = 0$. This heuristic observation came from the fact that after obtaining the expression for $T_1$ in (\ref{eq:s1_p_w1_v4}) we know that $\mathcal{N}(\LL = 0) = T_2$. Therefore, the $0$-mode with $\LL=0$ has regularity $T_2$ at $Z_2$. Roughly speaking this mode can be understood as the difference between two self-similar solution associated with the same parameters $(d,\ell,r)$ but different $\kappa$ that are very close ($\Delta\kappa \ll 1$). For this reason, we expected that the regularity of the 0-mode with $\LL$ is equal to the regularity of these solutions; namely, $\nu$. This led to the conjecture that $T_2 = \nu$ and the confirmation by the brute force calculation described above.


	\section{Appendix: Numerical Methods\vspace{3mm}}\label{sec:Appendix_Numerical_Methods}
	
	\subsection{Construction of Blow-Up Profiles and Linear Perturbations}
	\label{sec:Appendix_Spectral_Methods}
	
	For the construction of blow-up profiles we followed several strategies mainly based on spectral methods\footnote{Similar techniques were used in \cite{BiasiFloquet,BiasiFloquet2} to construct time-periodic solutions in anti-de Sitter space-time.}. In this section we describe the most successful procedure. Our different approaches to this problem show that the direct resolution of $(\s(Z),\w(Z))$ is more accurate than the construction of $\w(\s)$ and after that $\s(Z)$ or the direct resolution of $(\hat{\rho},\hat{u})$. The main challenges that one must face are:
	\begin{itemize}
		\item $(Z\ll Z_2)$: $\s\to\infty$ when $Z\to0$.
		\item $(Z> Z_2)$: $(\s,\w)$ have a relatively sharp structure at $Z\gtrsim Z_2$.
		\item $(Z\gg Z_2)$: $(\s,\w)$ have a polynomial decay at infinity in powers of $Z^r$.
	\end{itemize}
	To deal with these difficulties the $Z$-coordinate is divided into three intervals
	\beq
	I_{1} := [0,Z_{p1}], \qquad I_{2} := [Z_{p1}, Z_{p2}], \qquad I_{3} := [Z_{p2},\infty),
	\eeq
	with $0<Z_{p1}<Z_{p2}<\infty$ (we usually work with $Z_{p1}=Z_2$). In order to remove singularities at the origin $\s(Z)$ is redefined in $I_1$ making use of its expansion at the origin
	\beq
	\s(Z) = \frac{\delta}{Z}+ s_1 Z + \tilde{\s}(Z) Z^3,
	\eeq
	where $\delta$ is a free parameter, $\tilde{\s}(Z)$ is the new unknown that satisfies $\tilde{\s}(0)=s_3$ and $s_{1},s_3$ come from the expansion of $\s(Z)$ at the origin. In $I_2$ we work with the standard unknowns $(\s,\w)$ but introduce a change of coordinate $Z(X)$, similar to (\ref{eq:ZX_coordinates}), to increase the density of the numerical grid around a specific point. In $I_3$, setting $Z_{p2}\gg Z_2$, $(\s,\w)$ are substituted by their asymptotic expansions at infinity. 
	
	At the numerical level, we use a pseudospectral method to construct $(\s,\w)$ in $I_{1,2}$. These intervals are discretized using the Gauss-Lobatto collocation points of Chebyshev polynomials \cite{BoydBook} (also denominated Chebyshev points of the second kind)
	\beq
	x_{n}^{(i)}=\cos\left(\pi\frac{n}{N_{i}}\right), \qquad n = 0,...,N_{i}
	\label{eq:collocation_points_x}
	\eeq
	which transferred to our intervals give
	\beq
	I_{1}\to Z_{n}^{(1)} = \frac{Z_{p1}}{2}(1+x_n^{(1)}), \qquad I_{2}\to Z_{n}^{(2)}(x_n).
	\eeq
	The values of our unknown functions on these grids are
	\beq
	\left(\s_{n}^{(i)},\w_{n}^{(i)}\right) := \left(\s(Z_{n}^{(i)}),\w(Z_{n}^{(i)})\right), \qquad \left(\s^{(3)},\w^{(3)}\right) := \left(\sum_{k=1}^{k_{cut}} \frac{\hat s_k(\kappa,\eta)}{\left(Z_{p2}\right)^{rk}}, \sum_{k=1}^{k_{cut}} \frac{\hat\w_k(\kappa,\eta)}{\left(Z_{p2}\right)^{rk}}\right)
	\eeq
	where the last equation represents the first $k_{cut}$ terms in the expansions of $(\s,\w)$ at infinity evaluated at $Z_{p2}$. We choose the set (\ref{eq:collocation_points_x}) of collocation points because they contain the limits of the intervals (note that $Z_{N_1}^{(1)}=0,\ Z_0^{(1)}=Z_{N_2}^{(2)}=Z_{p1}$ and $Z_{0}^{(2)}=Z_{p2}$), easing the imposition of boundary conditions
	\beq
	\left(\tilde{\s}_{N_1}^{(1)},\w_{N_1}^{(1)}\right) = (s_3, \w_e), \quad \left(\s_{0}^{(1)},\w_{0}^{(1)}\right) = \left(\s_{N_2}^{(2)},\w_{N_2}^{(2)}\right), \quad \left(\s_{0}^{(2)},\w_{0}^{(2)}\right) = \left(\s^{(3)},\w^{(3)}\right).
	\label{eq:boundary_conditions_Spectral}
	\eeq
	The discrete differentiation matrices $D^{(k)}$ that allow us to calculate derivatives at the collocation points $x_n$, namely,
	\beq
	d_x^{k}f_i:=\frac{d^{k}f(x)}{dx^k}\big{|}_{x=x_i} = \sum_{j=0}^{N}D_{ij}^{(k)}f_j 
	\eeq
	have the form \cite{BoydBook}
	\beq
	D_{ij}^{(1)} =
	\begin{cases}
		\frac{2N^2+1}{6} & \text{if } i=j=0\\
		-\frac{2N^2+1}{6} & \text{if } i=j=N\\
		-\frac{x_{j}}{2(1-x_j)^2} & \text{if } i=j\neq 0,N\\
		\frac{c_{i}}{c_j}\frac{(-1)^{i+j}}{x_i-x_j} & \text{if } i\neq j
	\end{cases}
	\qquad D^{(k)} = \left(D^{(1)}\right)^{k}
	\eeq
	with $c_{0}=c_N=2$ and $c_{i}=1$ otherwise. Differentiation matrices on $Z_{n}^{(i)}$ follow from the chain rule. Hence, the discrete version of system (\ref{eq:System_W_S}) is an algebraic nonlinear system of $2N_i+2$ equations
	\beq
	Z_{n}^{(i)}d_z\w_n^{(i)} = - \frac{\Delta_1(\w_n^{(i)},\s_n^{(i)})}{\Delta(\w_n^{(i)},\s_n^{(i)})}, \qquad Z_{n}^{(i)}d_z\s_n^{(i)} = - \frac{\Delta_2(\w_n^{(i)},\s_n^{(i)})}{\Delta(\w_n^{(i)},\s_n^{(i)})},
	\label{eq:discrete_system_W_S}
	\eeq
	with unknowns $(\delta,\tilde{\s}_n^{(1)},\w_n^{(1)})$ in $I_1$ ($2N_1+3$) and $(\s_n^{(2)},\w_n^{(2)})$ in $I_2$ ($2N_2+2$).
	To solve these nonlinear problem we replace some equations by boundary conditions, fix one of our unknowns in $I_1$, and apply the Newton-Raphson method (working with numerical precision higher than standard double-numbers). This method requires an initial seed that must be quite close to the solution to guarantee the convergence to this solution. Once $(\s,\w)$ are constructed $(\hat\rho,\hat{u})$ follow from (\ref{eq:Emden_transform}) and we apply the same discretization to obtain the linear modes $(\LL,\alpha,\beta)$. Finally, all these profiles are transferred to the continuous variable Z using their expansions in $I_3$ and the barycentric interpolation \cite{Barycentric} in $I_{1,2}$
	\beq
	P_{f^{(i)}}(Z) = \frac{\sum\limits_{n=0}^{N_i}\frac{c_n}{Z-Z_n^{(i)}}f_n^{(i)}}{\sum\limits_{n=0}^{N_i}\frac{c_n}{Z-Z_n^{(i)}}} \qquad c_0=\frac{1}{2},\ c_{N_i}=\frac{(-1)^{N_i}}{2},\ c_{n\neq 0,N_i} = (-1)^{n}.
	\label{eq:barycentric_interpolation}
	\eeq
	
	
	\subsection{Numerical Search of Smooth Solutions}
	\label{sec:Appendix_Search_rk}
\begin{figure}[h!]
		\centering		
		\begin{subfigure}[b]{0.5\textwidth}
			\centering		\hspace{2cm}\includegraphics[width=7.5cm]{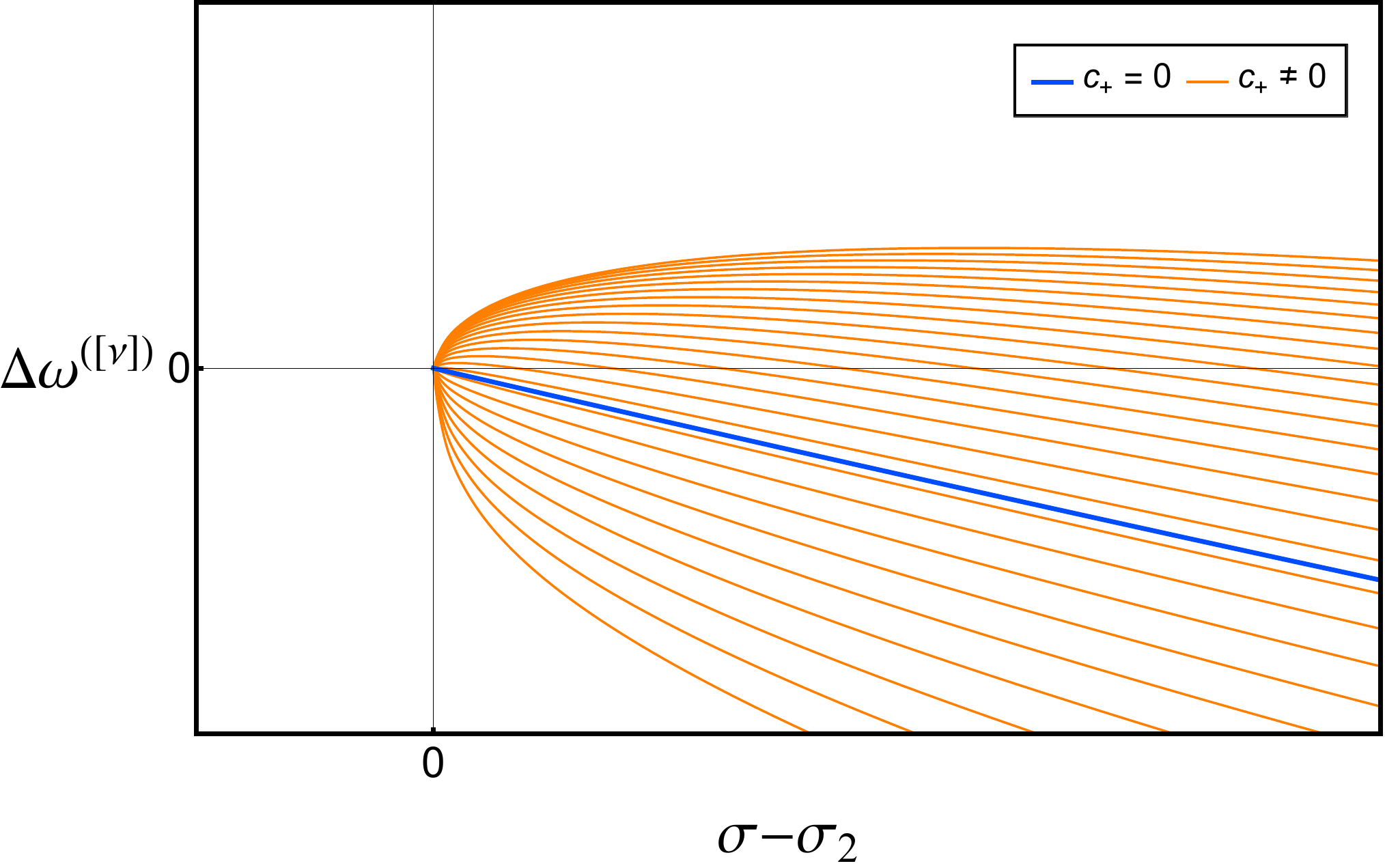}
			\caption{Search $c_+(r) = 0$}
		\end{subfigure}%
		\begin{subfigure}[b]{0.5\textwidth}
			\centering		\hspace{2cm}\includegraphics[width=7.5cm]{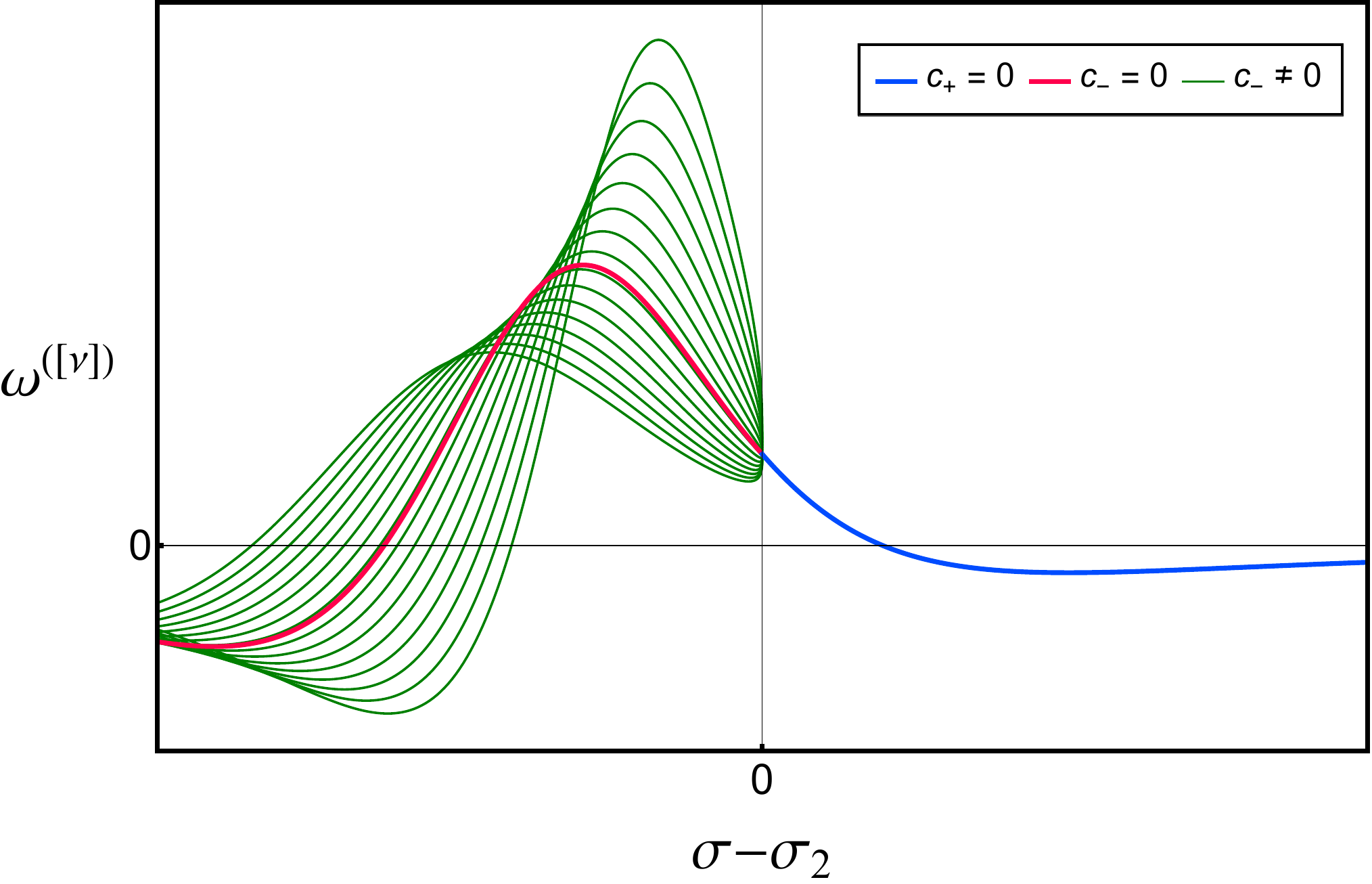}
			\caption{Search $c_-(\kappa) = 0$}
		\end{subfigure}%
		
		\begin{subfigure}[b]{0.5\textwidth}
			\centering
			\includegraphics[width=7.5cm]{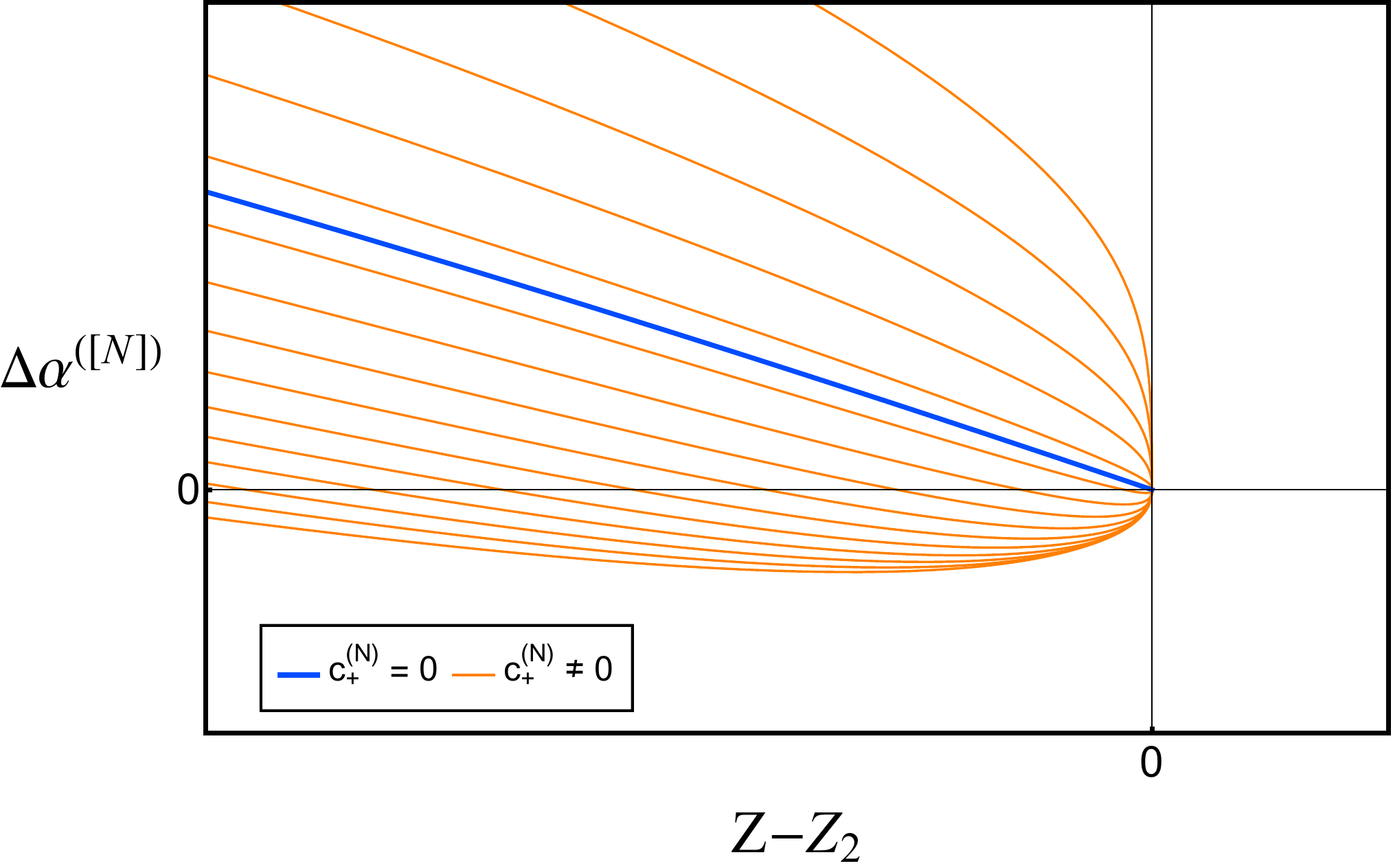}
			\caption{Search $c_+^{(\mathcal{N})}(\LL) = 0$}
		\end{subfigure}%
		\begin{subfigure}[b]{0.5\textwidth}
			\centering
			\includegraphics[width=7.5cm]{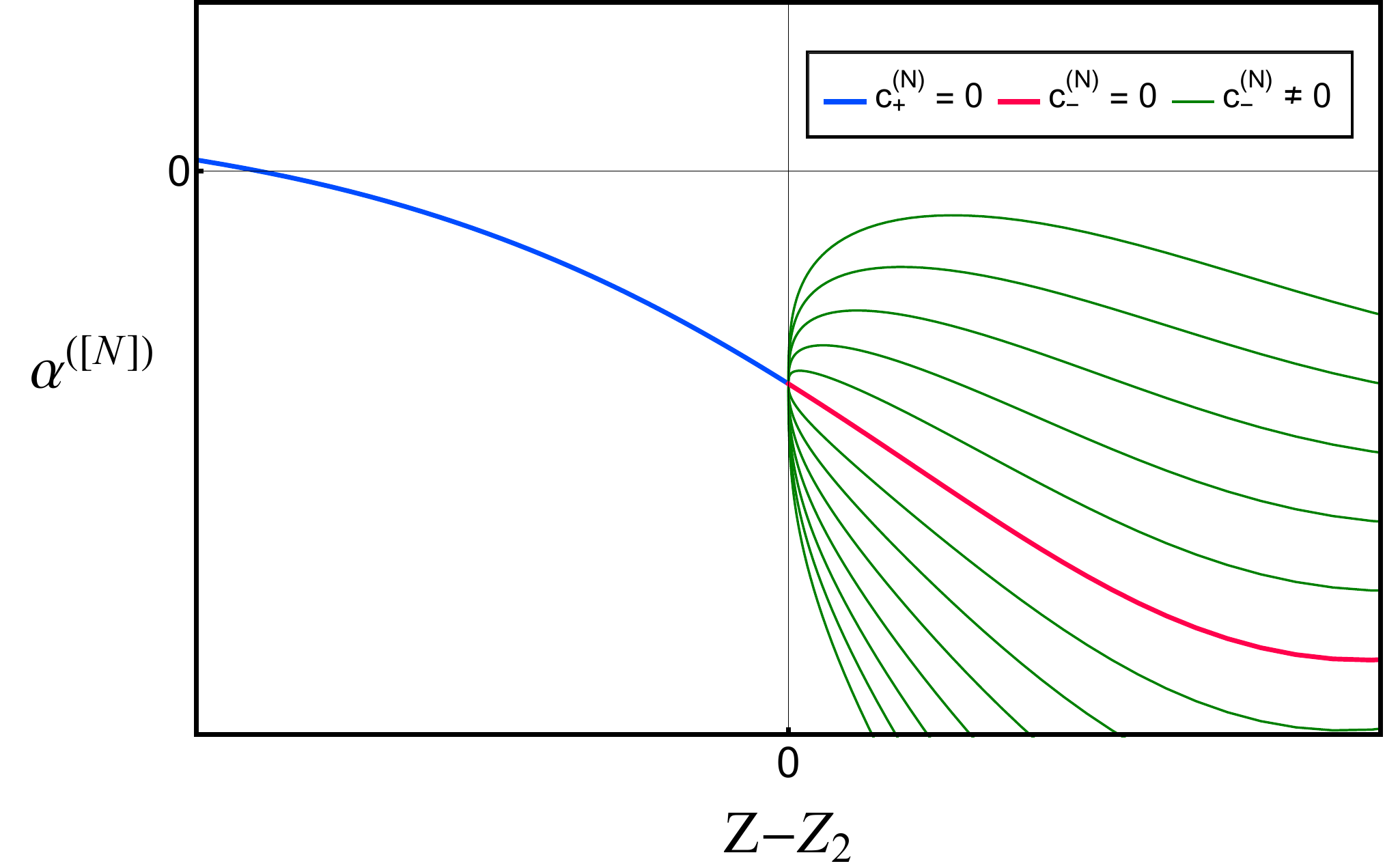}
			\caption{Search $c_-^{(\mathcal{N})}(\theta) = 0$}
		\end{subfigure}%
		\caption{\small Search of smooth solutions ($c_{\pm}=0$) and smooth linear modes ($c_{\pm}^{\mathcal{N}}=0$). (a) shows function (\ref{eq:DeltaW_nu}) for trajectories $P_2$-$P_6$ with different $r$. Orange lines are trajectories with $c_+\neq0$ while the blue one has $c_+=0$, the smooth one. (b) shows the $[\nu]$-derivative of $\w(\sigma)$ for the smooth trajectory $P_2$-$P_6$ in blue, for trajectories $P_2$-$P_4$ with $c_-\neq0$ in green and for the one with $c_-=0$ in red; namely, the blue and red trajectories have a smooth connection at $P_2$. (c) shows function (\ref{eq:DeltaW_nu}) adapted to $\alpha(Z)$; namely, $\Delta\alpha^{([\mathcal{N}])}(Z) := \alpha^{([\mathcal{N}])}(Z) - \alpha^{([\mathcal{N}])}(Z_2)$. Orange lines are associated with linear modes on the interior of the cone with $c_{+}^{(\mathcal{N})}\neq0$ while the blue one has $c_{+}^{(\mathcal{N})}=0$, the smooth one. (d) shows the $[\mathcal{N}]$-derivative of the smooth linear mode on the interior of the cone in blue, the modes on the exterior with $c_{-}^{(\mathcal{N})}\neq0$ in green and the smooth one, $c_{-}^{(\mathcal{N})}=0$ in red; namely, the blue and red curves have a smooth connection on the acoustic cone.}
		\label{fig:Search_Omega_k_1-modes_new_search_cp_cm}
	\end{figure}

	In this appendix we describe the numerical implementation of our strategy to construct SSs given in section~\ref{sec:Smooth_Solutions}. The method is the same for the search of $c_{\pm}(\kappa)=0$ and SLMs (under minor modifications) as we can see in fig.~\ref{fig:Search_Omega_k_1-modes_new_search_cp_cm}. For this reason we only describe the steps to determine $c_{+}(r)=0$.
	\begin{enumerate}
		\item Given $(d,\ell)$ we choose a low integer $1,2,3,...$ and determine the interval $I_{[\nu]}$ of $r$ where the integer part of $\nu$ given in (\ref{eq:expansion_w_in_powers_sigma_Z2}), $[\nu]$, is equal to this integer but $\nu\neq [\nu]$. For example, if we choose $[\nu]=3$ then
		\beq
			I_{3} = \{\ r\in (1,\reye)\quad /\quad  [\nu(r)] = 3,\ \nu(r)\neq[\nu(r)]\ \}.
		\eeq
		\item We construct trajectories $P_2$-$P_6$ for values of $r$ in this interval following the spectral methods described in appendix~\ref{sec:Appendix_Spectral_Methods} adapted to $\omega(\s)$. To increase the accuracy of our method we can perform a change of coordinates $\s(x)$ that concentrates a high density of points close to $P_2$.
		\item After differentiating $\omega(\s)$ $[\nu]$-times and subtracting its value at $\s_2$ we get ($\xi := \s - \s_2$, $\Delta\omega^{([\nu])}(\xi) := \omega^{([\nu])}(\xi) - \omega^{([\nu])}(0)$)
		\beq
			\Delta\omega^{([\nu])}(\xi) \sim \tilde{c}_+ \xi^{\nu-[\nu]} + a_1 \xi + a_2 \xi^2 + \tilde{c}_1 \xi^{\nu-[\nu]+1}+...
			\label{eq:DeltaW_nu}
		\eeq 
		where the coefficients $\tilde{c}_+,a_1,a_2,\tilde{c}_1,...$ are the coefficients in (\ref{eq:expansion_w_in_powers_sigma_Z2_plusminus}) redefined to include extra factors coming from differentiation. 
		\item Fitting the values of $\Delta\omega^{([\nu])}(\xi)$ for $\xi \ll 1$ using enough terms we obtain $c_+$.
		
		\item Iterating these steps for different values of $r$ we can determine if $c_+$ transits from negative to positive values. In this case we focus our search on the neighborhood of this transition to extract the value of $r$ such that $c_+(r)=0$ by interpolation.
		
		\item The precision of our results can be estimated using $a_1$. We can compare the numerical value of $a_1$ when $c_+=0$ with its analytic value (obtained from the expansion of $\omega(\xi\sim0)$ imposing $c_+=0$).
	\end{enumerate} 
	This method is very powerful to determined $c_+(r)=0$ when $[\nu]$ is a low integer, but unfortunately, its accuracy rapidly decreases when $[\nu]$ grows because it relies on $[\nu]$-numerical differentiation and the interval $I_{[\nu]}$ shrinks, see fig.~\ref{fig:Search_rk_dimensions}. For this reason our exploration of SSs was restricted to $[\nu]$ from 1 to 7. Other difficulties that we find are regions where $\nu\sim[\nu]$ or $[\nu]+1$. From a numerical point of view, in these situations it is very difficult to distinguish between $\xi^{\nu}$ and $\xi^{[\nu]}$ or $\xi^{[\nu]+1}$ respectively. In $d=2$ we find that zeros of $c_+(r)$ are dangerously located in these regions and the errors are larger than in higher dimensions. Other problematic regions in the space of parameters where we do not have accurate access are of course situations like $r\ll1$, $\ell\ll 1$, $\ell\gg1$, $d\gg 1$ and the limits of $\kappa\sim \kappa_{\min},\ \kappa_{\max}$. Finally, we want to remark that the construction of SSs is highly sensitive to small deviations from the exact parameters. First of all, once $(d,\ell)$ are fixed these solutions are isolated points in the 2-dimensional space $(r,\kappa)$; therefore, any deviation is non-smooth. However, we mean that this is sensitive in a practical sense because a small $c_{\pm}$ is translated into errors that rapidly increase with the number of derivatives. In the space of parameters, the region around the exact SSs where a numerical self-similar solution has $k$ extra continuous derivatives decreases exponentially when $k$ grows. Then, in practice even if $[\nu+1]$-numerical-derivatives are continuous, higher derivatives may not enjoy this property. Note that these limitations are also applicable to the construction of SLMs.

	\subsection{Time-Evolution Scheme}
	
\begin{figure}[h!]
		\centering	
		\begin{subfigure}[b]{0.5\textwidth}
			\centering
			\includegraphics[width=7.5cm]{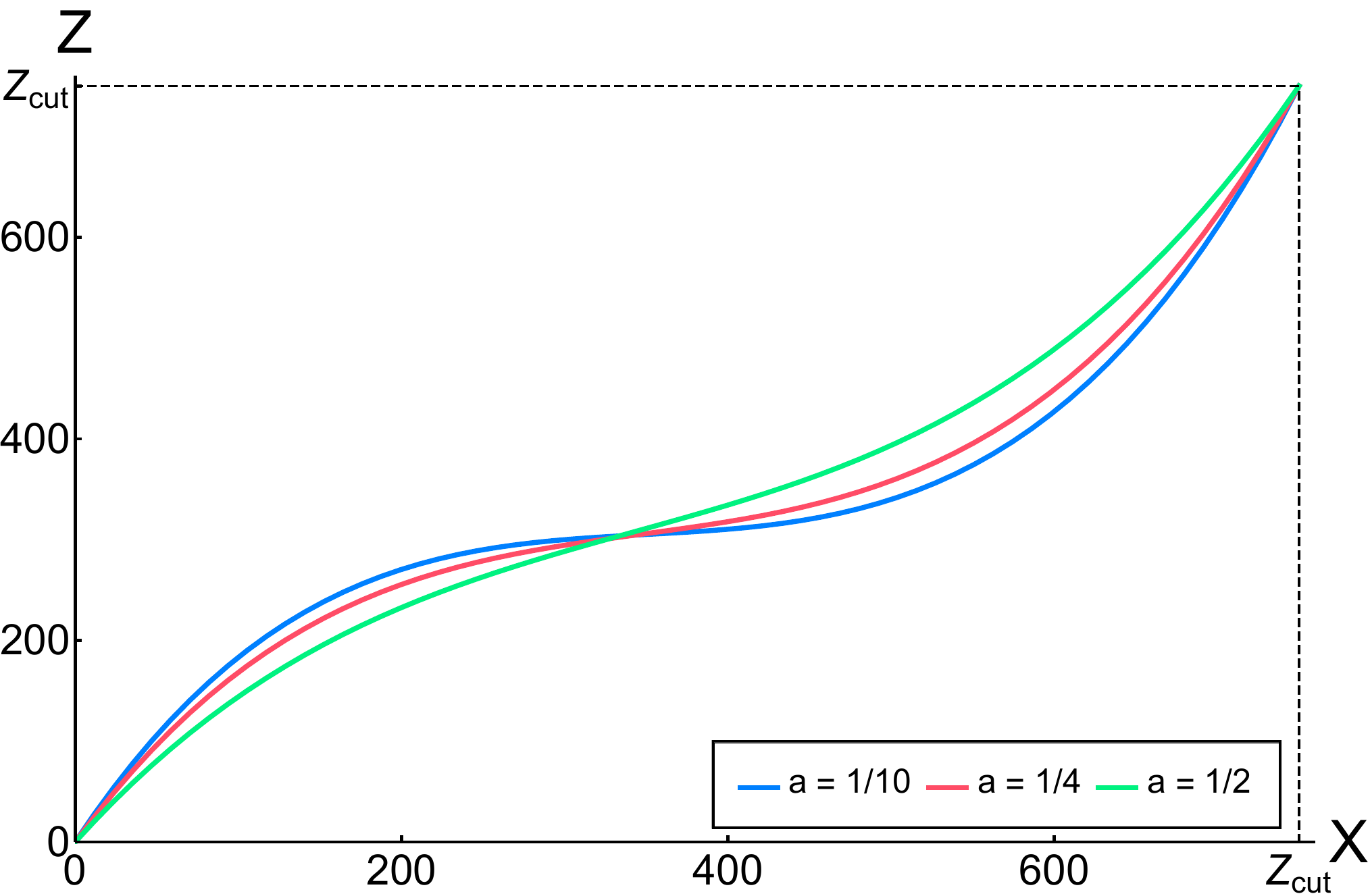}
			\label{fig:ZX_coordinates}
		\end{subfigure}%
		\begin{subfigure}[b]{0.5\textwidth}
			\centering		\hspace{2cm}\includegraphics[width=7.5cm]{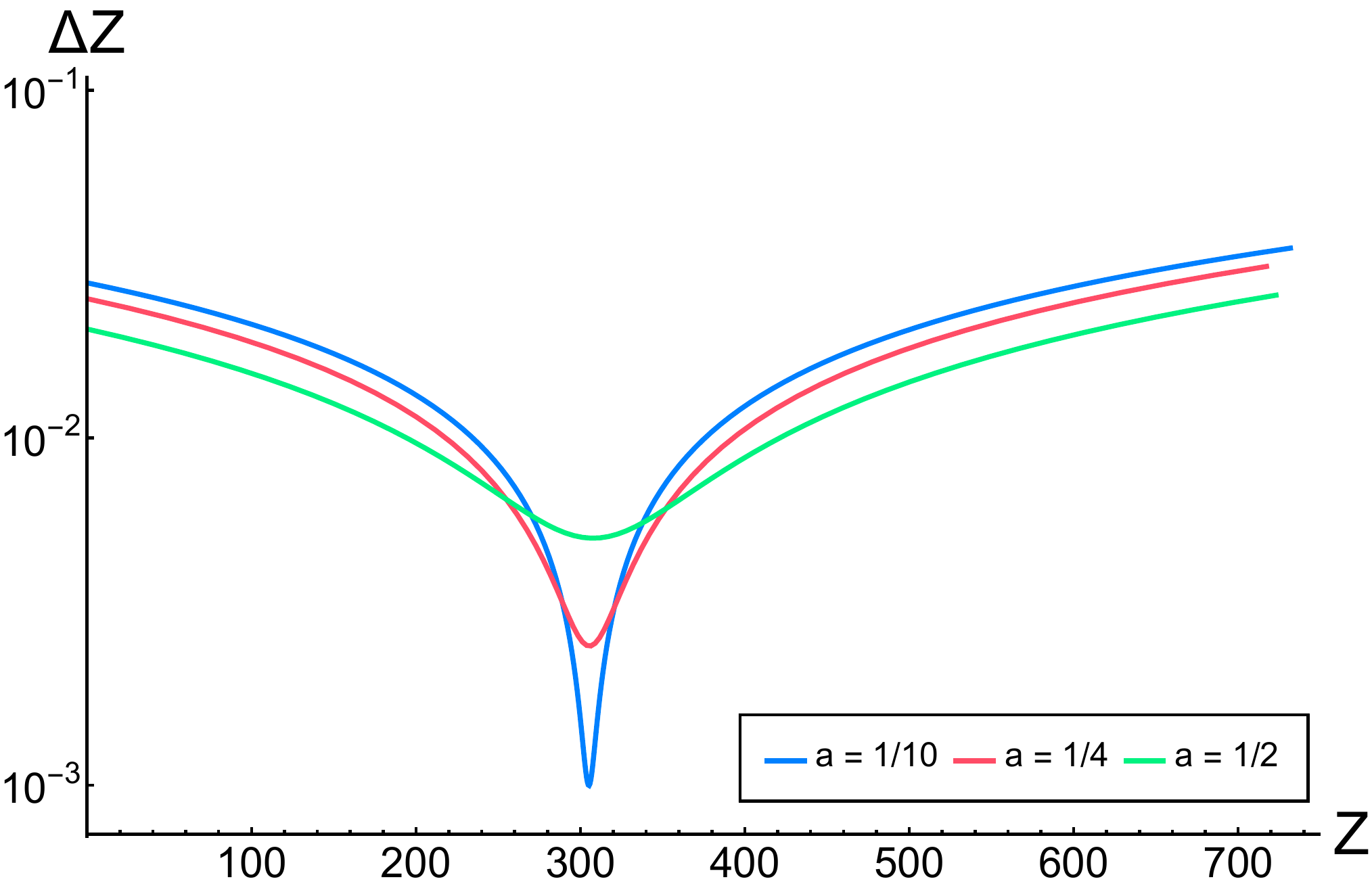}
			\label{fig:ZX_coordinates_DZ}
		\end{subfigure}%
		\caption{\small (a) three representations of the relation $Z(X)$ given in (\ref{eq:ZX_coordinates}) with the plateau located at the same $Z$. (b) representation of $\Delta Z$ when $X$ is discretized using a grid with constant steps size, $Z_{cut}=750$ and $2^{16}$ points. It results in a region with a localized high density of points $(1/\Delta Z)$.}
		\label{fig:ZX_coordinates_full}
	\end{figure}

	Time-evolution is carried out making use of a numerical scheme with the following features:
	\begin{itemize}
		\item Time-evolution: 
		\begin{itemize}
			\item Scheme: An explicit 4th order Runge-Kutta scheme is used to advance in time. The relation ``speed-accuracy" is appropriate for our purpose.
			\item Dissipation: Spurious oscillations coming from the edges of the grid are tamed with the implementation of numerical dissipation.
			\item Variable time-increments: When we deal with rapidly evolving structures, as the loss of regularity in finite time, variable time-increments \cite{SaitoSasaki} provide a better description of the process. For this problem we use the relation
			\beq
			\Delta \tau_{n+1} = \min\left(1,\frac{\max_{z}{|\partial_z\hat{\rho}(0)|}}{\max_{z}{|\partial_z\hat{\rho}(\tau_{n})|}}\right) \Delta \tau_0 \qquad \tau_{n+1} = \sum\limits_{i=0}^{n}\Delta \tau_i.
			\eeq 
			
			\item Maximum time: When instabilities are present the numerical noise imposes a natural impediment to perform large-time simulations. Although unstable modes are not contained in the initial data, the noise triggers these modes. Then, there is a {\em numerical ramp} determined by the line $\epsilon_{\text{noise}}e^{\LL_{max}\tau}$ that contaminates the simulation; see fig.~\ref{fig:numerical_ramp}.
			\begin{figure}[t!]
				\centering
				\includegraphics[width=8cm]{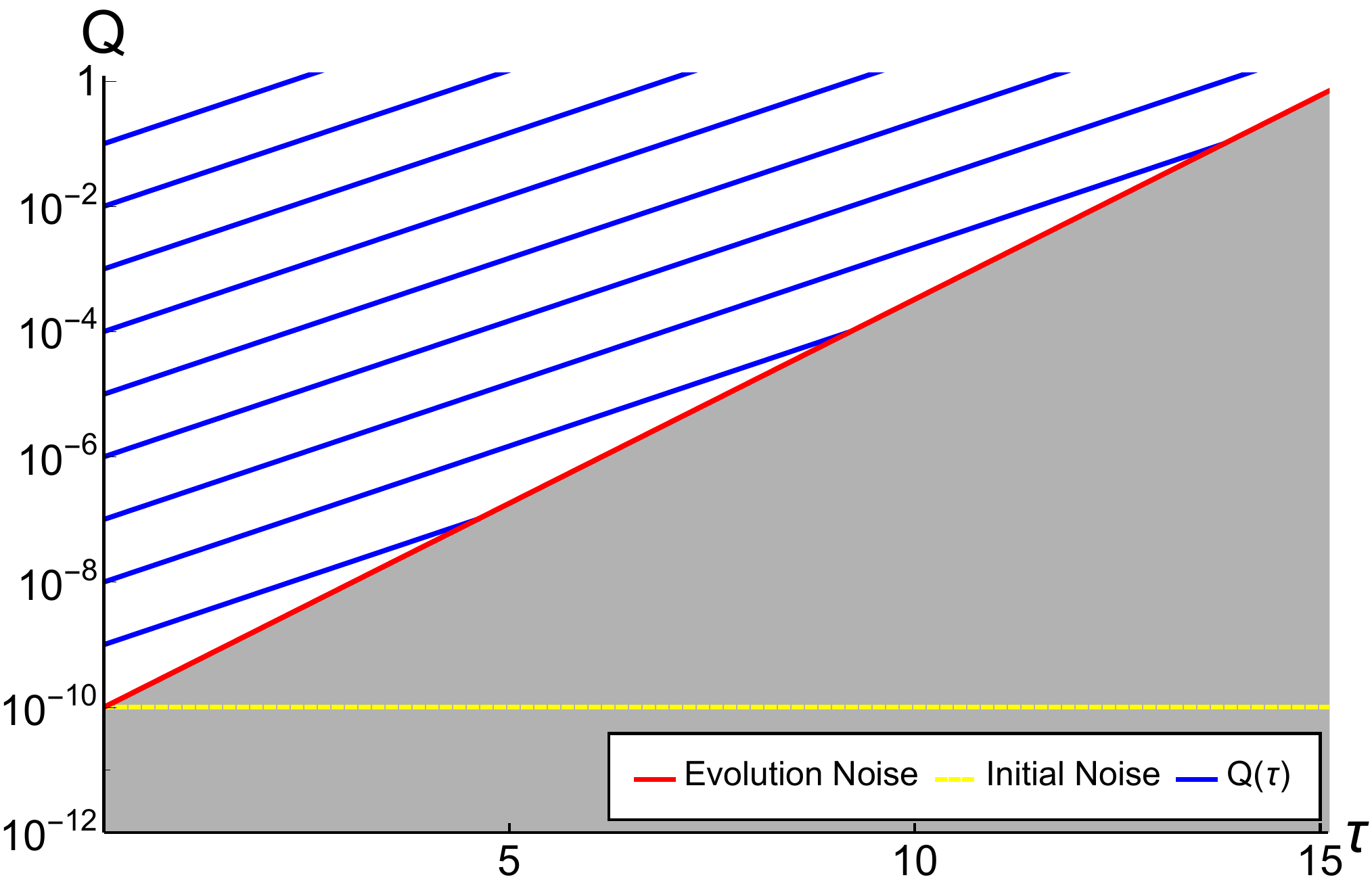}
				\caption{\small Blue lines represent a quantity of our interest that grows like $Q(\tau) \sim \epsilon e^{\LL \tau}$. The yellow dashed line represents the initial level of noise, $\epsilon_{\text{noise}}$ placed at $10^{-10}$ in this example, and the red line the {\em numerical ramp}, $\epsilon_{\text{noise}}e^{\LL_{\max}\tau}$. When $Q$ is close to the numerical ramp the contamination is too high and the simulation must be stopped.}
				\label{fig:numerical_ramp}
			\end{figure}
		\end{itemize}
		\item Radial-coordinate: 
		\begin{itemize}
			\item Truncation: The radial coordinate is truncated at $Z_{cut}\gg 1$. Despite blow-up profiles have a slow decay at infinity we do not observe that this truncation has effect on our results. This is a consequence that the tail is almost static and we do not impose boundary conditions at $Z_{cut}$. Instead of that we use a backward differentiation scheme close to $Z_{cut}$.
			\item Differentiation: A 6th-order central finite difference scheme is applied on the grid except close to $Z_{cut}$ where we switch to a backward scheme of the same order. In order to exploit the symmetries of $(\hat\rho,\hat{u})$ at the origin we extend the grid some point to the left of $Z=0$.
			\item Regularization: In (\ref{eq:Rho_dtau}), the term $\frac{d-1}{Z}$ may be regularized using the L'Hopital rule in the neighborhood of the origin.
			\item Sharp structures: depending on the needs of the problem we use two different grids. The first one is the $Z$-coordinate with the standard constant step size which works very well when profiles do not develop sharp structures. When we deal with the formation of singularities, this grid is not appropriate because localized structures are not well described or we waste computational power in regions of the grid that are well behaved. In this situation we increase the density of points in a particular location of $Z$ using a constant step size in $X$ 
			\beq
			Z(X) = C \left((1-a)\left(\frac{2}{Z_{cut}+\alpha}X-1\right)^q+\frac{2a}{Z_{cut}+\alpha} X+(1-a)\right)
			\label{eq:ZX_coordinates}
			\eeq
			with $q=1,3,5,...$, $0<a\leq 1$, $-R<\alpha<R$ and constant
			\beq
			C = \frac{Z_{cut}(Z_{cut}+\alpha)^q}{2aZ_{cut}(Z_{cut}+\alpha)^{q-1}+(1-a)(Z_{cut}-\alpha)^q+(1-a)(Z_{cut}+\alpha)^q}.
			\eeq
			This relation has the following properties: $Z(0)=0$, $Z(Z_{cut})=Z_{cut}$, $Z(X)=X$ when $a=1$, $Z(X)\sim c X$ for $X\ll1$ and $Z(X)$ is almost constant at a particular value of $Z$. It translates a constant step size in $X$ into a high density of points around a particular value of $Z$. An example of $Z(X)$ can be found in fig.~\ref{fig:ZX_coordinates_full}.
			
		\end{itemize} 
	\end{itemize}
	
\end{document}